\documentclass[11pt]{amsart}
\usepackage{amsmath,amscd,amssymb,amsthm,amsfonts}

\theoremstyle{plain}
\newtheorem{thm}[equation]{Theorem}
\newtheorem{lem}[equation]{Lemma}

\newtheorem{prop}[equation]{Proposition}

\theoremstyle{definition}
\newtheorem{defi}[equation]{Definition}
\newtheorem{rmk}[equation]{Remark}

\theoremstyle{remark}

\numberwithin{equation}{section}

\newcommand{\bib}{\bibitem}
\newcommand{\ra}{\rightarrow}
\newcommand{\lra}{\longrightarrow}
\newcommand{\nin}{\noindent}
\newcommand{\bsk}{\bigskip}

\newcommand{\ssk}{\smallskip}
\newcommand{\sst}{\scriptstyle}
\newcommand{\contrac}{\mathop{\raise.45ex\hbox{$\underline{\hskip7pt}$} \raise.5ex\hbox{$\mkern-2mu\scriptstyle|$}}\nolimits}

\newcommand{\al}{\alpha}
\newcommand{\be}{\beta}
\newcommand{\ga}{\gamma}
\newcommand{\Ga}{\Gamma}
\newcommand{\de}{\delta}
\newcommand{\De}{\Delta}
\newcommand{\la}{\lambda}
\newcommand{\La}{\Lambda}
\newcommand{\Th}{\Theta}
\newcommand{\si}{\sigma}
\newcommand{\Si}{\Sigma}
\newcommand{\om}{\omega}
\newcommand{\Om}{\Omega}
\newcommand{\ka}{\kappa}

\newcommand{\vep}{\varepsilon}
\newcommand{\vphi}{\varphi}
\newcommand{\vth}{\vartheta}
\newcommand{\na}{\nabla}

\newcommand{\CA}{\mathcal A}
\newcommand{\CB}{\mathcal B}

\newcommand{\CCD}{\mathcal D}

\newcommand{\CG}{\mathcal G}
\newcommand{\CH}{\mathcal H}

\newcommand{\CL}{\mathcal L}
\newcommand{\CM}{\mathcal M}

\newcommand{\CP}{\mathcal P}

\newcommand{\CS}{\mathcal S}

\newcommand{\CZ}{\mathcal Z}

\newcommand{\BP}{\mathbb P}
\newcommand{\BT}{\mathbb T}
\newcommand{\BR}{\mathbb R}
\newcommand{\BZ}{\mathbb Z}
\newcommand{\BC}{\mathbb C}

\newcommand{\te}{\text{e}}
\newcommand{\ti}{\text{i}}

\newcommand{\bq}{\mathbf{q}}

\newcommand{\bc}{\mathbf{c}}
\newcommand{\bb}{\mathbf{b}}
\newcommand{\bh}{\mathbf{h}}

\newcommand{\bw}{\mathbf{w}}

\newcommand{\bell}{\boldsymbol{\ell}}

\newcommand{\CLH}{\hat{\mathcal L}}
\newcommand{\pX}{\partial X}
\newcommand{\XC}{X^{cut}}
\newcommand{\pXC}{\partial X^{cut}}
\newcommand{\Det}{\text{Det}}

\begin{document}

\title{Abelian Chern-Simons theory}
\author{Mihaela Manoliu}
\address{Department of Mathematics\\
University of Texas\\
Austin, TX-78712}
\email{miha@math.utexas.edu}

\begin{abstract}
We give a construction of the abelian Chern-Simons gauge theory from the point of view of a $2+1$ dimensional topological quantum field theory. 
The definition of the quantum theory relies on geometric quantization ideas which have been previously explored in connection to the nonabelian Chern-Simons theory \cite{JW,ADW}.
We formulate the topological quantum field theory in terms of the category of extended 2- and 3-manifolds introduced by Walker \cite{Wa}
and prove that it satisfies the axioms of unitary topological quantum field theories formulated by Atiyah \cite{A1}.
\end{abstract}

\maketitle
\newpage
\tableofcontents
\newpage

\section{Introduction} \label{s:intro}

A $2+1$ dimensional topological quantum field theory (TQFT) with Lagrangian the nonabelian Chern-Simons invariant of connections was introduced by Witten several years ago \cite{Wi1}.
It was shown that this theory generates invariants of framed oriented links in arbitrary 3-manifolds. It also leads to topological invariants of closed oriented 3-manifolds endowed with a 2-framing \cite{A3}.
Witten's approach \cite{Wi1} to defining these invariants is based 
on Feynman's path integral, results from two-dimensional conformal field theory and the axioms of TQFT formulated in \cite{A1}. 
A mathematically rigorous derivation of Witten's 3-manifold invariants is given in \cite{RT,Wa}.
Moreover, Walker \cite{Wa} gives a complete construction of the $SU(2)$ Chern-Simons theory, proving that it satisfies the axioms of a TQFT.
A simpler model of a $2+1$ dimensional TQFT is the Chern-Simons theory with finite gauge group treated in full detail in \cite{FQ}.

In line with the above mentioned developments relating to the Chern-Simons gauge theory, this paper presents a construction of the {\em abelian} Chern-Simons theory as a $2+1$ dimensional TQFT.
Given the group $\BT$ of complex numbers of unit modulus and an (even) integer $k$, called the  level, we aim to associate to every closed oriented 2-dimensional manifold $\Si$ a finite dimensional vector space $\CH(\Si)$ and to every compact oriented 3-manifold $X$ a vector $Z_X$ in the vector space $\CH(\pX)$ functorially associated to the boundary $\pX$ of $X$.

The first assignment, $\Si \ra \CH(\Si)$, results from the geometric quantization of the moduli space $\CM_{\Si}$ of flat $\BT$-connections on $\Si$, having in view that $\CM_{\Si}$ is identified to the symplectic torus $H^1(\Si;\BR)/H^1(\Si;\BZ)$ with standard symplectic form $\om_{\Si}$.
There are various models for the quantization of $\CM_{\Si}$ according to  the type of polarization one chooses on this space.
The holomorphic (K\"{a}hler) quantization of $\CM_{\Si}$ is discussed in \cite{ADW,A2,G}, where it is shown that the vector spaces constructed  for various choices of complex structures on $\CM_{\Si}$ are  projectively identified.
In this paper we choose to quantize $\CM_{\Si}$ by using {\em real} polarizations. 
To each rational Lagrangian subspace $L$ in $H^1(\Si;\BR)$, i.e. a Lagrangian subspace with the property that $L \cap H^1(\Si;\BZ)$ generates $L$ as a vector space, there corresponds an invariant real polarization $\CP_L$ of the torus $\CM_{\Si} \cong H^1(\Si;\BR)/H^1(\Si;\BZ)$.
The quantization of $\CM_{\Si}$, at level $k$, in this real polarization constructs a $k^g$-dimensional inner product space $\CH(\Si,L)$, where
$g$ denotes the genus of $\Si$. 
The quantization of symplectic tori in a real polarization is discussed in detail in \cite{Ma} and the results obtained there are directly applicable to the case of the moduli space $\CM_{\Si}$. 
According to the results from \cite{Ma}, the vector spaces associated to $\Si$ and different choices of rational Lagrangian subspaces in $H^1(\Si;\BR)$ are all projectively identified, similarly to the holomorphic quantization case.
The projective factor, an 8-th root of unity, is expressible in terms of the Maslov-Kashiwara index $\tau(\cdot,\cdot,\cdot)$ of a triple of Lagrangian subspaces of $H^1(\Si;\BR)$. 
More precisely, for any two rational Lagrangian subspaces $L_1,L_2 \subset H^1(\Si;\BR)$ there is a canonically defined unitary operator $F_{L_2 L_1} : \CH(\Si,L_1) \ra \CH(\Si,L_2)$; if $L_3$ is another such Lagrangian subspace in $H^1(\Si;\BR)$, then the unitary operators relating the vector spaces corresponding to $\Si$ and each of these Lagrangian subspaces satisfy the composition law:
$F_{L_1 L_3} \circ F_{L_3 L_2} \circ F_{L_2 L_1} \, = \, \te^{- \frac{\pi \ti}{4} \tau(L_1,L_2,L_3)} \, I $.
The idea of using real polarizations in quantizing the moduli space $\CM_{\Si}$ and the method of construction of the vector space $\CH(\Si,L)$ was inspired by \cite{JW}, where, in the context of the nonabelian Chern-Simons theory, the authors address the problem of quantizing the moduli space of flat $SU(2)$ connections on a closed oriented 2-manifold in a particular real polarization of that space.

For the second assignment, $X \ra Z_X$, we start from the known fact that, for $X$ a compact oriented 3-manifold with boundary, the image of the restriction map $H^1(X;\BR) \ra H^1(\pX;\BR)$ defines a Lagrangian subspace $L_X$ in the  symplectic vector space $(H^1(\pX;\BR), \om_{\pX})$.
Then we provide a canonical construction which defines $Z_X$ as an element of the vector space $\CH(\pX, L_X)$ associated to the closed 2-manifold $\pX$ and the rational Lagrangian subspace $L_X$.
The construction of $Z_X$ relies on geometric quantization ideas.
Following the Chern-Simons line construction exposed in \cite{F}, we introduce a prequantum line bundle $\CL_{\pX}$ over the moduli space of flat connections $\CM_{\pX}$. 
The line bundle $\CL_{\pX}$ carries a natural connection with curvature the symplectic form $k \om_{\pX}$.
The rational Lagrangian subspace $L_X \subset H^1(\pX;\BR)$ defines an invariant real polarization on $\CM_{\pX}$. 
One of the leaves of this polarization is the image $\La_X$ in $\CM_{\pX}$ of the moduli space $\CM_X$ of flat $\BT$-connections on $X$ under the map $r_X : \CM_X \ra \CM_{\pX}$ determined by restricting connections on $X$ to the boundary $\pX$. 
We show that the leave $\La_X$ is a Bohr-Sommerfeld leave, that is, the restriction of the line bundle $\CL_{\pX}$ to $\La_X$ admits nonzero covariantly constant global sections.
The vector $Z_X$ is defined as a section of the line bundle $\CL_{\pX} \bigr|_{\La_X}$ times a section of the bundle of half-densities on $\La_X$. 
The former section is defined in terms of the Chern-Simons functional of flat connections on $X$ and the latter in terms of the Reidemeister torsion invariant of the 3-manifold $X$.
For $X$ a closed manifold the construction gives for $Z_X$  a complex number.
The definition of $Z_X$, that is, finding the appropriate ingredients which should enter into this definition, was inspired in part by the results for the nonabelian Chern-Simons theory from \cite{JW}, as well as by the path integral formulation of the closed 3-manifold $SU(2)$ invariant given in \cite{FG,Wi1}.

We find that a natural way to incorporate the assignments $\Si \ra \CH(\Si)$ and $X \ra Z_X$,  defined as outlined above, into a TQFT is to make use of the category of {\em extended} 2- and 3-manifolds introduced by Walker in his treatment \cite{Wa} of the nonabelian Chern-Simons theory.
In this paper an extended 2-manifold is a pair $(\Si,L)$, with $\Si$ a closed oriented 2-manifold and $L$ a rational Lagrangian subspace in $H^1(\Si;\BR)$.
An extended 3-manifold is a triple $(X,L,n)$, with $X$ a compact oriented 3-manifold, $L$ a rational Lagrangian subspace in $H^1(\pX;\BR)$ and $n \in \BZ/ 8\BZ$.
Then the theory assigns to $(\Si,L)$ the finite dimensional inner product space $\CH(\Si,L)$ and to $(X,L,n)$ the vector $Z_{(X,L,n)} = \te^{\frac{\pi \ti}{4} n} F_{L L_X}(Z_X)$ belonging to the vector space $\CH(\pX,L)$.
We prove that these assignments satisfy the axioms \cite{A1} of a unitary TQFT, that is, the functoriality, orientation, disjoint union and gluing properties.

This paper is organized as follows.
In Sect.\ref{s:moduli} we review some fairly standard material on the moduli space of flat $\BT$-connections.
We outline properties of this moduli space, essential to the constructions of the subsequent sections, for the cases when connections are on 2- and 3-dimensional manifolds.
In Sect.\ref{s:CSfunct} we recall first the definitions of induced principal bundles and induced connections. Then we introduce the Chern-Simons functional of $\BT$-connections on a 
3-manifold with and without boundary and present properties of this functional which are relevant for the sections to follow.
In Sect.\ref{s:line} we give the construction of the prequantum line bundle over the moduli space of flat $\BT$-connections on a closed 2-manifold.
Then we show that, if $X$ is a 3-manifold with boundary, the pullback of the prequantum line bundle $\CL_{\pX}$ over $\CM_{\pX}$ to $\CM_X$ has a covariantly constant section defined in terms of the Chern-Simons functional.
The constructions of this section are based almost entirely on the material in (\cite{F},\S 2).
Sect.\ref{s:qth} defines the quantum theory. We construct the finite dimensional Hilbert space $\CH(\Si,L)$ associated to a closed oriented 2-manifold $\Si$ and a rational Lagrangian subspace $L \subset H^1(\Si;\BR)$.
Then, for a compact oriented 3-manifold $X$, we construct the vector $Z_X$ belonging the Hilbert space $\CH(\pX,L_X)$.
Sect.\ref{s:TQFT} contains the definition of the abelian Chern-Simons TQFT. Following \cite{Wa} we introduce the notions of extended 2- and 3-manifolds, extended morphisms and gluing of extended 3-manifolds.
Using the results of Sect. \ref{s:qth}, we define then a $2+1$ dimensional TQFT based on this category of extended manifolds and prove that it satisfies the required axioms.
Sect.\ref{s:pathint} relates the definition of the vector $Z_X$ associated to the 3-manifold $X$, given in Sect.\ref{s:qth}, to results obtained from the path integral approach to the Chern-Simons gauge theory. 

Throughout this paper manifolds, bundles, sections, maps are assumed smooth.

\bsk


\section{The moduli space of flat $\BT$-connections} \label{s:moduli}

\subsection{The space of $\BT$-connections}

Let $M$ be a smooth manifold.
The space $\CG_M = \text{Map}(M,\BT)$ forms a group under pointwise multiplication. $\CG_M$ is the  {\em group of gauge transformations} on $M$.
The space of components $\pi_0(\CG_M)$ is isomorphic to $H^1(M;\BZ)$ \cite{AB}.
For any principal $\BT$-bundle $\pi : P \ra M$, the group $\CG_M$ can be identified with the group $\CG_P = \text{Aut}(P)$ of bundle automorphisms of $P$, that is, the group of $\BT$-equivariant maps $\phi : P \ra P$ covering the identity on $M$.
An element $u$ of $\CG_M$ defines a bundle map $\phi_u : P \ra P$ by $\phi(p) = p \cdot u(\pi(p))$ and, 
conversely, an automorphism $\phi : P \ra P$ defines a map $u_{\phi} : M \ra \BT$. 
We call $\CG_P$ the group of gauge transformations of $P$.

Let $\CA_M$ denote the space of $\BT$-connections on $M$.
An element $\Th$ in $\CA_M$ is a connection on a principal $\BT$-bundle $P \overset{\pi}{\ra} M$.
Thus $\CA_M$ is equal to the union 
\begin{equation*}
\CA_M = \underset{P}{\bigsqcup} \, \CA_P
\end{equation*}
over all principal $\BT$-bundles $P$ on $M$.
For each bundle $P$ the space $\CA_P$ of connections on $P$ is an affine space with vector space $2 \pi \ti \, \Om^1(M;\BR)$, where  $\Om^1(M;\BR)$ is the space of 1-forms on $M$ (the Lie algebra of $\BT$ is identified with $2 \pi \ti \, \BR$).

A bundle isomorphism between two principal $\BT$-bundles $P$ and $P'$ over $M$ is a $\BT$-equivariant map $\phi : P' \ra P$ which covers the identity on $M$.
Two elements $\Th$ and $\Th'$ in $\CA_M$ are called {\em gauge equivalent\/} if there exists an isomorphism $\phi : P' \ra P$ such that $\Th' = \phi^* \Th$, that is, the connection $\Th'$ on $P'$ is the pullback under $\phi$ of the connection $\Th$ on $P$.
This defines an equivalence relation $\Th' \sim \Th$ on $\CA_M$.
We let $\CA_M/\negmedspace \sim$ denote the space  of gauge equivalence classes of $\BT$-connections on $M$.

For each $\BT$-bundle $P$ there is a right-action of $\CG_M$ on the space of connections $\CA_{P}$.
If $u : M \ra \BT$ is an element of $\CG_M$ with associated bundle automorphism $\phi_u : P \ra P$, the action of $u$ on $\Th \in \CA_P$ is described by
\begin{equation} \label{e:ggaction}
\Th \cdot u \, = \, \phi_u^* \Th \, = \, \Th + (u \cdot \pi)^* \vartheta \, ,
\end{equation}
where $\vartheta$ is the Maurer-Cartan form of $\BT$, that is, $\frac{\vth}{2 \pi \ti}$ generates $H^1(\BT;\BZ) \subset H^1(\BT;\BR)$.

A principal $\BT$-bundle $P \ra M$ has flat connections if and only if $c_1(P) \in \text{Tors} \, H^2(M;\BZ)$, that is, its first Chern class is a torsion class.
The subspace $\CA^f_M \subset \CA_M$ of flat $\BT$-connections on $M$ is given by the union
\begin{equation*}
\CA_M^f = \underset{\substack{P \\ c_1(P) \in \text{Tors} \, H^2(M;\BZ)}}{\bigsqcup} \CA_P^f \, ,
\end{equation*}
where $\CA^f_P$ is the space of connections $\Th$ on $P$ with curvature $F_{\Th} = d \Th =0$.
For each $\BT$-bundle $P$ with first Chern class torsion, the {\em moduli space of flat connections\/} on $P$ is the quotient space $\CM_P = \CA_P^f /\CG_P$.

If $P$ and $P'$ are $\BT$-bundles over $M$ with $c_1(P) = c_1(P') \in \text{Tors} \, H^2(M;\BZ)$, then there is a canonical isomorphism $\CM_P \cong \CM_{P'}$.
To see this choose a bundle isomorphism $\phi : P' \ra P$. 
The bundle map $\phi$ determines an isomorphism $\phi^* : \CA_P \ra \CA_{P'}$ given by the pullback of a connection $\Th$ in $\CA_{P}$ to a connection $\phi^* \Th$ in $\CA_{P'}$.
The map $\phi^*$ pushes down to an isomorphism between the quotients by the groups of gauge transformations $\phi^* : \CA_{P}/\CG_P \ra \CA_{P'}/\CG_{P'}$.
The induced isomorphism $\CA_{P}/\CG_P \cong \CA_{P'}/\CG_{P'}$ does not depend on the choice of bundle isomorphism $\phi : P' \ra P$.

For each torsion class $p \in \text{Tors} \, H^2(M;\BZ)$ let 
\begin{equation*}
\CA_{M,p}^f = \underset{\substack{P \\ c_1(P)=p}}{\bigsqcup} \CA_P^f
\end{equation*}
The space $\CM_{M,p} = \CA_{M,p}^f / \negmedspace \sim$ is the moduli space of flat connections on principal $\BT$-bundles over $M$ with first Chern class equal to $p$.
For any principal $\BT$-bundle $P \ra M$ with $c_1(P)=p$, there is a natural isomorphism $\CM_{M,p} \cong \CM_P$.
The moduli space $\CM_M = \CA_M^f /\negmedspace \sim$ of gauge equivalence classes of flat $\BT$-connections on $M$ is equal to the disjoint union 
\begin{equation*}
\CM_M \, = \, \underset{p \in \text{Tors} H^2(M;\BZ)}{\bigsqcup} \CM_{M,p}
\end{equation*}
and we have
\begin{prop} \label{p:modtor}
Let $M$ be a smooth manifold. \\
(i) There is a natural identification
\begin{equation*}
\CM_M \, = \, H^1(M;\BT)
\end{equation*}
(ii) $\pi_0(\CM_M) \cong \mathrm{Tors} \, H^2(M;\BZ)$ and each connected component of $\CM_M$ is diffeomorphic to the torus $H^1(M;\BR) / H^1(M;\BZ)$.
\end{prop}

\begin{proof}
The first assertion is a consequence of a standard result in the theory of connections (\cite{KN}, ch.II) which provides the natural identification
\begin{equation*}
\CM_M = \underset{\al}{\prod} \text{Hom}(\pi_1(M_{\al},*), \BT) \, ,
\end{equation*}
where the index $\al$ labels the connected components $M_{\al}$ of $M$.
The relation between $\pi_1$ and the first homology group then gives (i).
For (ii), consider the exact sequence of groups
\begin{equation*}
0 \lra \BZ \lra \BR \xrightarrow[\exp 2 \pi \ti (\cdot)]{} \BT \lra 1
\end{equation*}
and the induced exact cohomology sequence
\begin{equation} \label{e:homexs}
\begin{split}
0 \ra H^1(M;\BZ) &\ra H^1(M;\BR) \ra \\
&\ra  H^1(M;\BT) \overset{\de}{\ra} H^2(M;\BZ) \overset{i}{\ra} H^2(M;\BR) \ra \cdots 
\end{split}
\end{equation}
All the maps in the above sequence are group homomorphisms and $\text{Im} \,\de = \text{Ker} \, i = \{ p \mid p \in \text{Tors} \, H^2(M;\BZ)\}$.
Thus, for each torsion class $p$ in $H^2(M;\BZ)$, we have $\CM_{M,p} = \de^{-1}(p) \cong \de^{-1}(0) \cong H^1(M;\BR)/H^1(M;\BZ)$.
\end{proof}

\bsk


\subsection{The moduli space of flat $\BT$-connections on a 2-manifold}

Let $\Si$ be a closed oriented 2-dimensional manifold.
We note that a principal $\BT$-bundle $Q \ra \Si$ has flat connections if and only if $c_1(Q) =0$, that is, the bundle $Q$ is trivializable.
For $\Si$ the exact sequence (\ref{e:homexs}) splits to
\begin{equation*}
0 \lra H^1(\Si;\BZ) \lra H^1(\Si;\BR) \lra H^1(\Si;\BT) \lra 0
\end{equation*}
Thus the moduli space of flat $\BT$-connections on $\Si$ is the torus
\begin{equation*}
\CM_{\Si} \, = \, H^1(\Si;\BT) \, \cong \, H^1(\Si,\BR)/ H^1(\Si,\BZ) \, .
\end{equation*}
$\CM_{\Si}$ carries a natural symplectic structure. 
The space $H^1(\Si;\BR)$ is a symplectic vector space with symplectic form $\om_{\Si}$ defined by cup product followed by evaluation on the fundamental cycle, or in de Rham cohomology by
\begin{equation*}
\om_{\Si}([\al],[\be]) \, = \, \int\limits_{\Si} \al \wedge \be \; , \qquad [\al],[\be] \in H^1(\Si;\BR) \, .
\end{equation*}
The integer lattice $H^1(\Si;\BZ) \subset H^1(\Si;\BR)$ is self-dual with respect to $\om_{\Si}$.
The symplectic form $\om_{\Si}$ on $H^1(\Si;\BR)$ descends to a symplectic form $\om_{\Si}$ on the quotient torus $\CM_{\Si} \cong H^1(\Si;\BR)/ H^1(\Si;\BZ)$.
We note that since the tangent space $T \CM_{\Si}$ is identified with $2 \pi \ti \, H^1(\Si;\BR)$ we have
\begin{equation} \label{e:sympf}
\om_{\Si}([\dot{\eta}],[\dot{\eta}']) \, = \, - \frac{1}{4 \pi^2} \, \int\limits_{\Si} \dot{\eta} \wedge \dot{\eta}' \, ,
\end{equation}
for any $[\dot{\eta}],[\dot{\eta}'] \in T \CM_{\Si}$. 
Moreover, $\om_{\Si}$ is normalized so that it gives $\CM_{\Si}$ total volume equal to 1, that is,
\begin{equation*}
\int\limits_{\CM_{\Si}} \, \frac{(\om_{\Si})^g}{g !} \, = \, 1 \, ,
\end{equation*}
where $g = \frac{1}{2} \dim H^1(\Si;\BR)$.
If $- \Si$ denotes the manifold $\Si$ with the opposite orientation then by $H^1(-\Si;\BR)$ we understand the symplectic vector space $H^1(\Si;\BR)$ with symplectic form $- \om_{\Si}$.
 
\bsk


\subsection{The moduli space of flat $\BT$-connections on a 3-manifold}

Let $X$ be a compact oriented 3-manifold with boundary $\partial X \neq \emptyset$.
The restriction of a connection over $X$ to the boundary $\partial X$ induces a map $r_X : \CM_X \ra \CM_{\partial X}$ from the moduli space of flat $\BT$-connections on $X$ into the moduli space of flat $\BT$-connections on $\partial X$.
In view of the identifications $\CM_X = H^1(X;\BT)$ and $\CM_{\partial X} = H^1(\partial X;\BT)$, the moduli spaces $\CM_X$ and $\CM_{\partial X}$ are compact abelian Lie groups and the restriction map $r_X$ is a group homomorphism.
We have

\begin{prop} \label{p:resmap}
(i) The image $\La_X = \mathrm{Im} \, \{r_X : \CM_X \ra \CM_{\partial X} \}$ is a Lagrangian submanifold of the symplectic torus $(\CM_{\partial X}, \om_{\partial X})$. \\
(ii) $r_X : \CM_X \ra \La_X$ is a principal fibre bundle with structure group a compact abelian Lie group $H$ which fits into the exact sequence 
\begin{equation*}
1 \lra \BT^q \lra H \lra \mathrm{Tors} \, H^2(X;\BZ) \lra 0
\end{equation*}
with $q = \dim H^1(X,\pX;\BR) - \dim H^0(\pX;\BR) + \dim H^0(X;\BR) -\dim H^0(X,\pX;\BR)$.
\end{prop}

\begin{proof}
There is a long exact cohomology sequence (\cite{Br},ch.V) associated to the pair of spaces 
$\partial X \subset X$
\begin{equation} \label{e:homrel}
\cdots \ra H^i(X,\partial X;G) \ra H^i(X;G) \ra H^i(\partial X;G) \ra H^{i+1}(X,\partial X;G) \ra \cdots \, ,
\end{equation}
where $G = \BZ, \BR$ or $\BT$.
At the same time the exact sequence of abelian groups $0 \ra \BZ \ra \BR \ra \BT \ra 1$ induces the long exact cohomology sequence (\cite{Br}, ch.V)
\begin{equation} \label{e:homg}
\cdots \lra H^i(A;\BZ) \lra H^i(A;\BR) \lra H^i(A;\BT) \lra H^{i+1}(A;\BZ) \lra \cdots \, ,
\end{equation}
where $A$ stands for either the space $\partial X$ or $X$, or $(X,\partial X)$ for the relative cohomology.
The above exact sequences fit into a commutative diagram:
\begin{equation*}
\begin{matrix}
   &   0   &  &   0    &  &   0   &  \\
   & \downarrow & & \downarrow & & \downarrow & \\ 
\cdots \lra & H^1(X,\pX;\BZ) & \lra & H^1(X;\BZ) & \lra & H^1(\pX;\BZ) & \lra \cdots \\
   & \downarrow & & \downarrow & & \downarrow & \\ 
\cdots \lra & H^1(X,\pX;\BR) & \lra & H^1(X;\BR) & \overset{\dot{r}_X}{\lra} & H^1(\pX;\BR) & \lra \cdots \\
   & \downarrow & & \downarrow & & \downarrow & \\
\cdots \lra & H^1(X,\pX;\BT) & \lra & H^1(X;\BT) & \overset{r_X}{\lra} & H^1(\pX;\BT) & \lra \cdots \\
   & \downarrow & & \downarrow & & \downarrow & \\
\cdots \lra & H^2(X,\pX;\BZ) & \lra & H^2(X;\BZ) & \lra & H^2(\pX;\BZ) & \lra \cdots \\
 & \downarrow & & \downarrow & & \downarrow & \\
   & \vdots &  & \vdots & & \vdots & 
\end{matrix}
\end{equation*}
where all the maps are group homomorphisms.
The differentiable and group structure of $\BT$ make
\begin{equation} \label{e:exfibr}
\begin{split}
1 \ra H^0(X,\pX;\BT) &\ra H^0(X;\BT) \ra H^0(\partial X;\BT) \ra \\
& \ra H^1(X,\partial X;\BT) \ra H^1(X;\BT) \overset{r_X}{\lra} \underset{\substack{\cap \\ H^1(\partial X;\BT)}}{\La_X} \ra 1
\end{split}
\end{equation}
an exact sequence of homomorphisms of compact abelian Lie groups.
$\La_X$ is a compact subgroup and, therefore, a submanifold of $\CM_{\partial X} = H^1(\partial X;\BT)$.
The relevant part of the sequence (\ref{e:homg}) for $A= (X,\partial X)$ is 
\begin{equation*}
\cdots \ra H^1(X,\partial X;\BT) \ra H^2(X,\partial X;\BZ) \ra H^2(X,\partial X;\BR) \ra \cdots 
\end{equation*}
which implies that $\pi_0(H^1(X,\partial X;\BT)) \cong \text{Tors} \, H^2(X,\partial X;\BZ)$.
The universal coefficient theorem (U.C.T.) (\cite{Br},ch.V) together with Poincar\'{e} duality (P.D.) give the isomorphisms
\begin{equation} \label{e:Tors}
\text{Tors} \, H^2(X,\partial X;\BZ) \, \underset{(\text{U.C.T.})}{\cong}\, \text{Tors} \, H_1(X,\partial X;\BZ) \, \underset{(\text{P.D.})}{\cong} \, \text{Tors} \, H^2(X;\BZ)
\end{equation}
Thus (ii) follows from (\ref{e:exfibr}) and (\ref{e:Tors}).
To prove (i), first recall that the tangent spaces to $\CM_X$ and $\CM_{\pX}$ are identified with the corresponding real cohomology groups by $T \CM_X = 2 \pi \ti \, H^1(X;\BR)$ and $T \CM_{\pX} = 2 \pi \ti H^1(\pX;\BR)$.
The differential $\dot{r}_X$ of $r_X$ maps $T \CM_X$ onto $T \La_X \subset T \CM_{\pX}$ and we have
\begin{align*}
\om_{\pX}(\dot{r}_X [\dot{\Th}],\dot{r}_X [\dot{\Th}']) \, &= \, - \frac{1}{4 \pi^2} \, \int\limits_{\pX} \dot{r}_X(\dot{\Th}) \wedge \dot{r}_X(\dot{\Th}') \\
&= \, - \frac{1}{4 \pi^2} \, \int\limits_X \, d\,(\dot{\Th} \wedge \dot{\Th}') \, = \, 0 \, ,\quad \text{for any } \, [\dot{\Th}],[\dot{\Th}'] \in T \CM_X \, .
\end{align*}
Since $\om_{\pX} \bigl|_{T \La_X} \equiv 0$, $\La_X$ is an isotropic submanifold of $\CM_{\pX}$.
The exact sequence (\ref{e:homrel}) for $G= \BR$ and the Poincar\'{e} duality isomorphism $H^{3-i}(X;\BR) \cong H_i(X,\pX;\BR)$ give the commutative diagram (\cite{Br},ch.VI):
\begin{equation*}
\begin{CD}
H^1(X;\BR) @>{\dot{r}_X}>> H^1(\pX;\BR) @>>> H^2(X,\pX;\BR) \\
@VV{\cong}V @VV{\cong}V  @VV{\cong}V \\
H_2(X,\pX;\BR) @>>> H_1(\pX;\BR) @>j>> H_1(X;\BR) \\
\end{CD}
\end{equation*}
where $j$ and $\dot{r}_X$ are adjoint maps.
Therefore $\dim \text{Im} \, \dot{r}_X = \dim H^1(\pX;\BR) - \dim \text{Ker} \, j$.
Combining this with the exact sequence condition $\text{Im} \, \dot{r}_X = \text{Ker} \, j$, we find that $\dim \text{Im} \, \dot{r}_X = \frac{1}{2}\dim H^1(\pX;\BR)$.
Therefore
\begin{equation*}
L_X \, = \, \text{Im} \,\{ \dot{r}_X : H^1(X;\BR) \lra H^1(\pX;\BR) \}
\end{equation*}
is a Lagrangian subspace of $H^1(\pX;\BR)$ and $\La_X$ a Lagrangian submanifold of $(\CM_{\pX}, \om_{\pX})$.
\end{proof}

\bsk


\section{The Chern-Simons functional} \label{s:CSfunct}

\subsection{Induced principal bundles and induced connections}

Consider the group $SU(2)$ of complex $2 \times 2$ unitary matrices of determinant equal to 1.
A maximal torus in $SU(2)$ is isomorphic to the group $\BT$ of complex numbers of unit modulus.
We let 
\begin{align*}
\rho : \BT \, &\hookrightarrow \, SU(2) \\
\te^{2 \pi \ti \vphi} \, &\mapsto \, \begin{pmatrix} \te^{2 \pi \ti \vphi} & 0 \\ 0 & \te^{-2 \pi \ti \vphi} \end{pmatrix} \notag
\end{align*}
be the inclusion homomorphism from the circle group $\BT$ into $SU(2)$.
The corresponding  inclusion of the Lie algebra $\text{Lie} \,\BT = 2 \pi \ti \, \BR$ into the Lie algebra $\mathfrak{su}(2) $ of traceless skew-hermitian $2 \times 2$ complex matrices is 
\begin{align*}
\rho_{*} : \text{Lie} \, \BT \, &\hookrightarrow \, \mathfrak{su}(2) \\
\al \, &\mapsto \, \begin{pmatrix} \al & 0 \\ 0 & -\al \end{pmatrix} \notag
\end{align*}

Let $Ad$ denote the adjoint action of $SU(2)$ on its Lie algebra and let $\hat{\vartheta}$ denote the Maurer-Cartan form of $SU(2)$.
On the Lie algebra $\mathfrak{su}(2)$ we choose the  symmetric $Ad$-invariant bilinear form
\begin{alignat*}{2}
\langle \cdot , \cdot \rangle^{\flat} \, : \, \mathfrak{su}(2) &\times \mathfrak{su}(2) & & \lra \BR \\
(a &\, , b) & & \longmapsto \langle a , b \rangle^{\flat} \, = \, \frac{1}{8 \pi^2} \, \text{Tr}(a b) \notag
\end{alignat*}
The bilinear form $\langle \cdot , \cdot \rangle^{\flat}$ is normalized so that the closed $3$-form $\frac{1}{6} \langle \hat{\vartheta} \wedge [ \hat{\vartheta} , \hat{\vartheta} ] \rangle^{\flat}$ represents an integral cohomology class in $H^3(SU(2); \BR)$. 
The commutator $[ \cdot , \cdot ]$ of two $\mathfrak{su}(2)$-valued forms is defined by taking the wedge product of forms and the natural Lie bracket in the $\mathfrak{su}(2)$ Lie algebra.
The bilinear form $\langle \cdot , \cdot \rangle^{\flat} $ on $\mathfrak{su}(2)$ restricts to a symmetric bilinear form $\langle \cdot , \cdot \rangle$ on $\text{Lie} \,\BT$:
\begin{alignat*}{2}
\langle \cdot , \cdot \rangle \, : \, \text{Lie} \BT &\times \text{Lie} \BT & &\lra \BR \\
(\al &, \be) &  & \longmapsto \langle \al , \be \rangle \, = \, \langle \rho_{*} \al , \rho_{*} \be \rangle^{\flat} \, = \, \frac{1}{4 \pi^2} \, \al \be\notag 
\end{alignat*}
As in the previous section we let $\vartheta$ denote the Maurer-Cartan form on $\BT$ and note that $\rho^*(\hat{\vartheta}) = \rho_{*} \vth$.

Let $M$ be a smooth manifold and consider a principal $\BT$-bundle $\pi : P \ra M$. We can use the inclusion homomorphism $\rho : \BT \hookrightarrow SU(2)$ to extend the $\BT$-bundle $P$ on $M$ to a $SU(2)$-bundle $\hat{P}$ on $M$.
The  induced $SU(2)$-bundle $\hat{P}$ is defined as the quotient of $P \times SU(2)$ by the right $\BT$-action given by $(p,a) \cdot \la = (p \cdot \la, \rho(\la)^{-1} a)$, for any $\la \in \BT$ and $(p,a) \in P \times SU(2)$.
Let $[p,a] = (p,a) \cdot \BT$ denote the $\BT$-orbit through the point $(p,a)$.
The natural right $SU(2)$-action on $P \times SU(2)$, $\, (p,a) \cdot a'= (p,a a')$, commutes with the $\BT$-action and passes therefore to the quotient $\hat{P} = P \times_{\BT} SU(2) = (P \times SU(2))/\BT$.
Hence, $\hat{\pi} : \hat{P} \ra M$ with projection $\hat{\pi}([p,a])= \pi(p)$ is a principal $SU(2)$-bundle.
We have a natural morphism of principal bundles
\begin{align*}
\rho_P : P &\lra \hat{P} \\
p &\longmapsto [p,e] \notag
\end{align*}
covering the identity map on $M$.

For any connection $\Th$ on the $\BT$-bundle $P$ there is an induced connection $\hat{\Th}$ on the $SU(2)$-bundle $\hat{P}$.
The latter is determined by the $\mathfrak{su}(2)$-valued 1-form on $P \times SU(2)$:
\begin{equation} \label{e:indconn}
\hat{\Th}_{(p,a)} \, = \, Ad_{a^{-1}} (\rho_{*}\, pr_1^* \, \Th_p) \, + \ pr_2^* \, \hat{\vartheta}_a \, ,
\end{equation}
where $pr_1 : P \times SU(2) \ra P$ and $pr_2 : P \times SU(2) \ra SU(2)$ are the natural projections.
One can easily prove that $\hat{\Th}$ is invariant under the $\BT$-action and vanishes along the $\BT$-orbits in $P \times SU(2)$.
Thus $\hat{\Th}$ pushes down to a $\mathfrak{su}(2)$-valued 1-form on $\hat{P}= P \times_{\BT} SU(2)$ and defines a connection on $\hat{P}$.
The pullback of $\hat{\Th}$ under the bundle map $\rho_P : P \ra \hat{P}$ is 
\begin{equation*}
\rho_P^* \hat{\Th} \, = \, \rho_{*} \Th \, .
\end{equation*}
Similarly, the curvature forms $F_{\Th} = d \Th$ and $F_{\hat{\Th}} = d \hat{\Th} + \frac{1}{2} [ \hat{\Th}, \hat{\Th}]$ are related by
\begin{equation*}
\rho^*_P F_{\hat{\Th}} \, = \, \rho_{*} F_{\Th} \, .
\end{equation*}

\bsk


\subsection{The Chern-Simons functional for a closed 3-manifold}

Let $X$ be a closed oriented 3-manifold.
As previously shown, for any principal $\BT$-bundle $P \ra X$ there is an induced $SU(2)$-bundle $\hat{P} = P \times_{\BT} SU(2) \ra X$ and for any $\BT$-connection $\Th$ on $P$ an induced $SU(2)$-connection $\hat{\Th}$ on $\hat{P}$.
Now, any principal $SU(2)$-bundle over a manifold of dimension $\leq 3$ is trivializable. Thus $\hat{P}$ admits a global section and we have:

\begin{defi} \label{d:CSnob}
The Chern-Simons functional of a $\BT$-connection $\Th$ on $P \ra X$ is defined by 
\begin{equation*}
S_{X,P}(\Th) \, = \, \int\limits_X \hat{s}^* \al(\hat{\Th}) \; \pmod{1}\, ,
\end{equation*}
\end{defi}
\nin where $\al(\hat{\Th}) \in \Om^3(\hat{P}; \BR)$ is the Chern-Simons form of the induced $SU(2)$-connection $\hat{\Th}$ on $\hat{P} = P \times_{\BT} SU(2)$,
\begin{equation*}
\al(\hat{\Th}) \, = \, \langle \hat{\Th} \wedge F_{\hat{\Th}} \rangle^{\flat} - \frac{1}{6} \langle \hat{\Th} \wedge [ \hat{\Th} , \hat{\Th}] \rangle
^{\flat} \, ,
\end{equation*}
and $\hat{s} : X \ra \hat{P}$ is a global section of the $SU(2)$-bundle $\hat{P}$.

We need to show that the definition of $S_{X,P}(\Th)$ does not depend on the choice of section $\hat{s} : X \ra \hat{P}$. 
This is a consequence of the following property of the $SU(2)$ Chern-Simons form $\al(\hat{\Th})$:

\begin{prop} \label{p:ssp}
If $\hat{s}, \hat{s}_1 : X \ra \hat{P}$ are two sections of the $SU(2)$-bundle $\hat{P}$ over $X$, then
\begin{equation*}
\int\limits_X \hat{s}^*_1 \al(\hat{\Th}) \, = \, \int\limits_X \hat{s}^* \al(\hat{\Th}) + \int\limits_X d \langle Ad_{a^{-1}} \, \hat{s}^* \hat{\Th} \wedge a^* \hat{\vartheta} \rangle^{\flat} - \int\limits_X \frac{1}{6} \, a^* \langle \hat{\vartheta} \wedge [ \hat{\vartheta} , \hat{\vartheta}] \rangle^{\flat} \, ,
\end{equation*} 
where $a: X \ra SU(2)$ is the map defined by $\hat{s}_1 (x) = \hat{s}(x) \cdot a(x)$, for any $x \in X$.
\end{prop}
\begin{proof}
The expression follows after basic computations from the standard relation
$\, \hat{s}^*_1 \, \hat{\Th} \, = \, Ad_{a^{-1}}( \hat{s} \, \hat{\Th}) + a^* \hat{\vartheta}$.
\end{proof}
On the right-hand side of the equation in (\ref{p:ssp}), the second integral vanishes by Stokes theorem and the assumption $\partial X = \emptyset$, while the last integral is an integer due to the normalization of the bilinear form $\langle  \cdot, \cdot \rangle^{\flat}$ on $\mathfrak{su}(2)$.
This proves that $S_{X,P}(\Th)$ is well-defined.

\begin{rmk}
If the bundle $P \ra X$ is trivializable, then let $s: X \ra P$ be a section and take $\hat{s} : X \ra \hat{P}$ to be $\hat{s} = \rho_P \circ s$.
Then we have
\begin{align*}
\hat{s}^* \al(\hat{\Th}) &= s^* \rho_P^* \al(\hat{\Th}) 
= s^* \langle \rho^*_P \hat{\Th} \wedge \rho_P^* F_{\hat{\Th}} \rangle^{\flat} -  \frac{1}{6} s^* \langle \rho_P^* \hat{\Th} \wedge [\rho_P^* \hat{\Th} ,\rho_P^* \hat{\Th}] \rangle^{\flat} \\
&= s^* \langle \rho_{*} \Th \wedge \rho_{*} F_{\Th} \rangle^{\flat} = s^* \langle \Th \wedge F_{\Th} \rangle
\end{align*}
and the Chern-Simons functional has as expected the expression
\begin{equation} \label{e:CStbdl}
S_{X,P}(\Th) \, = \, \int\limits_X s^* \langle \Th \wedge F_{\Th} \rangle \; \pmod{1}
\end{equation}
\end{rmk}

\begin{thm} \label{t:propCSnob}
The Chern-Simons functional $S_{X,P} : \CA_P \ra \BR/\BZ$ defined for a closed oriented 3-manifold $X$ and a principal $\BT$-bundle $P \ra X$ has the following properties:

\nin (a)  Functoriality \\
If $\phi : P' \ra P$ is a morphism of principal $\BT$-bundles covering an orientation preserving diffeomorphism $\bar{\phi} : X' \ra X$ and if $\Th$ is a connection on $P$, then
\begin{equation*}
S_{X',P'}(\phi^* \Th) \, = \, S_{X,P}(\Th)
\end{equation*}

\nin (b)  Orientation \\
If $-X$ denotes the manifold $X$ with the opposite orientation, then
\begin{equation*}
S_{-X,P}(\Th) \, = \, - S_{X,P}(\Th)
\end{equation*}

\nin (c)  Disjoint union \\
Let $X$ be the disjoint union $X =X_1 \sqcup X_2$ and $P =P_1 \sqcup P_2$ a principal $\BT$-bundle over $X$.
If $\Th_i$ are connections on $P_i \ra X_i$, then
\begin{equation*}
S_{X_1 \sqcup X_2, P_1 \sqcup P_2}(\Th_1 \sqcup \Th_2) \, = \, S_{X_1,P_1}(\Th_1) \, + \, S_{X_2,P_2}(\Th_2)
\end{equation*}
\end{thm}

\begin{proof}
(a) To the principal $\BT$-bundles $P$ and $P'$ there correspond the induced $SU(2)$-bundles $\hat{P} = P \times_{\BT} SU(2)$ and $\hat{P}' = P' \times_{\BT} SU(2)$, respectively.
The morphism of $\BT$-bundles $\phi : P' \ra P$ induces a morphism of $SU(2)$-bundles $\hat{\phi} : \hat{P}' \ra \hat{P}$ covering $\bar{\phi}: X' \ra X$.
It is defined by $\hat{\phi}([p',a]) = [\phi(p'),a]$.
For any section $\hat{s}' : X' \ra \hat{P}'$ we have a section $\hat{s} = \hat{\phi} \circ \hat{s}' \circ \bar{\phi}^{-1} : X \ra \hat{P}$.
Then we get 
\begin{equation*}
\begin{split}
S_{X,P}(\Th)\, &= \, \int\limits_X \hat{s}^* \al(\hat{\Th}) \; \pmod{1} \, = \, \int\limits_X (\bar{\phi}^{-1})^* \, \hat{s}^{\prime *} \, \hat{\phi}^* \, \al(\hat{\Th}) \; \pmod{1} \\
&= \, \int\limits_{X'} \hat{s}^{\prime *} \al(\hat{\phi}^* \hat{\Th}) \; \pmod{1}  \, = \, S_{X',P'}(\phi^* \Th) \, .
\end{split}
\end{equation*}
The last equality follows from the definition (\ref{d:CSnob}) and the fact that $\hat{\phi}^* \hat{\Th}$ is the $SU(2)$-connection on $\hat{P}'$ induced from the $\BT$-connection $\phi^* \Th$ on $P'$, as can be seen from the relation
\begin{equation*}
\begin{split}
\hat{\phi}^* \hat{\Th}_{(p,a)} \, &= \, Ad_{a^{-1}} \, \rho_{*}(\hat{\phi}^* \, pr_1^* \, \Th_p) + \hat{\phi}^* \, pr_2^* \, \hat{\vartheta}_a \\
 &= \, Ad_{a^{-1}} \, \rho_{*}(pr_1^* \, \phi^*  \Th_p) +  pr_2^* \, \hat{\vartheta}_a \, .
\end{split}
\end{equation*}

\nin (b) is a direct consequence of the definition of integration of differential forms on oriented manifolds.

\nin (c) Let $\hat{\Th} = \hat{\Th}_1 \sqcup \hat{\Th}_2$ be the extension of the $\BT$-connection $\Th = \Th_1 \sqcup \Th_2$ on $P = P_1 \sqcup P_2$ to a $SU(2)$-connection on the induced $SU(2)$-bundle $\hat{P} = \hat{P}_1 \sqcup \hat{P_2}$, where $\hat{P}_i = P_i \times_{\BT} SU(2)$.
Take sections $\hat{s}_i : X_i \ra \hat{P}_i$ and let $\hat{s} = \hat{s}_1 \sqcup \hat{s}_2 : X \ra \hat{P}$.
Then we have
\begin{equation*}
\begin{split}
S_{X,P}(\Th) \, &= \, \int\limits_X \hat{s}^* \al(\hat{\Th}) \; \pmod{1} \, = \, \int\limits_{X_1} \hat{s}_1^* \al(\hat{\Th}_1) + \int\limits_{X_2} \hat{s}_2^* \al(\hat{\Th}_2) \; \pmod{1} \\
&= \, S_{X_1,P_1}(\Th_1) + S_{X_2,P_2}(\Th_2)
\end{split}
\end{equation*}
\end{proof}
As a particular case of (\ref{t:propCSnob}(a))  we have:
\begin{prop} \label{p:ginv}
The functional $S_{X,P} : \CA_P \ra \BR/\BZ$ for a closed 3-manifold $X$ is invariant under the group of gauge transformations, that is,
\begin{equation*}
S_{X,P}(\Th \cdot u) \, = \, S_{X,P}(\Th) \, ,
\end{equation*}
for any $u \in \CG_X $. Hence $S_{X,P}$ descends to a functional on the quotient space $\CA_P/\CG_P$ of gauge equivalence classes of connections on $P$.
\end{prop}
\nin We also have the property:
\begin{prop} \label{p:statnb}
The stationary points of the functional $S_{X,P} : \CA_P \ra \BR/\BZ$ are the flat connections, that is,
$d \, S_{X,P}(\Th) =0$ if and only if $F_{\Th} = d \Th =0$.
\end{prop}
\begin{proof}
Since $\CA_P$ is an affine space, it suffices to consider the variation of $S_{X,P}$ along lines $\Th_t =\Th + t \om$, $\, t \in \BR$, with $\om \in 2 \pi \ti \Om^1(X;\BR)$.
According to Lemma (\ref{l:diffCS}) and since $\pX= \emptyset$, we have
\begin{equation*}
S_{X,P}(\Th_t) =  S_{X,P}(\Th)  + 2 t\int\limits_X \langle F_{\Th} \wedge \om \rangle + t^2\int\limits_X \langle \om \wedge d\om \rangle \; \pmod{1} .
\end{equation*}
Thus $\frac{d}{dt} \bigr|_{t=0} S_{X,P}(\Th_t) = 2 \int\limits_X \langle F_{\Th} \wedge \om \rangle$ which proves the proposition.
\end{proof}

\bsk


\subsection{The Chern-Simons functional for a 3-manifold with boundary}

Let $X$ be a compact oriented 3-manifold with $\partial X \neq \emptyset$.
Let $ P \ra X$ be a $\BT$-bundle with first Chern class torsion.
Then the $\BT$-bundle $\partial P \ra \partial X$ obtained by restricting $P$ to the boundary $\partial X$ is trivializable.
As before, we consider the induced $SU(2)$-bundle $\hat{P} = P \times_{\BT} SU(2)$ on $X$ and the bundle morphism $\rho_P : P \ra \hat{P}$.
For any section $s : \partial X \ra \partial P$, we obtain a section $\rho_P \circ s : \partial X \ra \partial \hat{P}$ of the restriction of the $SU(2)$-bundle $\hat{P}$ to $\partial X$.
Since the group $SU(2)$ has $\pi_0(SU(2)) = \pi_1(SU(2)) = \pi_2(SU(2)) =0$, we can extend the section $\rho_P \circ s$ over $\partial X$ to a global section $\hat{s} : X \ra \hat{P}$. 
We make the following definition.

\begin{defi} \label{d:CSb}
The Chern-Simons functional of a $\BT$-connection $\Th$ on $P \ra X$ and a section $s : \partial X \ra \partial P$ is defined by
\begin{equation*}
S_{X,P}(s,\Th) \, = \, \int\limits_X \hat{s}^* \al(\hat{\Th}) \; \pmod{1} \,,
\end{equation*}
where $\al(\hat{\Th})$ is the Chern-Simons form of the induced $SU(2)$-connection $\hat{\Th}$ on the bundle $\hat{P} = P \times_{\BT} SU(2)$.
\end{defi}

\nin We have to check that the definition of $S_{X,P}(s,\Th)$ does not depend on the choice of extension $\hat{s}$ of the section $\rho_P \circ s : \partial X \ra \partial \hat{P}$.
So let $\hat{s}_1 : X \ra \hat{P}$ be another extension.
Then there is a map $a : X \ra SU(2)$ such that $\hat{s}_1 (x) = \hat{s}(x) \cdot a(x)$ and $a(x) = e$, for all $x \in \pX$, where $e$ denotes the identity element in $SU(2)$.
From (\ref{p:ssp}) we get
\begin{equation} \label{e:well}
\int\limits_X \hat{s}^*_1 \al(\hat{\Th}) - \int\limits_X \hat{s}^* \al(\hat{\Th}) \, = \, \int\limits_{\partial X} \langle Ad_{a^{-1}} \hat{s}^* \hat{\Th} \wedge a^* \hat{\vartheta} \rangle^{\flat} - \int\limits_X \frac{1}{6} \, a^* \langle \hat{\vartheta} \wedge [ \hat{\vartheta} , \hat{\vartheta}] \rangle^{\flat} 
\end{equation}
That $S_{X,P}(s,\Th)$ is well-defined by (\ref{d:CSb}) follows now from the fact that, since $a\bigr|_{\partial X} = e$, the first integral on the r.h.s of (\ref{e:well}) vanishes and the second term is an integer.
To justify the last assertion we note the following:

\begin{rmk}
For each $a: X \ra SU(2)$, let $W(a)$ denote the functional
\begin{equation*}
W_{\partial X}(a) \, = \, - \int\limits_X \frac{1}{6} \, a^* \langle \hat{\vartheta} \wedge [ \hat{\vartheta} , \hat{\vartheta}] \rangle^{\flat}
\end{equation*}
Then, following the argument in (\cite{F}, \S 2), we let $X'$ be a compact oriented 3-manifold with $\partial X' = \partial X$ and $a' : X' \ra SU(2)$ a map such that $a' \bigr|_{\partial X'} = a \bigr|_{\partial X}$.
Let $\tilde{X} = X \cup (-X')$ be the closed 3-manifold obtained by gluing $X$ and $X'$ along their common boundary.
The maps $a$ and $a'$ patch together into a map $\tilde{a} : \tilde{X} \ra SU(2)$ and we have
\begin{equation} \label{e:difW}
W_{\partial X}(a) - W_{\partial X'}(a') \, = \, - \int\limits_{\tilde{X}}  \frac{1}{6} \, \tilde{a}^* \langle \hat{\vartheta} \wedge [ \hat{\vartheta} , \hat{\vartheta}] \rangle^{\flat} \, \in \BZ
\end{equation}
Since $\tilde{X}$ is closed the r.h.s. is an integer due to the normalization of $\langle \cdot, \cdot \rangle^{\flat}$.
Thus $W_{\partial X}(a)$ depends modulo integers only on the restriction of $a : X \ra SU(2)$ to $\partial X$.
\end{rmk}

\nin If $a \bigr|_{\partial X} = e$, we can choose $a' : X' \ra SU(2)$ with $a' \bigr|_{\partial X'} = a \bigr|_{\partial X}$ to be the constant map $a'(x') = e$ for all $x' \in X'$.
Then $W_{\partial X'}(a') =0$  and therefore $W_{\partial X}(a) \in \BZ$.

\begin{rmk}
As in (\ref{e:CStbdl}), if $P \ra X$ is a trivializable $\BT$-bundle and $s:X \ra P$ a global section, then for any connection $\Th$ on $P$ the Chern-Simons functional $S_{X,P}(s,\Th)$ is given by the expression
\begin{equation} \label{e:CSbtbdl}
S_{X,P}(s,\Th) \, = \, \int\limits_X s^* \langle \Th \wedge F_{\Th} \rangle \; \pmod{1}
\end{equation}
\end{rmk}

The following proposition describes the dependence of the Chern-Simons functional $S_{X,P}(s,\Th)$ on the section $s: \partial X \ra \partial P$.

\begin{prop} \label{p:sect}
Let $s, s_1 : \partial X \ra \partial P$ be two sections of the restriction of the $\BT$-bundle $P \ra X$ to $\partial X$ and let $\Th$ be a connection on $P$.
Then
\begin{equation*}
S_{X,P}(s_1,\Th) \, = \, S_{X,P}(s,\Th) + \int\limits_{\partial X} \langle s^* \Th \wedge \la^* \vartheta \rangle \; \pmod{1} \, ,
\end{equation*}
where $\la : \partial X \ra \BT$ is the map defined by $s_1(x) = s(x) \cdot \la(x)$, for all $x \in \pX$.
\end{prop}

\begin{proof}
The sections $\rho_P \circ s , \, \rho_P \circ s_1 : \partial X \ra \partial \hat{P}$ extend to sections $\hat{s}, \, \hat{s}_1 : X \ra \hat{P}$ of the induced $SU(2)$-bundle $\hat{P} = P \times_{\BT} SU(2)$.
The map $a : X \ra SU(2)$ defined by $\hat{s}_1(x) = \hat{s}(x) \cdot a(x)$, for any $x \in X$, restricts to $a \bigr|_{\partial X} = \rho \circ \la$ over the boundary of $X$.
The difference $S_{X,P}(s_1,\Th) - S_{X,P}(s,\Th)$ is given by an expression identical to the r.h.s. of the equation (\ref{e:well}).
The second term $W_{\partial X}(a)$ is also in this case an integer.
To see that choose in (\ref{e:difW}) a 3-manifold $X'$ with $\partial X' = \partial X$, for which the map $\la : \partial X \ra \BT$ extends to a map $\la' : X' \ra \BT$.
Let $a' = \rho \circ \la' : X' \ra SU(2)$.
Then, since $\rho^* \hat{\vartheta} = \rho_{*} \vartheta$, we get that
\begin{equation*}
W_{\partial X'}(a') = -\int\limits_{X'}  \frac{1}{6} \, a^{\prime *} \langle \hat{\vartheta} \wedge [ \hat{\vartheta} , \hat{\vartheta}] \rangle^{\flat}
= -\int\limits_{X'}  \frac{1}{6} \, \la^{\prime *} 
\langle \vartheta \wedge [ \vartheta , \vartheta] \rangle =0
\end{equation*}
Thus (\ref{e:difW}) implies that $W_{\partial X}(a) \in \BZ$.
The first term on the r.h.s. of (\ref{e:well}) gives
\begin{equation*}
\begin{split}
\int\limits_{\partial X} \langle Ad_{a^{-1}} \hat{s}^* \hat{\Th} \wedge a^* \hat{\vartheta} \rangle^{\flat} \, &= \, 
\int\limits_{\partial X} \langle Ad_{(\rho \circ \la)^{-1}} \, s^* \rho_P^* \hat{\Th} \wedge \la^* \rho^* \hat{\vartheta} \rangle^{\flat} \\
&= \, \int\limits_{\partial X} \langle \rho_{*} \, (Ad_{\la^{-1}} 
s^* \Th) \wedge \rho_{*} \, (\la^* \vartheta) \rangle^{\flat} \\
&= \, \int\limits_{\partial X} \langle s^* \Th \wedge \la^* \vartheta \rangle
\end{split}
\end{equation*}
which proves the proposition.
\end{proof}

\nin Similarly to Theorem (\ref{t:propCSnob}) we have

\begin{thm} \label{t:propCSb}
The Chern-Simons functional $S_{X,P}(\cdot,\cdot)$, defined for a compact oriented 3-manifold $X$ with nonempty boundary and a principal $\BT$-bundle $P \ra X$ with first Chern class torsion, has the following properties:

\nin (a)  Functoriality \\
Let $\phi : P' \ra P$ be a morphism of principal $\BT$-bundles covering an orientation preserving diffeomorphism $\bar{\phi} : X' \ra X$.
If $\Th$ is a connection on $P$ and $s' : \partial X' \ra \partial P'$ a section of $P'$ over the boundary of $X'$, then
\begin{equation*}
S_{X',P'}(s',\phi^* \Th) \, = \, S_{X,P}(\partial \phi \circ s' \circ \partial \bar{\phi}^{-1} ,\Th) \, ,
\end{equation*}
where $\partial \phi: \partial P' \ra \partial P$ and $ \partial \bar{\phi} : \pX' \ra \pX$ are restrictions to the boundary.

\nin (b)  Orientation \\
If $-X$ denotes the manifold $X$ with the opposite orientation and $s: \partial X \ra \partial P$, then
\begin{equation*}
S_{-X,P}(s,\Th) \, = \, - S_{X,P}(s,\Th)
\end{equation*}

\nin (c)  Disjoint union \\
Let $X$ be the disjoint union $X =X_1 \sqcup X_2$ and $P =P_1 \sqcup P_2$ a principal $\BT$-bundle over $X$ with first Chern class torsion.
If $\Th_i$ are connections on $P_i \ra X_i$ and $s_i :\partial X_i \ra \partial P_i$ sections over the boundary, then
\begin{equation*}
S_{X_1 \sqcup X_2, P_1 \sqcup P_2}(s_1 \sqcup s_2,\Th_1 \sqcup \Th_2) \, = \, S_{X_1,P_1}(s_1,\Th_1) \, + \, S_{X_2,P_2}(s_2,\Th_2)
\end{equation*}
\end{thm}

\begin{proof}
(a) Let $\hat{P} = P \times_{\BT} SU(2)$ and $\hat{P}' = P' \times_{\BT} SU(2)$ be the induced $SU(2)$-bundles, $\hat{\phi} : \hat{P}' \ra \hat{P}$ the induced bundle map and $\hat{\Th}$ the induced connection on $\hat{P}$.
The sections $s' : \partial X' \ra \partial P'$ and $s= \partial \phi \circ s'\circ \partial \bar{\phi}^{-1} : \partial X \ra \partial P$ determine sections $\rho_{P'} \circ s' : \partial X' \ra \partial \hat{P}'$ and $\rho_{P} \circ s : \partial X \ra \partial \hat{P}$.
Let $\hat{s}' : X' \ra \hat{P}'$ be an extension of $\rho_{P'} \circ s'$.
Then, using the fact that $\hat{\phi} \circ \rho_{P'} = \rho_{P} \circ \phi$, we can show that $\hat{s} = \hat{\phi} \circ \hat{s}' \circ \bar{\phi}^{-1} : X \ra \hat{P}$ is an extension of $\rho_{P} \circ s$.
Moreover, the $SU(2)$-connection on $\hat{P}'$ induced from the $\BT$-connection $\phi^* \Th$ on $P'$ is $\hat{\phi}^* \hat{\Th}$.
Applying the definition (\ref{d:CSb}) we get
\begin{equation*}
\begin{split}
S_{X,P}(s,\Th) \, &= \, \int\limits_X \hat{s}^* \al(\hat{\Th}) \; \pmod{1} \, = \, \int\limits_X (\bar{\phi}^{-1})^* \hat{s}^{\prime *} \hat{\phi}^* \al(\hat{\Th}) \; \pmod{1} \\
&= \, \int\limits_{X'} \hat{s}^{\prime *} \al(\hat{\phi}^* \hat{\Th}) \; \pmod{1} \, = \, S_{X',P'}(s', \phi^* \Th)
\end{split}
\end{equation*}
(b) and (c) are obvious.
\end{proof}
As a consequence of (\ref{t:propCSb}) and (\ref{p:sect}) we have:

\begin{prop} \label{p:gdep}
The functional $S_{X,P}(\cdot,\cdot)$ changes under the group of gauge transformations as
\begin{equation*}
S_{X,P}(s, \Th \cdot u) \, = \, S_{X,P}(s,\Th) + \int\limits_{\partial X} \langle s^* \Th \wedge u^* \vartheta \rangle \; \pmod{1} \, ,
\end{equation*}
for any $u \in \CG_X = \mathrm{Map}(X,\BT)$.
\end{prop}

\begin{proof}
Let $\phi_u : P \ra P$ denote the bundle automorphism determined by $u$. Then
\begin{alignat*}{2}
S_{X,P}(s,\Th \cdot u) \, &= \, S_{X,P}(s, \phi^*_u \Th) \, = \, S_{X,P}(\partial \phi_u \circ s, \Th) \qquad & &(\text{by} \, (\ref{t:propCSb}) (a)) \\
&= \, S_{X,P}(s,\Th) + \int\limits_{\partial X} \langle s^* \Th \wedge u^* \vartheta \rangle \; \pmod{1} \qquad & &(\text{by} \, (\ref{p:sect}))
\end{alignat*}
\end{proof}
We also note the following lemma:

\begin{lem} \label{l:diffCS}
Let $X$ be a compact oriented 3-manifold and $P \ra X$ a $\BT$-bundle with first Chern class torsion.
If $\Th$ and $\Th'$ are connections on $P$ and $s : \partial X  \ra \partial P$ a section over the boundary of $X$, then
\begin{align*}
S_{X,P}(s,\Th') - S_{X,P}(s,\Th) \, = \, &- \int\limits_{\partial X} s^* \langle \Th \wedge \om \rangle \\
&+ 2 \int\limits_X \langle F_{\Th} \wedge \om \rangle + \int\limits_X \langle \om \wedge d\om \rangle \; \pmod{1} \, ,
\end{align*}
where $\om = \Th' - \Th$.
\end{lem}

\begin{proof}
Let $\hat{\Th}$ and $\hat{\Th}'$ be the $SU(2)$-connections on
$\hat{P} = P \times_{\BT} SU(2)$ induced from the $\BT$-connections $\Th$ and $\Th'$ on $P$.
It follows from (\ref{e:indconn}) that
the $\mathfrak{su}(2)$-valued 1-form $\hat{\om} = \hat{\Th}' - \hat{\Th}$ on $\hat{P}$ is related to $\om$ by
\begin{equation} \label{e:om}
\hat{\om}_{[p,a]} \, = \, Ad_{a^{-1}}(\rho_{*} \, pr^*_1 \, \om_p) \,.
\end{equation}
and 
$\rho_P^* \hat{\om} \, = \, \rho_P^* (\hat{\Th}' - \hat{\Th}) \, = \, \rho_{*}(\Th' - \Th) \, = \, \rho_{*} \om $.
Routine algebraic computations give us for the difference of the Chern-Simons forms of the $SU(2)$-connections $\hat{\Th}'$ and $\hat{\Th}$ the expression:
\begin{equation*}
\begin{split}
\al(\hat{\Th}') - \al(\hat{\Th})\, = \, & \, -d \langle \hat{\Th} \wedge \hat{\om} \rangle^{\flat} + \, 2 \langle F_{\hat{\Th}} \wedge \hat{\om} \rangle^{\flat} \\
&+ \, \langle \hat{\om} \wedge (d \hat{\om} \, +\, [\hat{\Th},\hat{\om}]) \rangle^{\flat} + \frac{1}{3} \langle \hat{\om} \wedge [ \hat{\om},\hat{\om}] \rangle^{\flat}
\end{split}
\end{equation*}
According to (\ref{e:om}) the $\mathfrak{su}(2)$-valued 1-form $\hat{\om}$ comes from a $\text{Lie}\, \BT$-valued form $\om$, so we have $[\hat{\om}, \hat{\om}]=0$.
Thus the above expression reduces to
\begin{equation*}
\al(\hat{\Th}') - \al(\hat{\Th})\, =\, -d \langle \hat{\Th} \wedge \hat{\om} \rangle^{\flat} + \, 2 \langle F_{\hat{\Th}} \wedge \hat{\om} \rangle^{\flat} \, + \,
\langle \hat{\om} \wedge d \hat{\om} \rangle^{\flat}
\end{equation*}
The 3-forms $\langle F_{\hat{\Th}} \wedge \hat{\om} \rangle^{\flat}$ on $\hat{P}$ and $\langle F_{\Th} \wedge \om \rangle$ on $P$ are both lifts to $\hat{P}$ and $P$, respectively, of the same 3-form on $X$.
The same statement applies to the 3-forms $\langle \hat{\om} \wedge d \hat{\om} \rangle^{\flat}$ on $\hat{P}$ and $\langle \om \wedge d \om \rangle$ on $P$.
Letting $\hat{s}: X \ra \hat{P}$ be the extension to $X$ of the section $\rho_P \circ s : \partial X \ra \partial \hat{P}$, we can write
\begin{equation*}
\int\limits_X \hat{s}^* d \langle \hat{\Th} \wedge \hat{\om} \rangle^{\flat}
= \int\limits_{\partial X} s^* \rho_P^* \langle \hat{\Th} \wedge \hat{\om} \rangle^{\flat}
= \int\limits_{\partial X} s^* \langle \rho_{*} \Th \wedge \rho_{*} \om \rangle^{\flat}
= \int\limits_{\partial X} s^* \langle \Th \wedge \om \rangle
\end{equation*}
Putting all these results together, we find that the difference
\begin{equation*}
S_{X,P}(s,\Th') - S_{X,P}(s,\Th) = \int\limits_X \hat{s}^* \al(\hat{\Th}') - \int\limits_X \hat{s}^* \al(\hat{\Th})
\end{equation*}
is given by the expression in (\ref{l:diffCS}).
\end{proof}

\begin{lem} \label{l:Stokes}
Let $\Th$ be a connection on the $\BT$-bundle $P \ra X$ with $c_1(P)$ torsion and $s: \partial X \ra \partial P$ a section over the boundary of $X$.
If $u: X \ra \BT$ is an element of $\CG_X$, then
\begin{equation*}
\int\limits_X \langle F_{\Th} \wedge u^* \vartheta \rangle \, = \, \int\limits_{\partial X} \langle s^* \Th \wedge u^* \vartheta \rangle \quad \pmod{1}
\end{equation*}
\end{lem}

\begin{proof}
The fact that the integral $\int\limits_X \langle F_{\Th} \wedge u^* \vartheta \rangle$ depends modulo $\BZ$ only on the boundary data can be seen as follows.
Let $X'$ be a compact oriented 3-manifold with $\partial X' = \partial X$ and such that the map $u: \partial X \ra \BT$ extends to a map $u' : X' \ra \BT$.
Then glue the manifolds $X$ and $X'$ along their common boundary into the closed 3-manifold $\tilde{X} = X \sqcup (-X')$ and extend the connection $\partial \Th$ on $\partial X$ to a connection $\Th'$ over $X'$.
Let $\tilde{\Th} = \Th \sqcup \Th'$ and let $\tilde{u}: \tilde{X} \ra \BT$ denote the map with restrictions $\tilde{u} \bigr|_{X} = u$ and $\tilde{u} \bigr|_{X'} = u'$.
The integral
\begin{equation} \label{e:integer}
\int\limits_{\tilde{X}} \langle F_{\tilde{\Th}} \wedge \tilde{u}^* \vartheta \rangle \, = \, \int\limits_X \langle F_{\Th} \wedge u^* \vartheta \rangle \, - \, \int\limits_{X'} \langle F_{\Th'} \wedge u^{\prime *} \vartheta \rangle \; \in \BZ \, ,
\end{equation}
since $\frac{F_{\tilde{\Th}}}{2 \pi \ti}$ and $\frac{\tilde{u}^* \vartheta}{2 \pi \ti}$ represent integral cohomology classes.
Thus $\int\limits_X \langle F_{\Th} \wedge u^* \vartheta \rangle$ depends modulo $\BZ$ only on the restrictions $\partial \Th$ and $u \big|_{\partial X}$ to the boundary $\partial X$.
To prove the lemma take $\Th'$ extending $\partial \Th$ over $X'$ to be a connection on the trivial bundle $\pi' : P'=X' \times \BT \ra X'$.
Then if $s' : X' \ra P'$ is a section we have from (\ref{e:integer})
\begin{equation*}
\begin{split}
\int\limits_X \langle F_{\Th} \wedge u^* \vartheta \rangle \, &= \, \int\limits_{X'} \langle F_{\Th'} \wedge u^{\prime *} \vartheta \rangle \; \pmod{1} 
= \int\limits_{X'} s^{\prime *} \langle d \Th' \wedge \pi^{\prime *} u^{\prime *} \vartheta \rangle \; \pmod{1} \\
&= \, \int\limits_{X'} d \, s^{\prime *} \langle  \Th' \wedge \pi^{\prime *} u^{\prime *} \vartheta \rangle \; \pmod{1} 
= \int\limits_{\partial X'} \langle s^{\prime *}  \Th' \wedge u^{\prime *} \vartheta \rangle \; \pmod{1} \\
&= \, \int\limits_{\partial X} \langle s^*  \Th \wedge u^* \vartheta \rangle \; \pmod{1} \, ,
\end{split}
\end{equation*}
where $s: \partial X \ra \partial P$ is any section over the boundary.
\end{proof}

\begin{prop} \label{p:ginvf}
Let $P \ra X$ be a $\BT$-bundle with first Chern class torsion.
For any section $s:\pX \ra \partial P$ the functional $S_{X,P}(s,\cdot) : \CA_P^f \ra \BR/\BZ$ is invariant under the group $\CG_X$ of gauge transformations, that is,
\begin{equation*}
S_{X,P}(s,\Th \cdot u) \, = \, S_{X,P}(s,\Th) \, ,
\end{equation*}
for any $u \in \CG_X$. 
Hence $S_{X,P}(s,\cdot)$ descends to a functional on the moduli space $\CM_P$ of flat connections on $P$.
\end{prop}

\begin{proof}
If $\Th$ is flat, then Lemma (\ref{l:Stokes}) gives 
\begin{equation*}
\int\limits_{\pX} \langle s^* \Th \wedge u^* \vth \rangle \, = \, 0 \; \pmod{1} \, ,
\end{equation*}
for any $u \in \CG_X$. The proposition follows now from Prop. (\ref{p:gdep}) together with the above result.
\end{proof}
 
\bsk


\section{The prequantum line bundle} \label{s:line}

\subsection{The line bundle over the moduli space of flat connections}

As shown in Sect.\ref{s:moduli} to each closed oriented 2-dimensional manifold $\Si$ there corresponds the symplectic torus $(\CM_{\Si}, \om_{\Si})$ of gauge equivalence classes of flat $\BT$-connections on $\Si$.

Given a positive even integer $k$, we construct in this section a hermitian line bundle $\CL_{\Si}$ over $\CM_{\Si}$ with a unitary connection with curvature $- 2 \pi \ti k \om_{\Si}$.
In the language of geometric quantization such a line bundle is called a {\em prequantum line bundle} for the symplectic manifold $(\CM_{\Si},k \om_{\Si})$.
The set of prequantum line bundles for $(\CM_{\Si},k \om_{\Si})$ is a principal homogeneous space for the cohomology group $H^1(\CM_{\Si};\BT)$. Since $H^1(\CM_{\Si};\BT)$ is nontrivial, 
by general geometric quantization arguments the choice of the prequantum line bundle is not unique.
However, the 'Chern-Simons line construction' from (\cite{F}, \S 2), which we are going to apply here, singles out a prequantum line bundle $\CL_{\Si}$.

Let $ Q \ra \Si$ be a trivializable $\BT$-bundle over the closed oriented surface $\Si$.
To each connection $\eta$ on $Q$ we are going to associate a hermitian line $L_{\eta}$ by the following construction.
Let $\Ga(\Si;Q)$ stand for the space of sections of the line bundle $Q \ra \Si$.
Any two sections $s, s_1 \in \Ga(\Si;Q)$ are related by an element $\la : \Si \ra \BT$ of the group of gauge transformations $\CG_{\Si}$, defined by $s_1(x) = s(x) \cdot \la(x)$.
Thus we let $L_{\eta}$ be the space of functions $f : \Ga(\Si;Q) \ra \BC$ satisfying the relation
\begin{equation} \label{e:fs}
f(s_1) \, = \, f(s) \: c_{\Si}(s^* \eta , \la) \, ,
\end{equation}
where $c_{\Si}(s^* \eta , \la)$ is the $\BT$-valued cocycle
\begin{align} \label{e:cocycle}
&c_{\Si} \;  : \; 2 \pi \ti \, \Om^1(\Si;\BR) \times \CG_{\Si} \lra \BT \\
& c_{\Si}(\al, \la) \, = \, \te^{ \pi \ti k \int\limits_{\Si} \langle \al \wedge \la^* \vth \rangle} \, = \, 
\te^{ -\pi \ti k \int\limits_{\Si} \frac{\al}{2 \pi \ti} \wedge \la^* (\frac{\vth}{2 \pi \ti}) } \, . \notag
\end{align}
$c_{\Si}$ satisfies the cocycle identity
\begin{equation} \label{e:cocond}
c_{\Si} (\al, \la_1 \la_2) \, = \, c_{\Si}(\al, \la_1) \, c_{\Si}(\al+ \la_1^* \vth, \la_2)
\end{equation}
if and only if 
\begin{equation*}
\frac{k}{2} \int\limits_{\Si} \la_1^*(\frac{\vth}{2 \pi \ti}) \wedge \la_2^*(\frac{\vth}{2 \pi \ti}) \, \in \, \BZ \, ,
\end{equation*}
for any $\la_1, \la_2 \in \CG_{\Si}$.
This imposes the restriction $k \in 2 \BZ$.
Hence, for the remainder of this paper, we are going to assume that $k$ is a positive even integer.

An element $f \in  L_{\eta}$ is uniquely determined by  its value at a single point $s \in \Ga(\Si;Q)$.
Therefore $\dim_{\BC} L_{\eta} = 1$.
We also note that a  section $s : \Si \ra Q$ determines a trivialization $L_{\eta} \cong \BC$ by mapping $f \mapsto f(s)$.
The complex line $L_{\eta}$ has a hermitian inner product determined by the standard inner product of $\BC : \; (f_1,f_2) = \overline{f_1(s)} f_2(s)$.

As in \cite{F} we have the property that, as $\eta$ varies over the space $\CA_Q$ of connections on $Q$, the lines $L_{\eta}$ fit together smoothly into a hermitian line bundle $\CLH_Q \ra \CA_Q$.
Moreover, we have:

\begin{thm} \label{t:pline}
Let $\Si$ be a closed oriented $2$-manifold.
The assignment
\begin{equation*}
\eta \mapsto L_{\Si,\eta}=L_{\eta} \, ,
\end{equation*}
for $\eta$ connections on trivializable $\BT$-bundles over $\Si$, is smooth and satisfies: \\
(a) Functoriality \\
Let $\psi : Q' \ra Q$ be a morphism between trivializable $\BT$-bundles covering an orientation preserving diffeomorphism $\bar{\psi} : \Si' \ra \Si$ and let $\eta$ be a connection on $Q$.
Then there is an induced isometry
\begin{alignat*}{2}
\psi^* \, : \, & L_{\eta} && \lra L_{\psi^* \eta} \\
               &  f && \longmapsto \psi^* f \, ,
\end{alignat*}
where $\psi^* f \in  L_{\psi^* \eta}$ is defined by $(\psi^* f)(s') = f(\psi \circ s' \circ \bar{\psi}^{-1})$, for any $s' \in \Ga(\Si';Q')$.

\nin (b) Orientation\\
If $-\Si$ denotes the manifold $\Si$ with the opposite orientation, then there is a natural isometry
\begin{equation*}
L_{-\Si,\eta} \, \cong \, \overline{L_{\Si, \eta}}
\end{equation*}
(c) Disjoint union \\
If $\Si = \Si_1 \sqcup \Si_2$ is a disjoint union and $\eta_i$ are connections on trivializable $\BT$-bundles over $\Si_i$, then there is a natural isometry
\begin{equation*}
L_{\eta_1 \sqcup \eta_2} \, \cong \, L_{\eta_1} \otimes L_{\eta_2}
\end{equation*}
\end{thm}

\begin{proof}
(a) We have to check that $\psi^* f$ does indeed belong to the line $L_{\psi^* \eta}$.
So let $\la' :\Si' \ra \BT$ and let $s',s_1' \in \Ga(\Si';Q')$ be sections related by $s_1'(x') = s'(x') \cdot \la'(x')$, for all $x' \in \Si'$.
Then $(\psi \circ s_1' \circ \bar{\psi}^{-1})(x) = (\psi \circ s' \circ \bar{\psi}^{-1})(x) \cdot \la(x)$ at any $x \in \Si$, where $\la = \la' \circ \bar{\psi}^{-1} : \Si \ra \BT$.
Therefore we obtain 
\begin{equation*}
\begin{split}
(\psi^* f)(s_1') \, &= \, f(\psi \circ s_1' \circ \bar{\psi}^{-1}) \\
&= \, f(\psi \circ s' \circ \bar{\psi}^{-1}) \; c_{\Si}((\bar{\psi}^{-1})^* s^{\prime *} \psi^* \eta , \la) \qquad (\text{by definition (\ref{e:fs})}) \\
&= \, (\psi^* f)(s') \; c_{\Si'}(s^{\prime *}(\psi^* \eta), \la')
\qquad \quad (\text{$\bar{\psi}$ is orientation preserving})
\end{split}
\end{equation*}
(b) follows from the fact that the $\BT$-valued cocycle $c_{\Si}$ changes to the complex conjugate expression when the orientation of $\Si$ is reversed. \\
(c) $\eta_i$ are connections on trivializable $\BT$-bundles $Q_i \ra \Si_i$.
Let $s_i : \Si_i \ra Q_i$ be some sections and $\la_i : \Si_i \ra \BT$ gauge transformations and set $s = s_1 \sqcup s_2, \, \eta = \eta_1 \sqcup \eta_2$ and $\la = \la_1 \sqcup \la_2$.
Then $c_{\Si}(s^* \eta, \la) = c_{\Si_1}(s_1^* \eta_1, \la_1) \, c_{\Si_2}(s_2^* \eta_2, \la_2)$ and this induces the isometry
\begin{alignat*}{2}
L_{\eta_1} & \otimes L_{\eta_2} && \lra L_{\eta_1 \sqcup \eta_2} \\
f_1 & \otimes f_2 && \longmapsto f \, ,
\end{alignat*}
where $f(s) = f_1(s_1) f_2(s_2)$, for any section $s = s_1 \sqcup s_2$.
\end{proof}

Each section $s : \Si \ra Q$ induces a trivialization of the line bundle $\CLH_Q$, given by the unitary section
\begin{equation*}
\ga_s : \CA_Q \lra \CLH_Q
\end{equation*}
with $\ga_s(\eta) \in L_{\eta}$ defined, for each $\eta$, by $(\ga_s(\eta))(s) =1$.
We define the connection $\hat{\na}$ on $\CLH_Q$ by
\begin{equation} \label{e:defconn}
\hat{\na} \ga_s \, = \, - 2 \pi \ti k \, \hat{\theta}_s \, \ga_s \, ,
\end{equation}
where $\hat{\theta}_s$ is the 1-form on $\CA_Q$ given by
\begin{equation} \label{e:defconf}
(\dot{\eta} \, \contrac \,\hat{\theta}_s )_{\eta} \, = \, - \frac{1}{2} \int\limits_{\Si} s^* \langle \eta \wedge \dot{\eta} \rangle \, ,
\end{equation}
for any $\eta \in \CA_Q$ and $\dot{\eta} \in T_{\eta} \CA_Q$.

Let $\BC^{\times} = \BC \smallsetminus \{ 0 \}$ and let $\CLH_Q^{\times} \ra \CA_Q$ be the principal $\BC^{\times}$-bundle associated to $\pi : \CLH_Q \ra \CA_Q$ (the bundle $\CLH_Q^{\times}$ is obtained from $\CLH_Q$ by removing the zero section).
The connection form $- 2 \pi \ti k \,\hat{\al}$ on $\CLH_Q^{\times}$ corresponding to $\hat{\na}$ is defined by the expression
\begin{equation} \label{e:alpha}
\hat{\al} \, = \, \pi^* \hat{\theta}_s - \frac{1}{2 \pi \ti k} \, \Phi^*_s \, pr_2^* \vth \, ,
\end{equation} 
where $\Phi_s : \CLH_Q \ra \CA_Q \times \BC$ is the trivializing map defined by $\Phi_s(\ga_s(\eta)) = (\eta,1)$ and $pr_2: \CA_Q \times \BC \ra \BC$ the natural projection.
Obviously we have $\ga_s^* \hat{\al} = \hat{\theta}_s$.
If $\psi \in \CG_Q = \text{Aut}(Q)$, the new section $\psi \circ s : \Si \ra Q$ induces a new section $\ga_{\psi \circ s} : \CA_Q \ra \CLH_Q$ related to $\ga_s$ by
\begin{equation} \label{e:relsec}
\ga_{\psi \circ s}(\eta) \, = \, c_{\Si}(s^* \eta, \la_{\psi})^{-1} \, \ga_s(\eta) \, ,
\end{equation}
where $\la_{\psi} : \Si \ra \BT$ is the map associated to $\psi$.
The connection form $\hat{\theta}_{\psi \circ s}$ is related to $\hat{\theta}_s$ by
\begin{equation*}
\begin{split}
(\dot{\eta} \, \contrac \, \hat{\theta}_{\psi \circ s})_{\eta} \, &= \,
- \frac{1}{2} \, \int\limits_{\Si} s^* \psi^* \langle \eta \wedge \dot{\eta} \rangle \, = \, 
- \frac{1}{2} \, \int\limits_{\Si} \langle (s^* \eta + \la_{\psi}^* \vth) \wedge s^* \dot{\eta} \rangle  \\
&= \, (\dot{\eta} \, \contrac \, \hat{\theta}_s)_{\eta} +
 \frac{1}{2} \, \int\limits_{\Si} \langle s^* \dot{\eta} \wedge \la_{\psi}^* \vth \rangle
\end{split}
\end{equation*}
Taking the differential of the cocycle $c_{\Si}(s^* \eta,\la_{\psi})$ with respect to $\eta$ and using the above expression, we find that
\begin{equation*}
(\hat{\theta}_{\psi \circ s})_{\eta} \, = \, (\hat{\theta}_s)_{\eta} + \frac{1}{2 \pi \ti k} \, \frac{d \, c_{\Si}(s^* \eta, \la_{\psi})}{c_{\Si}(s^* \eta, \la_{\psi})}
\end{equation*}
Together with (\ref{e:relsec}) this implies that
\begin{equation*}
\pi^* \hat{\theta}_{\psi \circ s} - \frac{1}{2 \pi \ti k} \Phi^*_{\psi \circ s}\, pr_2^* \vth \, = \, \pi^* \hat{\theta}_s - \frac{1}{2 \pi \ti k} \Phi^*_s \,pr_2^* \vth \, ,
\end{equation*}
which shows that the definition (\ref{e:alpha}) of $\hat{\al}$ does not depend on the choice of trivialization of $\CLH_Q \ra \CA_Q$.

Since $\ga_s$ is a unitary section and $\hat{\theta}_s$ is real, the connection $\hat{\na}$ is compatible with the hermitian structure on $\CLH_Q$.
The curvature of $\hat{\na}$ is $- 2 \pi \ti k \om_{\Si}$, where
\begin{equation} \label{e:curv}
\om_{\Si}(\dot{\eta}, \dot{\eta}') \, = \, (d \hat{\theta}_s)(\dot{\eta}, \dot{\eta}') \, = \, - \int\limits_{\Si} \langle \dot{\eta} \wedge \dot{\eta}' \rangle 
\end{equation}
is the standard symplectic form on $T \CA_Q \cong 2 \pi \ti \, \Om^1(\Si;\BR)$.

It follows  from (\ref{t:pline} (a)) that the action of the group of gauge transformations $\CG_Q$ on $\CA_Q$ lifts to the line bundle $\CLH_Q$ preserving the hermitian metric.
The lift of $\CG_Q$ to $\CLH_Q$ preserves also the connection $\hat{\na}$.
To prove this we consider the induced action on sections of $\CLH_Q \ra \CA_Q$:
\begin{equation*}
(\psi^* \cdot \ga_s)(\eta) \, = \, \psi^* \ga_s(\psi^{* -1} \eta) \, .
\end{equation*}
A simple application of (\ref{t:pline} (a)) and the definition (\ref{e:fs}) gives that
\begin{equation} \label{e:actsect}
\psi^* \cdot \ga_s \, = \, \ga_{\psi^{-1} \circ s} \, .
\end{equation}
The action of $\CG_Q$ on $\CLH_Q$ preserves $\hat{\na}$ if 
\begin{equation*}
\Big[ \hat{\na}_{\psi^* \dot{\eta}} ( \psi^* \cdot \ga_s) \Big] \, = \, \psi^* \Big[ \big( \hat{\na}_{\dot{\eta}} \ga_s \big) (\eta) \Big] \, .
\end{equation*}
Using (\ref{e:defconn}), (\ref{e:defconf}) and (\ref{e:actsect}), a routine check shows that the above equation is indeed satisfied.
The $\CG_Q$-action on $\CA_Q$ preserves the symplectic form $\om_{\Si}$.
This follows from the fact that $-2 \pi \ti k \om_{\Si}$ is the curvature of $\hat{\na}$.

Since the action of $\CG_Q$ on $\CA_Q$ is symplectic and lifts to the line bundle $\pi : \CLH_Q \ra \CA_Q$ preserving the hermitian metric and the connection $\hat{\na}$, there exists a moment map for this action. 
Following the general construction in \cite{K}, we give below the explicit form of this moment map.
The $\CG_Q$-action associates to each element $\xi : \Si \ra \text{Lie} \BT$ in $\text{Lie} \,\CG_Q$ real vector fields $X_{\xi}$ on $\CA_Q$ and $V_{\xi}$ on $\CLH_Q$, with $\pi_* V_{\xi} = X_{\xi}$.
The vector field $V_{\xi}$ preserves the connection form $- 2 \pi \ti k \,\hat{\al}$ on $\CLH_Q^{\times}$, i.e. the Lie derivative $L_{V_{\xi}} \hat{\al} =0$.
Together with the fact that $V_{\xi}$ preserves the hermitian metric, this implies that the function $\hat{\mu}_{\xi}$ on $\CA_Q$, defined by $\hat{\mu}_{\xi} \circ \pi = V_{\xi} \, \contrac \, (k \hat{\al})$, is real-valued and satisfies $X_{\xi} \, \contrac \, (k \om_{\Si}) + d \hat{\mu}_{\xi} = 0$.
The map $\mu : \CA_Q \ra (\text{Lie} \,\CG_Q)^*$, $\, \eta \mapsto \mu(\eta)$ with $(\mu(\eta))(\xi) = \hat{\mu}_{\xi}(\eta)$, is the {\em moment map} of the $\CG_Q$-action on $\CA_Q$.
Using (\ref{e:alpha}) we find:
\begin{equation*}
\begin{split}
\hat{\mu}_{\xi} \circ \pi \, = \, V_{\xi} \, \contrac \, (k \hat{\al}) \, &= \, V_{\xi} \, \contrac \, \pi^* (k \hat{\theta}_s) - \frac{1}{2 \pi \ti } V_{\xi} \, \contrac \, \Phi^*_s \, pr_2^* \vth \\
&= \, k (X_{\xi} \contrac \hat{\theta}_s) \circ \pi - \frac{1}{2 \pi \ti } V_{\xi} \, \contrac \, \Phi^*_s \, pr^*_2 \vth \, .
\end{split}
\end{equation*}
A direct computation gives
\begin{equation*}
(X_{\xi})_{\eta} \, = \, \frac{d}{dt} \Big|_{t=0} \eta \cdot \te^{t \xi} \, = \, d \xi
\end{equation*}
and
\begin{equation*}
V_{\xi} \, \contrac \, (\Phi^*_s \, pr^*_2 \vth) \, = \, \Big[ (pr_{2*} \Phi_{s*} V_{\xi}) \, \contrac \, \vth \Big] \circ pr_2 \circ \Phi_s \, = \, -2 \pi \ti k (d \xi \, \contrac \, \hat{\theta}_s) \circ \pi \, .
\end{equation*}
Therefore we obtain for the moment map the expression
\begin{equation*}
\hat{\mu}_{\xi}(\eta) \, = \, 2 k (d\xi \, \contrac \, \hat{\theta}_s)_{\eta} 
= \, - k \int\limits_{\Si} \langle s^* \eta \wedge d \xi \rangle \, = \, 
 - k\int\limits_{\Si} \langle s^* d\eta \wedge  \xi \rangle \, .
\end{equation*}

The preimage $\mu^{-1}(0)$, which is $\CG_Q$-invariant, is the space $\CA_Q^f$ of flat connections.
Hence the space $\CM_Q$ of equivalence classes of flat connections is the symplectic quotient $\CM_Q = \CA_Q / \!/ \CG_Q = \mu^{-1}(0) / \CG_Q$.
Since $\CG_Q$ preserves $\om_{\Si}$ and $\om_{\Si}(X_{\xi}, \dot{\eta}) =0$, for any element $\xi$ in $\text{Lie} \, \CG_Q$ and any $\dot{\eta} \in T \CA_Q^f$, the symplectic form $\om_{\Si}$ on $\CA_Q$ pushes down to a symplectic form on $\CM_Q$, which we continue to denote by $\om_{\Si}$. 

The group $\CG_Q$ does not act freely on $\CA_Q$.
However, the stabilizer at a point $\eta \in \CA_Q$ is the subgroup $Z \subset \CG_Q$ of constant gauge transformations 
$Z \cong \{ \la : \Si \ra \BT \mid \la = \text{constant} \} \cong \BT$ and $Z$ acts trivially on $L_{\eta}$.
Therefore the line bundle $\CLH_Q \ra \CA_Q^f=\mu^{-1}(0)$ pushes down to a line bundle $\CL_Q \ra \CM_Q$.
We have the identification
\begin{equation} \label{e:invsect}
\Ga(\CM_Q; \CL_Q) \, = \, \Ga(\CA_Q^f ; \CLH_Q)^{\CG_Q}
\end{equation}
of the space of sections of $\CL_Q \ra \CM_Q$ with the space of $\CG_Q$-invariant sections of $\CLH_Q \ra \CA_Q^f$.
The hermitian metric and connection on $\CLH_Q$ also push down to $\CL_Q$.
A $\CG_Q$-invariant section $\ga$ of $\CLH_Q \ra \CA_Q^f$ satisfies $\hat{\na}_{X_{\xi}} \ga =0$, for all $\xi : \Si \ra \text{Lie} \BT$.
Having in view the identification (\ref{e:invsect}), the push-down connection $\na$ on the line bundle $\CL_Q \ra \CM_Q$ is defined by
\begin{equation*}
\na_X \ga \, = \, \hat{\na}_{\hat{X}} \ga \, , 
\end{equation*}
for any $X \in T \CM_Q$, where $\hat{X} \in T \CA_Q^f$ is any vector which maps onto $X$ under the quotient map $\mu^{-1}(0) \ra \CM_Q$.

\begin{rmk}
Let us assume for the moment that in constructing the hermitian line bundle $\CLH_Q \ra \CA_Q$ we take the cocycle $c_{\Si}$ to have, instead of (\ref{e:cocycle}), the more general form
\begin{equation*}
c_{\Si}(\al, \la) \, = \, \vep_{\Si}([\la]) \te^{\pi \ti k \int\limits_{\Si} \langle \al \wedge \la^* \vth \rangle } \, , 
\end{equation*}
where $[\la]$ denotes the homotopy class of the map $\la : \Si \ra \BT$.
In order to satisfy the cocycle condition (\ref{e:cocond}) the $\BT$-valued multiplier $\vep_{\Si}$ must have the property
\begin{equation*}
\vep_{\Si}([\la_1][\la_2]) \, = \, \vep_{\Si}([\la_1]) \, \vep_{\Si}([\la_2])
\,\te^{\pi \ti k \int\limits_{\Si} \langle \la_1^* \vth \wedge \la_2^* \vth \rangle }
\end{equation*}
and $k$ can be any integer.
All what has been said up to this point, except (\ref{t:pline} (a)), remains true.
We get a hermitian line bundle over $\CA_Q$ with a unitary connection with curvature $- 2 \pi \ti k \om_{\Si}$ and the action of $\CG_Q$ lifts preserving the metric and the connection.
The restriction of this line to $\CA_Q^f$ pushes down to a prequantum line bundle over $\CM_Q$.
The initial choice (\ref{e:cocycle}) for $c_{\Si}$, in which $\vep_{\Si}$ is taken to be the trivial multiplier (which forces $k$ to be even), is dictated by the requirement to satisfy also the functoriality property (\ref{t:pline} (a)).
This can be seen immediately if we follow the proof given for (\ref{t:pline} (a)).
Only for this choice do we get a prequantum line bundle over the moduli space $\CM_{\Si}$ for which there is a lift of the group of orientation preserving diffeomorphisms of $\Si$.
We encountered the same situation in \cite{Ma} in the quantization of symplectic tori, where a prequantum line bundle was picked out by the requirement to have a lift of the group of symmetries of the torus.
In the present case the prequantum line bundle over $\CM_{\Si}$ picked out by the above stated reasons is also closely connected to the Chern-Simons theory in 3 dimensions.
This connection will become clear later in this section.
\end{rmk}

Given any two trivializable $\BT$-bundles $Q$ and $Q'$ over $\Si$, let $\psi : Q' \ra Q$ be a bundle isomorphism.
It follows from (\ref{t:pline} (a)) that the induced isomorphism between the spaces of connections lifts to an isomorphism of line bundles
\begin{equation} \label{e:cdlbdl}
\begin{CD}
\CLH_Q @>{\psi^*}>> \CLH_{Q'} \\
@VVV @VVV \\
\CA_Q @>{\psi^*}>> \CA_{Q'}
\end{CD}
\end{equation}
To each gauge transformation $\phi \in \CG_Q$ there corresponds a gauge transformation $\psi^{-1} \circ \phi \circ \psi \in \CG_{Q'}$, so that for any connection $\eta \in \CA_Q$ we have the commutative diagram
\begin{equation} \label{e:cdlines}
\begin{CD}
L_{\eta}  @>{\psi^*}>> L_{\psi^* \eta} \\
@V{\phi^*}VV @VV{(\psi^{-1} \circ \phi \circ \psi)^*}V \\
L_{\phi^* \eta} @>{\psi^*}>> L_{\psi^* \phi^* \eta}
\end{CD}
\end{equation}
of line isomorphisms defined by (\ref{t:pline}(a)).
The line bundle isomorphism (\ref{e:cdlbdl}) pushes down to an isomorphism between the quotient line bundles by the groups of gauge transformations
\begin{equation}
\begin{CD}
\CL_Q = \CLH_Q/ \CG_Q @>>> \CL_{Q'} = \CLH_{Q'}/ \CG_{Q'} \\
@VVV @VVV \\
\CA_Q / \CG_Q @>>> \CA_{Q'}/ \CG_{Q'}
\end{CD}
\end{equation}
It follows from (\ref{e:cdlines}) that the above isomorphism is canonical, that is, independent of the choice of $\BT$-bundle isomorphism $\psi: Q' \ra Q$.
Restricting the above considerations to the subspaces of flat connections  we obtain that, for any two trivializable $\BT$-bundles $Q$ and $Q'$ over $\Si$, there is a canonical isomorphism of hermitian line bundles with connections
\begin{equation} \label{e:equivlines}
\begin{CD}
\CL_Q @>>> \CL_{Q'} \\
@VVV @VVV \\
\CM_Q @>>> \CM_{Q'}
\end{CD}
\end{equation}
Since $\CM_{\Si} \cong \CM_Q$ for any trivializable $\BT$-bundle $Q \ra \Si$, the above shows that we get a hermitian line bundle 
\begin{equation*}
\CL_{\Si} \ra \CM_{\Si}
\end{equation*}
with a unitary connection with curvature $- 2 \pi \ti k \om_{\Si}$, with $\om_{\Si}$ the standard symplectic form on $\CM_{\Si}$.

\bsk


\subsection{The Chern-Simons section}

We turn now our attention to the 3-dimensional Chern-Simons theory.
Let $X$ be a compact oriented 3-manifold with nonempty boundary and $P \ra X$ a principal  $\BT$-bundle with $c_1(P)$ a torsion class.
If $\Th$ is a connection on $P$, then Prop. (\ref{p:sect}) shows that the function
\begin{alignat*}{2}
\te^{\pi \ti k S_{X,P}(\Th)} \, :  & \: \Ga(\pX;Q) && \lra \, \BC \\
&\quad  s && \longmapsto \, \te^{\pi \ti k S_{X,P}(\Th)}(s) = \te^{\pi \ti k S_{X,P}(s,\Th)} \notag
\end{alignat*}
is an element of norm 1 in the line $L_{\partial \Th}$ attached to the restriction of $\Th$ to $\pX$
\begin{equation} \label{e:exps}
\te^{\pi \ti k S_{X,P}(\Th)} \, \in \, L_{\partial \Th}
\end{equation}
If $\pX = \emptyset$, then $\te^{\pi \ti k S_{X,P}(\Th)} \, \in \, \BC$, so we set $L_{\emptyset} = \BC$ to account for this case too.

\begin{thm} \label{t:CSsect}
Let $X$ be a compact oriented 3-manifold.
The assignment
\begin{equation*}
\Th \lra \te^{\pi \ti k S_{X,P}(\Th)} \, \in \, L_{\partial \Th} \, ,
\end{equation*}
for $\Th$ connections on $\BT$-bundles $P \ra X$ with $c_1(P)$ torsion, is smooth and satisfies: \\
(a) Functoriality \\
Let $\phi : P' \ra P$ be a morphism of principal $\BT$-bundles covering an orientation preserving diffeomorphism $\bar{\phi} : X' \ra X$.
If $\Th$ is a connection on $P$, then
\begin{equation*}
(\partial \phi)^* \te^{\pi \ti k S_{X,P}(\Th)} \, = \, \te^{\pi \ti k S_{X',P'}(\phi^* \Th)} 
\end{equation*}
under the induced isometry $(\partial \phi)^* : L_{\partial \Th} \ra L_{\partial \phi^* \Th}$.
In particular if $\pX = \emptyset$, then
\begin{equation*}
\te^{\pi \ti k S_{X,P}(\Th)} \, = \, \te^{\pi \ti k S_{X',P'}(\phi^* \Th)} 
\end{equation*}
(b) Orientation \\
If $-X$ is the manifold $X$ with the opposite orientation, then
\begin{equation*}
\te^{\pi \ti k S_{-X,P}(\Th)} \, = \, \overline{\te^{\pi \ti k S_{X,P}(\Th)}} \, \in \, \overline{L_{\partial \Th}}
\end{equation*}
(c) Disjoint union \\
Let $X$ be the disjoint union $X = X_1 \sqcup X_2$.
If $\Th_i$ are connections on $\BT$-bundles $P_i \ra X_i$, then we have the identification
\begin{equation*}
\te^{\pi \ti k S_{X_1 \sqcup X_2,P_1 \sqcup P_2}(\Th_1 \sqcup \Th_2)} \, = \,
\te^{\pi \ti k S_{X_1,P_1}(\Th_1)}  \otimes \te^{\pi \ti k S_{X_2,P_2}(\Th_2)} 
\end{equation*}
under the isomorphism $L_{\partial \Th_1 \sqcup \partial \Th_2} \cong L_{\partial \Th_1} \otimes L_{\partial \Th_2}$.

\nin (d) Gluing \\
Let $X$ be a compact oriented 3-manifold and let $\XC$ be the 3-manifold obtained by cutting $X$ along a closed oriented 2-dimensional submanifold $\Si$.
Then $\pXC = \pX \sqcup (-\Si) \sqcup \Si$.
Let $\Th$ be a connection on a $\BT$-bundle $P \ra X$, with $\eta$ the restriction of $\Th$ to $\Si$, and let $\Th^{cut}$ denote the induced connection on the pullback bundle $P^{cut} = g^* P$ under the gluing map $g : \XC \ra X$.
If
\begin{equation*}
\mathrm{Tr}_{\eta} \, : \, L_{\partial \Th^{cut}} \cong L_{\partial \Th} \otimes \overline{L_{\eta}} \otimes L_{\eta} \lra L_{\partial \Th}
\end{equation*}
is the contraction map using the hermitian inner product in $L_{\eta}$, then
\begin{equation} \label{e:tr}
\mathrm{Tr}_{\eta} \Big( \te^{\pi \ti k S_{\XC,P^{cut}}(\Th^{cut})} \Big) \, = \, \te^{\pi \ti k S_{X,P}(\Th)}
\end{equation}
\end{thm}

\begin{proof}
(a) follows from (\ref{t:pline} (a)) and (\ref{t:propCSb} (a)).\\
(b) follows from (\ref{t:pline} (b)) and (\ref{t:propCSb} (b)).\\
(c) follows from (\ref{t:pline} (c)) and (\ref{t:propCSb} (c)).\\ 
(d) Let $g : P^{cut} \ra P$ be the bundle map from the pullback bundle $P^{cut} = g^* P$ to $P$, covering the gluing map $g : \XC \ra X$.
Then $\Th^{cut} = g^* \Th$ and $\partial \Th^{cut} = \partial \Th \sqcup \eta \sqcup \eta$.
Under the isometry $L_{\partial \Th^{cut}} \ra L_{\partial \Th} \otimes \overline{L_{\eta}} \otimes L_{\eta}$ the element $\te^{\pi \ti k S_{\XC,P^{cut}}(\Th^{cut})} \in L_{\partial \Th^{cut}}$ is mapped to 
\begin{equation} \label{e:3lines}
\frac{\te^{\pi \ti k S_{\XC,P^{cut}}(s^{cut},\Th^{cut}) } }{\te^{\pi \ti k S_{X,P}(s,\Th)}  \overline{f(\si)} f(\si) } \, \Big(\te^{\pi \ti k S_{X,P}(\Th)} \otimes \bar{f} \otimes f \Big) \, \in \,  
L_{\partial \Th} \otimes \overline{L_{\eta}} \otimes L_{\eta} \, ,
\end{equation}
where $s, \si$ and $s^{cut}$ are sections $s: \pX \ra \partial P, \, \si: \Si \ra P \bigr|_{\Si}$ and $s^{cut} = s \sqcup \si \sqcup \si : \pXC \ra \partial P^{cut}$.
To prove (\ref{e:tr}) we have to show that 
\begin{equation} \label{e:eq}
\te^{\pi \ti k S_{\XC,P^{cut}}(s^{cut},\Th^{cut})} \, = \, \te^{\pi \ti k S_{X,P}(s,\Th)} \, .
\end{equation}
Let $\rho_P : P \ra \hat{P}=P \times_{\BT} SU(2)$ and $\rho_{P^{cut}} : P^{cut} \ra \hat{P}^{cut}=P^{cut} \times_{\BT} SU(2) = g^* \hat{P}$
be the bundle maps to the induced $SU(2)$-bundles.
There is a $SU(2)$-bundle morphism $\hat{g} : \hat{P}^{cut} \ra \hat{P}$ so that $\rho_P \circ g = \hat{g} \circ \rho_{P^{cut}}$.
Let $\hat{s} : X \ra \hat{P}$ be the extension of $\rho_P \circ s : \pX \ra \partial \hat{P}$ and $\hat{s}^{cut} : X \ra \hat{P}^{cut}$ the pullback of the section $\hat{s}$.
Then $ \hat{g} \circ \hat{s}^{cut} = \hat{s} \circ g$ and $\hat{s}^{cut} \bigr|_{\pXC} = \rho_{P^{cut}} \circ s^{cut}$.
If $\hat{\Th}^{cut} = \hat{g}^* \hat{\Th}$ then $\rho_{P^{cut}}^* \hat{\Th}^{cut} = \rho_* \Th^{cut}$, so that $\hat{\Th}^{cut}$ is the $SU(2)$-connection induced from $\Th^{cut}$.
Therefore, using the definition (\ref{d:CSb}) of the Chern-Simons functional, we can write
\begin{equation*}
\begin{split}
S_{\XC,P^{cut}}(s^{cut},\Th^{cut})  
&=  \int\limits_{\XC} \hat{s}^{cut *} \al(\hat{\Th}^{cut})  \pmod{1}  
=  \int\limits_{\XC} \hat{s}^{cut *} \al(\hat{g}^* \hat{\Th})  \pmod{1} \\
&=  \int\limits_{\XC} g^* \hat{s}^* \al(\hat{\Th})  \pmod{1}  
=  \int\limits_{X} \hat{s}^* \al(\hat{\Th})  \pmod{1} \\
&=  \, S_{X,P}(s,\Th) 
\end{split}
\end{equation*}
which proves (\ref{e:eq}).
Then, taking in (\ref{e:3lines}) the hermitian inner product in $L_{\eta}$, i.e. $ (f,f)=\overline{f(\si)} f(\si)$, we get the announced result (\ref{e:tr}).
\end{proof}

Consider now the restriction map $r: \CA_P \ra \CA_{\partial P}$ which sends a connection on the $\BT$-bundle $P \ra X$ to its restriction over $\pX$.
As previously shown, to the trivializable $\BT$-bundle $\partial P \ra \pX$ there corresponds a hermitian line bundle $\pi : \CLH_{\partial P} \ra \CA_{\partial P}$ with connection $\hat{\na}$ and these pull back under $r$ to a hermitian line bundle $r^* \CLH_{\partial P} \ra \CA_P$ with connection $r^* \hat{\na}$.
The restriction of $r^* \CLH_{\partial P}$ to the subspace $\CA^f_P$ of flat connections is a {\em flat} line bundle since
\begin{equation*}
(r^* \om_{\Si})(\dot{\Th}, \dot{\Th}')
= \om_{\Si}(r_* \dot{\Th}, r_* \dot{\Th}')
= - \int\limits_{\pX} \langle \partial \dot{\Th} \wedge \partial \dot{\Th}' \rangle 
= - \int\limits_X d \langle \dot{\Th} \wedge \dot{\Th}' \rangle =0
\end{equation*}
for any $\dot{\Th}, \dot{\Th}' \in T\CA^f_P$.

As $\CG_P \subset \CG_{\partial P}$, the action of the group of gauge transformations $\CG_P$ on $\CA_P$ lifts to $r^* \CLH_{\partial P}$.
Since $r^* \CLH_{\partial P} = \{ (\Th,f) \in  \CA_P \times \CLH_{\partial P} \mid r(\Th) = \pi(f) \}$,
we see from (\ref{e:exps}) and (\ref{t:CSsect} (a)) that 
\begin{align} \label{e:CSsect}
\hat{\si}_P \, : \, \CA_P & \lra r^* \CLH_{\partial P} \\
\Th & \longmapsto (\Th, \te^{\pi \ti k S_{X,P}(\Th)}) \notag
\end{align}
is a nowhere-zero $\CG_P$-invariant section, i.e. for any $\phi \in \CG_P$,
\begin{equation*} 
\begin{split}
\hat{\si}_P(\Th) \cdot \phi \, &= \, (\phi^* \Th, (\partial \phi)^*  \te^{\pi \ti k S_{X,P}(\Th)} ) \\
&= \, (\phi^* \Th,  \te^{\pi \ti k S_{X,P}(\phi^* \Th)} ) \, = \, \hat{\si}_P(\phi^* \Th)
\end{split}
\end{equation*}
Moreover we have:

\begin{prop}
The section $\hat{\si}_P : \CA_P^f \ra r^* \CLH_{\partial P}$ is covariantly constant.
\end{prop}

\begin{proof}
Let $s: \pX \ra \partial P$ be a section.
The induced section $\ga_s : \CA_{\partial P} \ra \CLH_{\partial P}$ pulls back to $(r^* \ga_s)(\Th) = (\Th, \ga_s(\partial \Th)) : \CA_P \ra r^* \CLH_{\partial P}$.
Since $ \te^{\pi \ti k S_{X,P}(\Th)} =  \te^{\pi \ti k S_{X,P}(s,\Th)} \ga_s(\partial \Th)$, we can write $\hat{\si}_P = \varphi \, (r^* \ga_s)$ with
$\varphi(\Th)=  \te^{\pi \ti k S_{X,P}(s,\Th)}$.
By differentiating $\varphi$ with respect to $\dot{\eta} \in T_{\Th} \CA^f_P$, we find
\begin{equation*}
\begin{split}
(\dot{\eta} \varphi)(\Th) \, &= \, \frac{d}{dt} \Big|_{t=0} \varphi(\Th + t \dot{\eta})  
= \, \frac{d}{dt} \Big|_{t=0} \te^{ \pi \ti k  S_{X,P}(s,\Th + t \dot{\eta}) } \\
&= \, \frac{d}{dt} \Big|_{t=0} \te^{ \pi \ti k [ S_{X,P}(s,\Th) - t \int\limits_{\pX} s^* \langle \Th \wedge \dot{\eta} \rangle + 2 t \int\limits_X \langle d \Th \wedge \dot{\eta} \rangle + t^2 \int\limits_X \langle \dot{\eta} \wedge d \dot{\eta} \rangle ] } \\
&= \, \Big[ - \pi \ti k \int\limits_{\pX} s^* \langle \Th \wedge \dot{\eta} \rangle \Big] \, \varphi(\Th) \, .
\end{split}
\end{equation*}
On the other hand
\begin{equation*}
(r^* \hat{\na})_{\dot{\eta}} (r^* \ga_s) \, = \, - 2 \pi \ti k \, (\dot{\eta} \, \contrac \, r^* \hat{\theta}_s) \, (r^* \ga_s) \, .
\end{equation*}
Combining these two results we obtain 
\begin{equation*}
\begin{split}
(r^* \hat{\na})_{\dot{\eta}} \hat{\si}_P \, &= \, (r^* \hat{\na})_{\dot{\eta}} (\varphi \, (r^* \ga_s)) \\
&= \, \big[ \dot{\eta}(\varphi) - 2 \pi \ti k (\dot{\eta} \, \contrac \, r^* \hat{\theta}_s) \varphi \big] \, (r^* \ga_s) \, = \, 0 \, .
\end{split}
\end{equation*}
\end{proof}
The flat line bundle $r^* \CLH_{\partial P} \ra \CA_P^f$ together with its metric and connection push down to the quotient by $\CG_P$ to a flat hermitian line bundle which is identified with the pullback  $r^* \CL_{\partial P} \ra \CM_P$  of the line bundle $\CL_{\partial P} \ra \CM_{\partial P}$ under the restriction map $r: \CM_P \ra \CM_{\partial P}$.
The $\CG_P$-invariant and covariantly constant section $\hat{\si}_P : \CA_P^f \ra r^* \CLH_{\partial P}$ corresponds to a covariantly constant section 
\begin{equation} \label{e:CSsdown}
\si_P : \CM_P \ra r^* \CL_{\partial P} \, .
\end{equation}
For any two principal $\BT$-bundles $P$ and $P'$ over $X$ with $c_1(P) = c_1(P') =p \in \text{Tors} H^2(X;\BZ)$, let $\phi : P' \ra P$ be a bundle isomorphism.
This induces an isomorphism of hermitian line bundles with connections $\phi^* : r^* \CLH_{\partial P} \ra r^* \CLH_{\partial P'}$ such that for the corresponding Chern-Simons sections (\ref{e:CSsect}) we have $\phi^* \hat{\si}_P(\Th) = \hat{\si}_{P'}(\phi^* \Th)$, for any $\Th \in \CA_P$.
Restricting the line bundles to the subspaces of flat connections and taking the quotient by the group of gauge transformations, we obtain a canonical isomorphism of flat hermitian lines
\begin{equation} \label{e:identif}
\begin{CD}
r^* \CL_{\partial P} @>>> r^* \CL_{\partial P'} \\
@VVV @VVV \\
\CM_P @>>> \CM_{P'}
\end{CD}
\end{equation}
Thus, for each $p \in \text{Tors} H^2(X;\BZ)$, we obtain a hermitian line bundle with a flat connection,
\begin{equation*} 
\CL_p \ra \CM_{X,p} \, ,
\end{equation*}
over the connected component $\CM_{X,p}$ of the moduli space $\CM_X$.
The line bundle $\CL_p \ra \CM_{X,p}$ is canonically identified with the restriction to $\CM_{X,p} \subset \CM_X$ of the pullback of $\CL_{\pX} \ra \CM_{\pX}$ under the map $r_X : \CM_X \ra \CM_{\pX}$.
In view of the canonical isomorphism (\ref{e:identif}), the section (\ref{e:CSsdown}) defines a nowhere-zero covariantly constant section
\begin{equation}
\si_{X,p} \, : \, \CM_{X,p} \lra \CL_p = r^*_X \CL_{\pX} \bigr|_{\CM_{X,p}} 
\end{equation}
which we refer to as the {\em Chern-Simons section}.

\bsk
 

\section{The quantum theory} \label{s:qth}

\subsection{The Hilbert space}

We consider again a closed oriented 2-dimensional manifold $\Si$ and let $g= \frac{1}{2} \dim H^1(\Si;\BR)$.
The space of gauge equivalence classes of flat $\BT$-connections on $\Si$ is a symplectic $2g$-dimensional manifold $(\CM_{\Si}, \om_{\Si})$.
The construction in Sect.\ref{s:line} provides $(\CM_{\Si},k \om_{\Si}), \, k \in 2 \BZ_{+}$, with a hermitian line bundle $\CL_{\Si}$ with a unitary connection $\na$ with curvature $- 2 \pi \ti k \om_{\Si}$.
According to the general geometric quantization scheme the other data needed for the quantization of the symplectic manifold $(\CM_{\Si}, k \om_{\Si})$ is a polarization.

Since $(\CM_{\Si}, \om_{\Si})$ is a symplectic torus, we are going to use the results in \cite{Ma} on the quantization with half-densities of symplectic tori in a real polarization.
Thus let $L \subset H^1(\Si;\BR)$ be a {\em rational} Lagrangian subspace of the symplectic vector space $(H^1(\Si;\BR),\om_{\Si})$.
The attribute rational refers to the fact that the intersection $L \cap H^1(\Si;\BZ)$ with the integer lattice $H^1(\Si;\BZ) \subset H^1(\Si;\BR)$ generates $L$ as a vector space.
Under the identification of the tangent space at any point of $\CM_{\Si}$ with $2 \pi \ti H^1(\Si;\BR)$, the Lagrangian subspace $L$ determines an  invariant real polarization $\CP_L$ of $(\CM_{\Si}, \om_{\Si})$.
The polarization $\CP_L$ is the tangent bundle along the leaves of an invariant Lagrangian foliation of $\CM_{\Si}$.

The bundle of half-densities on $\CP_L$ is the line bundle $\left| \text{Det} \CP_L^* \right|^{\frac{1}{2}}$  over $\CM_{\Si}$.
Its restriction to any leaf $\La$ of $\CP_L$ has a canonical flat connection $\na^{\CP_L}$ defined as follows \cite{Wo}.
First let us recall that, since the quotient map $\pi : \CM_{\Si} \ra \CM_{\Si}/\CP_L$ is a smooth fibration, for each point in $\CM_{\Si}/\CP_L$ there is a local neighborhood $U$ such that $\CP_L \big|_{\pi^{-1}(U)}$ is spanned by Hamiltonian vector fields $[\dot{x}_1], \dots, [\dot{x}_g]$.
Then, on $\pi^{-1}(U)$, the covariant derivative $\na^{\CP_L}_{[\dot{w}]} \mu$ of a section $\mu$ of $\left| \text{Det} \CP_L^* \right|^{\frac{1}{2}}$, along a vector field $[\dot{w}] \in \CP_L \big|_{\pi^{-1}(U)}$, is defined by
\begin{equation*}
\Big( \na^{\CP_L}_{[\dot{w}]} \mu \Big) ([\dot{x}_1], \dots, [\dot{x}_g]) \, = \, [\dot{w}] \big( \mu ([\dot{x}_1], \dots, [\dot{x}_g]) \big) \, .
\end{equation*}

Since the leaves of $\CP_L$ are diffeomorphic to $g$-dimensional tori, the distribution $\CP_L$ has a canonical density $\ka$ invariant under the Hamiltonian vector fields in $\CP_L$ and which assigns to each integral manifold of $\CP_L$ the volume 1.
The square root of $\ka$ trivializes the line bundle $\left| \text{Det} \CP_L^* \right|^{\frac{1}{2}}$ and is a covariantly constant section of $\big(\left| \text{Det} \CP_L^* \right|^{\frac{1}{2}} \big) \big|_{\La}$ for each leaf $\La$ of $\CP_L$.

The Hilbert space of quantization for $(\CM_{\Si},k \om_{\Si})$ will be defined in terms of sections of the line bundle $\CL_{\Si} \otimes \left| \text{Det} \CP_L^* \right|^{\frac{1}{2}}$ obtained by tensoring the prequantum line bundle $\CL_{\Si}$ with the bundle of half-densities on $\CP_L$.
The restriction of this line bundle to a leaf $\La$ of $\CP_L$ has a flat connection defined by
\begin{equation*}
\na_{[\dot{w}]} (\si \otimes \mu) \, = \, \na_{[\dot{w}]} \si \otimes \mu + \si \otimes \na^{\CP_L}_{[\dot{w}]} \mu \, ,
\end{equation*}
for any $[\dot{w}] \in \CP_L \big|_{\La}$ and $\si \otimes \mu \in \Ga(\La; \CL_{\Si} \otimes \left| \text{Det} \CP_L^* \right|^{\frac{1}{2}})$.
The union of those leaves of $\CP_L$ for which the holonomy group of this connection is trivial defines the Bohr-Sommerfeld set $\CB \CS_{\CP_L}$ of the polarization $\CP_L$.
For each leaf $\La$ belonging to $\CB \CS_{\CP_L}$ we let $S_{\La}$ denote the one-dimensional vector space of covariantly constant sections of $\big( \CL_{\Si} \otimes \left| \text{Det} \CP_L^* \right|^{\frac{1}{2}} \big) \bigr|_{\La}$.
Then the Hilbert space associated to the symplectic manifold $(\CM_{\Si},k \om_{\Si})$ with real polarization $\CP_L$ is defined to be the complex vector space
\begin{equation*}
\CH(\Si,L) \, = \, \underset{\La \subset \CB\CS_{\CP_L} }{\oplus} S_{\La}
\end{equation*}
with inner product
\begin{equation*}
\langle \si \otimes \mu , \si' \otimes \mu' \rangle \, = \, 
\begin{cases}
0 \, , & \text{if } \, \si \otimes \mu  \in S_{\La},\, \si' \otimes \mu' \in S_{\La'}, \, \La \neq \La' \\
\int\limits_{\La} (\si,\si') \mu * \mu' \, , & \text{if } \, \si \otimes \mu , \si' \otimes \mu' \in S_{\La}
\end{cases}
\end{equation*}
$(\si,\si')$ is the function on $\La$ obtained by taking the hermitian inner product in the fibre of $\CL_{\Si}$ and $\mu * \mu'$ the density on $\La$ defined by
\begin{equation*}
(\mu * \mu') ([\dot{x}_1], \dots, [\dot{x}_g]) \, = \, \mu([\dot{x}_1], \dots, [\dot{x}_g]) \, \mu'([\dot{x}_1], \dots, [\dot{x}_g]) \, ,
\end{equation*}
for some vector fields $[\dot{x}_1], \dots, [\dot{x}_g]$ spanning $T\La$.

Now recall that we identify $\CM_{\Si} = \frac{2 \pi \ti H^1(\Si;\BR)}{2 \pi \ti H^1(\Si;\BZ)}$.
Let $[\dot{w}_1], \dots, [\dot{w}_g]$ be a basis for $2 \pi \ti (L \cap H^1(\Si;\BZ))$.
Then, as shown in (\cite{Ma}, \S 3), the Bohr-Sommerfeld leaves of $\CP_L$ on $\CM_{\Si}$ are determined by the requirement that their preimages in the linear symplectic space $2 \pi \ti H^1(\Si;\BR)$ covering $\CM_{\Si}$ satisfy the condition
\begin{equation} \label{e:BScondition}
\te^{2 \pi \ti k \om_{\Si}([\dot{w}_i], [x]) } \, = \, 1 \, , \qquad i=1, \dots,g \, .
\end{equation}
For each $ \bq = (q_1,\dots, q_g) \in (\BZ/k \BZ)^g$ the linear equations
\begin{equation} \label{e:BScondition2}
k \,\om_{\Si}([\dot{w}_i], [x]) \, = \, q_i \pmod{k} \, , \qquad i = 1,\dots,g \, ,
\end{equation}
for $[x] \in 2 \pi \ti H^1(\Si;\BR)$, define a family of parallel Lagrangian planes projecting onto a Bohr-Sommerfeld leaf $\La_{\bq}$ under the quotient map $ 2 \pi \ti H^1(\Si;\BR) \ra \CM_{\Si}$.
Thus the dimension of the Hilbert space $\CH(\Si,L)$ is $k^g$.

In conclusion, to each closed oriented surface $\Si$ together with a rational Lagrangian subspace $L \subset H^1(\Si;\BR)$ we associate by the previous construction a finite dimensional Hilbert space $\CH(\Si,L)$.
We recall the following results from \cite{Ma}.

If $L_1 , L_2 \subset H^1(\Si;\BR)$ are two rational Lagrangian subspaces, then there is a canonical unitary operator
\begin{equation} \label{e:intisom}
F_{L_2 L_1} \, : \, \CH(\Si,L_1) \lra \CH(\Si,L_2)
\end{equation}
induced by the Blattner-Kostant-Sternberg (BKS) pairing 
\begin{equation*}
\langle \langle \cdot , \cdot \rangle \rangle : \CH(\Si,L_2) \times \CH(\Si,L_1) \lra \BC
\end{equation*}
The BKS pairing between the Hilbert spaces $\CH(\Si,L_1) =  \underset{\La_1 \subset \mathcal{BS}_{\CP_{L_1}}}{\oplus} \; S_{\La_1}$ and $\CH(\Si,L_2) = \underset{\La_2 \subset \mathcal{BS}_{\CP_{L_2}}}{\oplus} \; S_{\La_2}$ is defined by setting
\begin{equation} \label{e:BKS}
\langle \langle s_2 \otimes \mu_2 , s_1 \otimes \mu_1  \rangle \rangle \, = \,  \int\limits_{\La_1 \cap \La_2} \: (s_2 , s_1) \; \mu_2 \ast \mu_1 \, ,
\end{equation}
for any Bohr-Sommerfeld leaves $\La_1 \subset \mathcal{BS}_{\CP_{L_1}}, \, \La_2 \subset \mathcal{BS}_{\CP_{L_2}}$ and for any sections $ s_1 \otimes \mu_1 \in S_{\La_1}, \, s_2 \otimes \mu_2 \in S_{\La_2}$. 
The density $\mu_2 \ast \mu_1$ on $\La_1 \cap \La_2$ is defined as follows \cite{Sn}.
For any point $[x] \in \La_1 \cap \La_2$, choose a symplectic basis $([\boldsymbol{\dot{v}_2}] , [\boldsymbol{\dot{w}}] ; [\boldsymbol{\dot{v}_1}] , [\boldsymbol{\dot{t}}] )$ for $(\CM_{\Si}, k \om_{\Si})$, that is a basis satisfying 
$k \om_{\Si}( [\dot{w}_i], [\dot{t}_j]) \, = \, \de_{ij}, \,
k \om_{\Si}( [\dot{v}_{2i}], [\dot{v}_{1j}]) \, = \, \de_{ij}$ and
$k \om_{\Si}( [\dot{v}_{1i}] ,  [\dot{t}_j]) \, = \, k \om_{\Si}
( [\dot{v}_{2i}] , [\dot{t}_j])  = 0$,
and with
$[\boldsymbol{\dot{w}}]$ a basis for $ T_{[x]} (\La_1 \cap  \La_2) $,
$( [\boldsymbol{\dot{v}_1}] , [\boldsymbol{\dot{w}}]) $ a basis for 
$T_{[x]} \La_1$ and
$ ( [\boldsymbol{\dot{v}_2}] , [\boldsymbol{\dot{w}}] ) $ a basis for 
$T_{[x]} \La_2$.
Then $\mu_{2} \ast \mu_{1}$ is the density on $\La_1 \cap \La_2$  defined by
\begin{equation*}
(\mu_{2} \ast \mu_{1}) ([\boldsymbol{\dot{w}}]) \, = \, 
\mu_2 ([\boldsymbol{\dot{v}_2}] , [\boldsymbol{\dot{w}}]) \: 
\mu_1 ([\boldsymbol{\dot{v}_1}] , [\boldsymbol{\dot{w}}]) .
\end{equation*}
The operator $F_{L_2 L_1}$ is determined by the relation
\begin{equation} \label{e:BKSisom}
\langle s_2 \otimes \mu_2 , F_{L_2 L_1} (s_1 \otimes \mu_1) \rangle
\, = \, \langle \langle  s_2 \otimes \mu_2 ,   s_1 \otimes \mu_1 \rangle \rangle \,.
\end{equation}
For the unitarity proof of $F_{L_2 L_1}$ we refer to (\cite{Ma},\S 4).

For any three rational Lagrangian subspaces $L_1, L_2, L_3 \subset H^1(\Si;\BR)$ the unitary operators relating the Hilbert spaces associated to $\Si$ and each of these Lagrangian subspaces compose transitively up to a projective factor expressible in terms of the Maslov-Kashiwara index $\tau(L_1,L_2,L_3)$:
\begin{equation} \label{e:trcomposition}
F_{L_1 L_3} \circ F_{L_3 L_2} \circ F_{L_2 L_1} \, = \, \te^{- \frac{\pi \ti}{4} \tau(L_1,L_2,L_3)} \, I
\end{equation}
For the proof of the above relation we refer to (\cite{Ma},\S 6). The definition of the index $\tau$ may be found in \cite{Go, LV, Ma}.

\bsk

\subsection{The vector}

We consider a compact oriented 3-manifold $X$ with nonempty boundary $\pX$.
Through the construction of Sect.\ref{s:line} we have a prequantum line bundle $\CL_{\pX}$ on the symplectic space $(\CM_{\pX},k \om_{\pX})$.
The Lagrangian map $r_X : \CM_X \ra \CM_{\pX}$ determined by the restriction of connections on $X$ to $\pX$ suggests a natural choice for the polarization needed in quantizing  $(\CM_{\pX},k \om_{\pX})$.
That is, we take the invariant real polarization $\CP_X$ on $\CM_{\pX}$ determined by the rational Lagrangian subspace $L_X \subset H^1(\pX;\BR)$, where
$L_X \, = \, \text{Im} \{ \dot{r}_X : H^1(X;\BR) \ra H^1(\pX;\BR) \}$.
To the pair $(\pX, L_X)$ corresponds the Hilbert space
\begin{equation*}
\CH(\pX,L_X) \, = \,\underset{\La \subset \CB\CS_{\CP_X}}{ \oplus} S_{\La}
\end{equation*}
defined as the direct sum of the one-dimensional vector spaces $S_{\La}$ of parallel sections of the line bundle $\CL_{\pX} \otimes \left| \text{Det} \CP_X^* \right|^{\frac{1}{2}}$, supported on the leaves $\La$ contained in the Bohr-Sommerfeld set $\CB\CS_{\CP_X}$ of the polarization $\CP_X$.

Our goal is to give a canonical construction of a vector $Z_X$ in $\CH(\pX,L_X)$ to be associated to the 3-manifold $X$.
First we note the following
\begin{prop}
The Lagrangian submanifold $\La_X = \mathrm{Im} \{ r_X : \CM_X \ra \CM_{\pX} \}$ of $\CM_{\pX}$ is contained in the Bohr-Sommerfeld set $\CB\CS_{\CP_X}$ and is connected.
\end{prop}
\begin{proof}
We start with the observation that for any flat connection $\Th$ on a $\BT$-bundle $P \ra X$ and any closed 1-form $\dot{\al}$ on $X$ for which the cohomology class $\frac{1}{2 \pi \ti} [\dot{\al}]$ in $H^1(X;\BR)$ is integral we have
\begin{equation*}
\int\limits_{\pX} \langle s^* \Th \wedge \dot{\al} \rangle \, = \, 0 \; \pmod{1} \, ,
\end{equation*}
for any section $s:\pX \ra \partial P$.
This follows from (\ref{l:Stokes}), since $F_{\Th} =0$ and since any integral class $[\dot{\al}]$ in $2 \pi \ti H^1(X;\BZ)$ is the cohomology class  $[u^* \vth]$ of the pullback of the Maurer-Cartan form $\vth$ of the group $\BT$ through some map $u:X \ra \BT$.
According to the definition (\ref{e:sympf}) of the symplectic form $\om_{\pX}$ on $\CM_{\pX} = \frac{2 \pi \ti H^1(\pX;\BR)}{2 \pi \ti H^1(\pX;\BZ)}$, we can rewrite the above equation as
\begin{equation} \label{e:integrality}
\om_{\pX}(\dot{r}_X [\dot{\al}], [s^* \Th]) \, = \, 0 \; \pmod{1} \, ,
\end{equation}
with $[s^* \Th] \in 2 \pi \ti H^1(\pX;\BR)$ the cohomology class of the 1-form $s^* \Th$ on $\pX$ and $\dot{r}_X [\dot{\al}]$ the image of $[\dot{\al}]$
under the restriction map $\dot{r}_X : H^1(X;\BR) \ra H^1(\pX;\BR)$.

Let $g = \frac{1}{2} \dim H^1(\pX;\BR)$ and choose a basis $[\dot{w}_1], \dots , [\dot{w}_g]$ for the subspace $2 \pi \ti (L_X \cap H^1(\pX;\BZ))$ of $2 \pi \ti H^1(\pX;\BR)$, with $[\dot{w}_i] = \dot{r}_X [\dot{\al}_i]$ for some integral classes $ [\dot{\al}_i] \in 2 \pi \ti H^1(X;\BZ)$.
Given any $[\eta] \in \La_X$, there exists $[\Th] \in \CM_X$ with $\Th$ a flat connection on a $\BT$-bundle $P \ra X$ such that its restriction to $\pX$ is $\partial \Th = \eta$.
For any section $s : \pX \ra \partial P$, the element $[s^* \eta] \in 2 \pi \ti H^1(\pX;\BR)$ projects to $[\eta]$ under the quotient map 
$ 2 \pi \ti H^1(\pX;\BR) \ra \CM_{\pX} = \frac{2 \pi \ti H^1(\pX;\BR)}{2 \pi \ti H^1(\pX;\BZ)}$.
Then it follows from (\ref{e:integrality}) that
\begin{equation*}
\om_{\pX}( [\dot{w}_i], [s^* \eta]) \, = \, \om_{\pX}( \dot{r}_X [\dot{\al}_i],[s^* \Th]) \, = \, 0 \; \pmod{1} \, , \quad i=1,\dots,g \, .
\end{equation*}
Now, since $\La_X$ is a Lagrangian submanifold of $(\CM_{\pX}, \om_{\pX})$ and since $T \La_X = \CP_X \big|_{\La_X}$, each connected component of $\La_X$ is a leaf of $\CP_X$.
A comparison of the above equation to the Bohr-Sommerfeld condition (\ref{e:BScondition2}) shows that the leaves of the polarization $\CP_X$ contained in $\La_X$ are all Bohr-Sommerfeld leaves.
Moreover, it shows that $\La_X$ contains only one Bohr-Sommerfeld leaf, that is, $\La_X$ is connected.
\end{proof}

The above proposition suggests that we might try to define the vector $Z_X$ in $\CH(\pX,L_X)$ as a parallel section of the line bundle $\CL_{\pX} \otimes \left| \text{Det} \CP_X^* \right|^{\frac{1}{2}}$ restricted to the Bohr-Sommerfeld leaf $\La_X$.
First we note that, since $\La_X$ is a connected Lagrangian submanifold of $\CM_{\pX}$, each connected component $\CM_{X,p}$ of the moduli space $\CM_X = \underset{p \in \text{Tors} H^2(X;\BZ)}{\sqcup} \CM_{X,p}$ maps {\em onto} $\La_X$ under the continuous map $r_X : \CM_X \ra \CM_{\pX} $.
Now recall from Sect.\ref{s:line} that, for each component $\CM_{X,p}$, the restriction $\CL_p = r^*_X \CL_{\pX} \bigr|_{\CM_{X,p}}$ of the pullback under $r_X$ of the prequantum line bundle $\CL_{\pX} \ra \CM_{\pX}$ has a nowhere-zero parallel section, the Chern-Simons section $\si_{X,p}$.
Then, since $\La_X$ is a Bohr-Sommerfeld leaf and since the map $r_X : \CM_{X,p} \ra \La_X$ is surjective, it follows that the section $\si_{X,p}$ is the pullback of a nowhere-zero covariantly constant section of $\CL_{\pX} \bigr|_{\La_X}$ which we continue to denote by $\si_{X,p}$.
We are thus led to define a section $\si_X$ of $\CL_{\pX} \bigr|_{\La_X}$ by
\begin{equation} \label{e:sectpql}
\si_X \, = \, \sum\limits_{p \in \text{Tors} H^2(X;\BZ)} \si_{X,p} \, .
\end{equation}
$\si_X$ is covariantly constant with respect to the connection $\na$ in $\CL_{\pX}$.

What we need now is a covariantly constant global section $\mu_X$ of the bundle of half-densities $\left| \text{Det} \CP_X^* \right|^{\frac{1}{2}}$ over $\La_X$.
We construct such a section using the Reidemeister torsion (R-torsion) invariant $T_X$ of the 3-manifold $X$.

The R-torsion $T_X$ of a compact manifold $X$ is defined as a norm on the determinant line $\left| \text{Det} H^{\bullet}(X;\BR) \right|$ of the cohomology of $X$.
We recall that $\left| \text{Det} H^{\bullet}(X;\BR) \right| = \otimes_q \left| \text{Det} H^q (X;\BR) \right|^{(-1)^{q+1}}$.
The definition and properties of the R-torsion may be found in \cite{RS1,RS2,Mu,V}.
The R-torsion $T_X$ is defined in terms of a smooth triangulation $K$ of $X$.
Let $(C^{\bullet}(K), d^c)$ be the cochain complex of $K$.
There is a preferred basis $\bc^{(q)}=(c^{(q)}_i)$ of $C^q(K)$ defined by 
$c^{(q)}_i(\si^{(q)}_j) =\de_{ij}$, where $(\si^{(q)}_j)$ are the $q$-simplices of $K$.
The space $C^q(K)$ has a natural inner product for which the basic cochains $(c^{(q)}_i)$ are orthonormal.
This inner product passes to the cohomology space $H^q(X;\BR) = H^q(C^{\bullet}(K))$.
Let $\bh^{(q)} = (h^{(q)}_1, \dots, h^{(q)}_{\be_q})$, with $\be_q = \dim H^q(X;\BR)$, be a basis for $H^q(C^{\bullet}(K))$ and let $\tilde{h}^{(q)}_i$ be representative cocycles in $C^q(K)$ for the elements $h^{(q)}_i$.
For each $q$, choose a basis $\bb^{(q)}$ for the image $d^c C^{q-1}(K) \subset C^q(K)$ and let $\tilde{\bb}^{(q)}$ be an independent set in $C^{q-1}(K)$ such that $d^c \tilde{\bb}^{(q)} = \bb^{(q)}$.
Then $(\bb^{(q)}, \tilde{\bb}^{(q+1)}, \tilde{\bh}^{(q)})$ is a basis for $C^q(K)$.
Let $D_q$ denote the matrix representing the change of basis from $(\bb^{(q)}, \tilde{\bb}^{(q+1)}, \tilde{\bh}^{(q)})$ to $\bc^{(q)}$.
The determinant of $D_q$ depends only on the choice of $\bb^{(q)}, \bb^{(q+1)}, \bh^{(q)}$.
The R-torsion of $X$ is the density $T_X \in \left| \text{Det} H^{\bullet}(X;\BR)^* \right|$ defined by
\begin{equation} \label{e:Rtors}
T_X \, = \, \overset{\dim X}{\underset{q=0}{\otimes}}
\left| (\det D_q)^{-1} \, (h^{(q)}_1 \wedge \cdots \wedge h^{(q)}_{\be_q}) \right|^{(-1)^q}
\end{equation}
The above expression is independent of the choice of bases $\bb^{(q)}$ and $\bh^{(q)}$.
The norm (\ref{e:Rtors}) is known to be a combinatorial invariant of $K$; hence any smooth triangulation of the compact manifold $X$ gives the same R-torsion norm on $\left| \text{Det} H^{\bullet}(X;\BR) \right|$.

For a compact oriented 3-manifold $X$ with nonempty boundary the R-torsion $T_X$ belongs to
$ \left| \text{Det} H^0(X;\BR) \right| \otimes
\left| \text{Det} H^1(X;\BR)^* \right| \otimes
\left| \text{Det} H^2(X;\BR) \right| \otimes
\left| \text{Det} H^3(X;\BR)^* \right| $.
By Poincare duality $ H^2(X;\BR) \cong H^1(X, \pX;\BR)^*$. 
Choosing canonical trivializations 
$\left| \text{Det} H^0(X;\BR) \right| \cong \BR_{+}$ 
and $\left| \text{Det} H^3(X;\BR) \right| \cong \BR_{+}$ with respect to orthonormal bases for $H^0(X;\BR)$ and $H^3(X;\BR)$, we can then regard $T_X$ as a density
\begin{equation} \label{e:Rt}
T_X \, \in \, \left| \text{Det} H^1(X;\BR)^* \right| \otimes \left| \text{Det} H^1(X, \pX;\BR)^* \right|
\end{equation}
We recall from Sect.\ref{s:moduli} that $\CM_X = H^1(X;\BT)$ and that we have an exact sequence (\ref{e:exfibr}) of compact abelian Lie groups.
The corresponding sequence between the tangent spaces at the identity of the groups in (\ref{e:exfibr}),
\begin{equation*} 
\begin{split}
0 \ra H^0(X,\pX;\BR) &\ra H^0(X;\BR) \ra H^0(\pX;\BR) \ra \\
& \ra H^1(X,\pX;\BR) \ra H^1(X;\BR) \ra T_e\La_X \ra 0 \, ,
\end{split}
\end{equation*}
induces the following isomorphism between the half-densities spaces:
\begin{equation} \label{e:tangLa}
\left| \text{Det} T_e^* \La_X \right|^{\frac{1}{2}} \cong 
\left| \text{Det} H^1(X;\BR)^* \right|^{\frac{1}{2}} \otimes  
\left| \text{Det} H^1(X,\pX;\BR) \right|^{\frac{1}{2}} \, ,
\end{equation}
where we used again canonical trivializations for the spaces $H^0(*;\BR)$.
Let us pick an arbitrary element $w \in \left| \text{Det} H^1(X, \pX;\BR)^* \right|$.
Since $H^1(X,\pX;\BR)$ is identified with the tangent space at the identity of the group $H^1(X,\pX;\BT)$, $w$ extends to an invariant density $\bw$ on $H^1(X,\pX;\BT)$.
If $w^{-1} \in \left| \text{Det} H^1(X, \pX;\BR) \right|$ stands for the dual of $w$, we note that $(T_X)^{\frac{1}{2}} \otimes w^{-1} \in \left| \text{Det} H^1(X;\BR)^* \right|^{\frac{1}{2}} \otimes  
\left| \text{Det} H^1(X,\pX;\BR) \right|^{\frac{1}{2}}$.
Thus we can define an invariant half-density $\mu_X$ on the abelian group $\La_X$, whose value at the identity is
\begin{equation} \label{e:hdensid}
(\mu_X)_e \,= \, \Bigg[ \int\limits_{H^1(X,\pX;\BT)} \bw \Bigg] \, 
(T_X)^{\frac{1}{2}} \otimes w^{-1}
\end{equation}
The identification between the l.h.s. and r.h.s. of the above equation is made via the isomorphism (\ref{e:tangLa}).
As $T \La_X = \CP_X \bigr|_{\La_X}$, the half-density $\mu_X$, being $\La_X$-invariant, defines therefore a covariantly constant section of the half-densities bundle $\big( \left| \text{Det} \CP_X^* \right|^{\frac{1}{2}}   \big) \bigr|_{\La_X}$.
We introduce the notation
\begin{equation} \label{e:hdens}
\mu_X \, = \, \int\limits_{H^1(X,\pX;\BT)} (T_X)^{\frac{1}{2}}
\end{equation}
with the implicit understanding of $\mu_X$ as the invariant extension of (\ref{e:hdensid}).

Using (\ref{e:sectpql}) and (\ref{e:hdens}) we associate then to the compact oriented 3-manifold $X$ with boundary $\pX$ the vector $Z_X$ in the Hilbert space $\CH(\pX,L_X)$, defined by the expression
\begin{equation} \label{e:Zb}
\begin{split}
Z_X \, &= \, \frac{k^{m_X}}{[ \# \text{Tors} H^2(X;\BZ) ]} \; \si_X \otimes \mu_X \\
&= \, \frac{k^{m_X}}{[ \# \text{Tors} H^2(X;\BZ) ]} \, \sum\limits_{p \in \text{Tors} H^2(X;\BZ)} \si_{X,p} \otimes \int\limits_{H^1(X,\pX;\BT)} (T_X)^{\frac{1}{2}}
\end{split}
\end{equation}
where
\begin{align} \label{e:mx}
m_X \, = \, \frac{1}{4} \, & \big( \dim H^1(X;\BR) + \dim H^1(X,\pX;\BR) \\
                     &- \dim H^0(X;\BR) - \dim H^0(X,\pX;\BR) \big) \,.\notag
\end{align}

If $X$ is a closed oriented 3-manifold then, according to (\ref{e:Rt}), the square root of the R-torsion $T_X$ is a density $(T_X)^{\frac{1}{2}} \in \left| \text{Det} H^1(X;\BR)^* \right|$.
Since the tangent space $T\CM_X \cong H^1(X;\BR)$, this defines an invariant density $(T_X)^{\frac{1}{2}}$ on the moduli space $\CM_X = H^1(X;\BT)$.
From (\ref{t:CSsect} (a)) it follows that we have a well-defined function $\si_X$ on $\CM_X$ whose restriction to each connected component $\CM_{X,p}$ of $\CM_X$ is given by $\si_{X,p}([\Th]) = \te^{\pi \ti k S_{X,P}(\Th)} $, for any $\BT$-bundle $P$ with $c_1(P)=p$.
As a result of (\ref{p:statnb}) $\si_X$ is constant on  components.
Thus to the closed oriented 3-manifold $X$ we associate the complex number
\begin{equation} \label{e:Znb}
\begin{split} 
Z_X \, &= \, \frac{k^{m_X}}{[ \# \text{Tors} H^2(X;\BZ) ]} \, \sum\limits_{p \in \text{Tors} H^2(X;\BZ)} \si_{X,p} \; \int\limits_{\CM_X} (T_X)^{\frac{1}{2}} \\
&= \, k^{m_X}  \sum\limits_{p \in \text{Tors} H^2(X;\BR)} \, \int\limits_{\CM_{X,p}} \si_{X,p} \, (T_X)^{\frac{1}{2}} \\
&= \, k^{m_X} \int\limits_{\CM_X} \si_X \, (T_X)^{\frac{1}{2}} 
\end{split}
\end{equation}
where
\begin{equation} \label{e:mxnb}
m_X \, = \, \frac{1}{2} \, \big( \dim H^1(X;\BR) - \dim H^0(X;\BR) \big) \,.
\end{equation}

\bsk


\section{The topological quantum field theory} \label{s:TQFT}

The results of the previous section show that we are very close to defining a 2+1 dimensional topological quantum field theory (TQFT) for the Lie group $\BT$ and an even integer $k$.
At this stage we know how to canonically associate to a closed oriented 2-manifold $\Si$ together with a rational Lagrangian subspace $L \subset H^1(\Si;\BR)$ a finite dimensional Hilbert space $\CH(\Si,L)$ and to a compact oriented 3-manifold $X$ a vector $Z_X$ in the Hilbert space $\CH(\pX,L_X)$ associated to $\pX$ and the Lagrangian subspace $L_X \subset H^1(\pX;\BR)$  determined by $X$.

A TQFT is supposed to satisfy the axioms listed in \cite{A1} which refer to functoriality, the orientation and disjoint union properties and the gluing of manifolds along their common boundaries.

In order to construct the Chern-Simons TQFT for the group $\BT$ we are going to use the notions introduced in \cite{Wa} of extended 2- and 3-manifolds, extended morphisms and gluing of extended 3-manifolds.
We give below the relevant definitions. As in \cite{Wa} we use the shortened notation with the prefix 'e-' standing for 'extended'.

\begin{defi} \label{d:e2}
An {\em e-2-manifold\/} is a pair $(\Si,L)$, with $\Si$ a closed oriented 2-dimensional manifold and $L$ a rational Lagrangian subspace of $H^1(\Si;\BR)$.
\end{defi}

\begin{defi} \label{d:e3}
(i) An {\em e-3-manifold\/} is a triple $(X,L,n)$, with $X$ a compact oriented 3-manifold, $L$ a rational Lagrangian subspace of $H^1(\pX;\BR)$ and $n \in \BZ/ 8\BZ$. \\
(ii) The boundary of an e-3-manifold is $\partial (X,L,n) = (\pX,L)$. \\
A closed e-3-manifold is just a pair $(X,n)$ with $X$ a closed oriented 3-manifold and $n \in \BZ/ 8\BZ$.
\end{defi}

\begin{defi} \label{d:emor}
An  extended morphism from an e-2-manifold $(\Si',L')$ to an e-2-manifold $(\Si,L)$ is  a pair $(h,m)$ with $h: \Si' \ra \Si$ an orientation preserving diffeomorphism and $m \in \BZ/ 8\BZ$.
We refer to extended morphisms between e-2-manifolds as {\em e-2-morphisms}.
\end{defi}

\begin{defi} \label{d:cemor2}
{\em Composition of e-2-morphisms.}\\
If $(h,m) : (\Si',L') \ra (\Si,L)$ and $(h',m') : (\Si'', L'') \ra (\Si',L')$ are e-2-morphisms, their composition is defined to be
\begin{equation*}
(h,m) (h',m') = \big( h h',\, m+m'+\tau(L'',h^{\prime *} L', (h h')^* L) \pmod{8} \big) \, ,
\end{equation*}
where $\tau$ is the Maslov-Kashiwara index \cite{Go,LV} of a triple of Lagrangian subspaces of a symplectic vector space.
\end{defi}

\nin The e-2-morphism $(id,0)$ acts as the identity and the inverse of an e-2-morphism $(h,m)$ is $(h,m)^{-1} = (h^{-1},-m)$.

\begin{defi} \label{d:e3mor}
An extended morphism from an e-3-manifold $(X,L,n)$ to an e-3-manifold $(X',L',n')$ is a pair $(\Phi, m)$ with  $\Phi : X' \ra X$ an orientation preserving diffeomorphism  and $m \in \BZ/ 8\BZ$ such that 
\begin{equation*}
n' = n+m+\tau(L_{X'},L', (\partial \Phi)^* L) \pmod{8}
\end{equation*}
where $\partial \Phi : \pX' \ra \pX$ is the induced diffeomorphism between the boundaries. 
We refer to extended morphisms between e-3-manifolds as {\em e-3-morphisms}.\\
To the e-3-morphism $(\Phi,m) : (X',L',n') \ra (X,L,n)$ there corresponds an induced e-2-morphism between the boundaries $(\partial \Phi,m) : (\pX',L') \ra (\pX,L)$. \\
Two closed e-3-manifolds $(X,n)$ and $(X',n')$ are isomorphic if there exists an orientation preserving diffeomorphism $\Phi : X' \ra X$ and $n'=n$.
\end{defi}

\begin{defi} \label{d:cemor3}
{\em Composition of e-3-morphisms.} \\
If $(\Phi,m) :(X',L',n') \ra (X,L,n)$  and $(\Phi',m') :(X'',L'',n'') \ra (X',L',n')$ are e-3-morphisms, then their composition is the e-3-morphism
\begin{equation*}
(\Phi,m) (\Phi',m') = \big( \Phi \Phi', \, m+m'+\tau(L'', (\partial \Phi')^* L', (\partial \Phi \partial \Phi')^* L) \pmod{8} \big)
\end{equation*}
From the cocycle and the symplectic invariance properties \cite{LV,Go} of the Maslov-Kashiwara index it follows that the equality 
\begin{equation*}
\begin{split}
n'' = n+m+m' &+ \tau(L'', (\partial \Phi')^* L', (\partial \Phi \partial \Phi')^* L) + \\
&+\tau(L_{X''}, L'', (\partial \Phi \partial \Phi')^* L) \pmod{8}
\end{split}
\end{equation*}
is indeed satisfied.
\end{defi}

\begin{defi} \label{d:edisjointu}
(i) If $(\Si_1,L_1)$ and $(\Si_2,L_2)$ are e-2-manifolds, their disjoint union is the e-2-manifold
\begin{equation*}
(\Si_1,L_1) \sqcup (\Si_2,L_2) \, = \, (\Si_1 \sqcup \Si_2, L_1 \oplus L_2) \, ,
\end{equation*}
where we have in view the identification $H^1(\Si_1 \sqcup \Si_2;\BR) = H^1(\Si_1;\BR) \oplus H^1(\Si_2;\BR)$. \\
(ii) If $(h_1,m_1)$ and $(h_2,m_2)$ are e-2-morphisms, then
\begin{equation*}
(h_1,m_1) \sqcup (h_2,m_2) \, = \, (h_1 \sqcup h_2, m_1+m_2)
\end{equation*}
(iii) If $(X_1,L_1,n_1)$ and $(X_2,L_2,n_2)$ are e-3-manifolds, their disjoint union is the e-3-manifold
\begin{equation*}
(X_1,L_1,n_1) \sqcup (X_2,L_2,n_2) \, = \, (X_1 \sqcup X_2, L_1 \oplus L_2, n_1 +n_2)
\end{equation*}
\end{defi}

\begin{defi} \label{d:gl3e}
{\em Gluing e-3-manifolds.} \\
Let $(\tilde{X}, \tilde{L}, \tilde{n})$ be an e-3-manifold such that $\partial (\tilde{X}, \tilde{L}, \tilde{n}) = (\partial \tilde{X}, \tilde{L}) = (Y,L) \sqcup (- \Si_1,L_1) \sqcup (\Si_2,L_2)$
and assume there exists an e-2-morphism $(h,m) : (\Si_1,L_1) \ra (\Si_2,L_2)$.
Under these assumptions, the gluing of $(\tilde{X}, \tilde{L}, \tilde{n})$ by $(h,m)$ is defined to be the e-3-manifold
\begin{equation*}
(\tilde{X}, \tilde{L}, \tilde{n})_{(h,m)} \, = \, \big( X , L, \tilde{n} +m+ \tau(\tilde{L},L_{\tilde{X}}, L_X \oplus L_h) \pmod{8} \big) \, ,
\end{equation*}
where $X$ is the 3-manifold obtained by gluing $\tilde{X}$ by the diffeomorphism $h : \Si_1 \ra \Si_2$ and where
$\tilde{L}, L_{\tilde{X}}$ and $L_X \oplus L_h$ are the Lagrangian subspaces of 
$H^1(\partial \tilde{X};\BR) = H^1(Y;\BR) \oplus H^1(-\Si_1;\BR) \oplus H^1(\Si_2;\BR)$ defined as follows: \\
$\bullet \; \tilde{L} \, = \, L \oplus L_1 \oplus L_2$ by the initial assumption on $(\tilde{X}, \tilde{L}, \tilde{n})$; \\
$\bullet \; L_{\tilde{X}} \, = \, \text{Im} \, \{ H^1(\tilde{X};\BR) \lra H^1(\partial \tilde{X};\BR) \} $; \\
$\bullet \; L_X \oplus L_h$ is the direct sum of  $L_X = \text{Im} \{ H^1(X;\BR) \ra H^1(Y;\BR) \} \subset H^1(Y;\BR)$ and of
the Lagrangian subspace $L_h \subset  H^1(-\Si_1;\BR) \oplus H^1(\Si_2;\BR)$ defined as the graph of the symplectic isomorphism $h^* : H^1(\Si_2;\BR) \ra H^1(\Si_1;\BR)$, that is, \\
$L_h \, = \, \{ (h^* [\dot{\eta}], [\dot{\eta}] ) \, \mid \, [\dot{\eta}] \in  H^1(\Si_2;\BR) \} $.
\end{defi}

Having introduced the necessary definitions, we are ready now to construct the full TQFT.
For each e-2-manifold $(\Si,L)$ we have a finite dimensional Hilbert space $\CH(\Si,L)$ by the construction of Sect.\ref{s:qth}.
Now let $(X,L,n)$ be an e-3-manifold.
Then to the e-2-manifold $(\pX,L)$ obtained as the boundary of $(X,L,n)$ there corresponds the Hilbert space $\CH(\pX,L)$. 
On the other hand, the 3-manifold $X$ determines the Lagrangian subspace $L_X \subset H^1(\pX;\BR)$ and therefore an e-2-manifold $(\pX,L_X)$ with corresponding Hilbert space $\CH(\pX,L_X)$.
Let
\begin{equation}
F_{L L_X} : \CH(\pX,L_X) \lra \CH(\pX,L)
\end{equation}
be the unitary isomorphism induced by the BKS pairing (\ref{e:BKS}).
Then we define the vector $Z_{(X,L,n)}$ in $\CH(\pX,L)$ associated to the e-3-manifold $(X,L,n)$ by 
\begin{equation} \label{e:eZb}
Z_{(X,L,n)} \, = \, \te^{\frac{\pi \ti}{4} n} F_{L L_X} (Z_X) \, ,
\end{equation}
where $Z_X$ is the standard vector in $\CH(\pX,L_X) $ defined by the expression (\ref{e:Zb}).

The following theorem shows that we have a unitary TQFT.

\begin{thm} \label{t:TQFT}
The assignments
\begin{alignat*}{2}
\text{e-2-manifold} & & (\Si,L)  & \; \longmapsto \; Hilbert \: space \; \CH(\Si,L) \\ 
\text{e-3-manifold} & & \; (X,L,n)  & \; \longmapsto \; vector \; Z_{(X,L,n)} \in \CH(\pX,L)
\end{alignat*}
satisfy: \\
(a) Functoriality \\
To each e-2-morphism $(h,m) : (\Si',L') \ra (\Si,L)$ there corresponds a unitary isomorphism
\begin{equation*}
U(h,m) : \CH(\Si,L) \lra \CH(\Si',L') 
\end{equation*}
and these compose properly. \\
Let $(\Phi,m) : (X',L',n') \lra (X,L,n)$ be an e-3-morphism between isomorphic e-3-manifolds and $(\partial \Phi,m) : (\pX',L') \lra (\pX,L)$ the induced e-2-morphism between the boundaries. Then
\begin{equation*}
U(\partial \Phi,m) \, Z_{(X,L,n)}  \, = \, Z_{(X',L',n')} 
\end{equation*}
(b) Orientation \\
There is a natural isomorphism of Hilbert spaces
\begin{equation*}
\CH(-\Si,L) \cong \overline{\CH(\Si,L) }
\end{equation*}
and
\begin{equation*}
Z_{(-X,L,n)} = \overline{Z_{(X,L,n)}}
\end{equation*}
(c) Disjoint union \\
If $(\Si,L) = (\Si_1,L_1) \sqcup (\Si_2,L_2) $ is a disjoint union of e-2-manifolds, then there is a natural unitary isomorphism
\begin{equation*}
\CH(\Si_1 \sqcup \Si_2, L_1 \oplus L_2) \cong \CH(\Si_1,L_1) \otimes \CH(\Si_2,L_2)
\end{equation*}
If $(X,L,n) = (X_1,L_1,n_1) \sqcup (X_2,L_2,n_2)$ is a disjoint union of e-3-manifolds, then
\begin{equation*}
Z_{(X_1 \sqcup X_2,L_1 \oplus L_2 ,n_1 +n_2)} \, = \, Z_{(X_1,L_1,n_1)} \otimes Z_{(X_2,L_2,n_2)}
\end{equation*}
(d) Cylinder axiom \\
If $\Si$ is a closed oriented 2-manifold, $I=[0,1]$ the unit interval and $L_{\Si} \subset H^1(\Si;\BR)$ a rational Lagrangian subspace, then to the e-3-manifold $(\Si \times I, L_{\Si} \oplus L_{\Si}, 0)$ there corresponds
\begin{equation*}
Z_{(\Si \times I, L_{\Si} \oplus L_{\Si}, 0)} \, = \, Id \, : \, \CH(\Si,L_{\Si}) \lra  \CH(\Si,L_{\Si}) 
\end{equation*}
(e) Gluing \\
Let $X$ be a compact connected oriented 3-manifold with boundary $\pX$ and let $\XC$  denote the manifold obtained by cutting $X$ along a codimension one closed oriented submanifold $\Si$.
Then $\pXC = \pX \sqcup (-\Si) \sqcup \Si$. 
Let $Z_{(X,L,n)}$ be the vector associated to the e-3-manifold $(X,L,n)$ and $Z_{ (\XC,L \oplus L_{\Si} \oplus L_{\Si} ,n^{cut} ) }$ the vector for the cut e-3-manifold $(\XC,L \oplus L_{\Si} \oplus L_{\Si} ,n^{cut} )$, with $L_{\Si} \subset H^1(\Si;\BR)$ an arbitrary rational Lagrangian subspace and 
\begin{equation*}
n^{cut} = n - \tau(L \oplus L_{\Si} \oplus L_{\Si}, L_{\XC} , L_X \oplus L_{\De}) \pmod{8} \,.
\end{equation*}
The Lagrangian subspace $L_{\De} \subset H^1(-\Si;\BR) \oplus H^1(\Si;\BR) $ is the diagonal. 
Then to the gluing of e-3-manifolds
\begin{equation*}
(\XC, L \oplus L_{\Si} \oplus L_{\Si} ,n^{cut} ) _{(id_{\Si},0)} \, = \, (X,L,n)
\end{equation*}
there corresponds the quantum gluing property
\begin{equation*}
Z_{(X,L,n)} \, = \, \mathrm{Tr}_{\Si} \Big[  Z_{ (\XC,L \oplus L_{\Si} \oplus L_{\Si} ,n^{cut} ) } \Big]
\end{equation*}
where the operator $\mathrm{Tr}_{\Si}$ is the contraction
\begin{equation*}
\mathrm{Tr}_{\Si} : \CH(\pXC, L \oplus L_{\Si} \oplus L_{\Si}) \cong 
\CH(\pX,L) \otimes \overline{\CH(\Si,L_{\Si})} \otimes \CH(\Si,L_{\Si}) \ra \CH(\pX,L)
\end{equation*}
using the hermitian inner product in $\CH(\Si,L_{\Si})$.
\end{thm}

\begin{proof}
(a) {\em Functoriality}.
The orientation preserving diffeomorphism $h: \Si' \ra \Si$ induces a map in cohomology $h^* : H^1(\Si;*) \ra H^1(\Si';*)$ and therefore a symplectic diffeomorphism $h^* : \CM_{\Si} \ra \CM_{\Si'}$ between the moduli spaces of flat $\BT$-connections.
Since there exist lifts of $h$ to bundle morphisms between trivializable $\BT$-bundles over $\Si'$ and trivializable $\BT$-bundles over $\Si$, it follows from (\ref{t:pline} (a)) and (\ref{e:equivlines}) that there is an induced isomorphism $h^* : \CL_{\Si} \ra \CL_{\Si'}$ of prequantum line bundles covering $h^* : \CM_{\Si} \ra \CM_{\Si'}$.
The Lagrangian subspace $L \subset H^1(\Si;\BR)$ is mapped onto the Lagrangian subspace $h^* L \subset H^1(\Si';\BR)$.
The map of line bundles
\begin{equation*}
\begin{CD}
\CL_{\Si} \otimes \left| \text{Det} \CP_L^* \right|^{\frac{1}{2}} @>{h^*}>>
\CL_{\Si'} \otimes \left| \text{Det} \CP_{h^* L}^* \right|^{\frac{1}{2}} \\
@VVV @VVV \\
\CM_{\Si} @>{h^*}>> \CM_{\Si'}
\end{CD}
\end{equation*}
induces a map of sections which, having in view the constructions of Sect.\ref{s:qth}, gives rise to a unitary map of Hilbert spaces
\begin{equation*}
h^* \, : \, \CH(\Si,L) \lra \CH(\Si',h^* L) \, .
\end{equation*}
The composition of this map with the isomorphism 
\begin{equation*}
F_{L',h^* L} \, : \, \CH(\Si',h^* L) \lra \CH(\Si',L')
\end{equation*}
arising from the BKS pairing between $\CH(\Si',L')$ and $\CH(\Si',h^* L)$ gives the unitary operator
$F_{L',h^* L} \circ h^* \, : \, \CH(\Si,L) \lra \CH(\Si',L')$.
Thus we are led to associate to the e-2-morphism $(h,m)$ the unitary operator
\begin{equation*}
U(h,m) \, = \, \te^{\frac{\pi \ti}{4} m} \, F_{L',h^* L} \circ h^* \, : \,
\CH(\Si,L) \lra \CH(\Si',L')
\end{equation*}
We note that if $L$ and $K$ are rational Lagrangian subspaces in $H^1(\Si;\BR)$, then it follows from the definition of the intertwining isomorphisms (\ref{e:intisom}) that the following diagram commutes
\begin{equation*}
\begin{CD}
\CH(\Si,L) @>{F_{K,L}}>> \CH(\Si,K) \\
@V{h^*}VV @VV{h^*}V \\
\CH(\Si',h^* L) @>{F_{h^* K,h^* L}}>> \CH(\Si',h^* K)
\end{CD}
\end{equation*}
Now let $(h',m') : (\Si'',L'') \ra (\Si',L')$ be another e-2-morphism.
Then, using  the composition property (\ref{e:trcomposition}), the above commutative diagram property and the multiplication law of e-2-morphisms, we obtain
\begin{equation*}
\begin{split}
U(h',m') \, U(h,m) \, &= \, 
\te^{\frac{\pi \ti}{4} (m+m')} \, F_{L'',h^{\prime *} L'} \circ h^{\prime *} \circ F_{L',h^* L} \circ h^* \\
&= \, \te^{\frac{\pi \ti}{4} (m+m')} \,  F_{L'',h^{\prime *} L'}  \circ F_{h^{\prime *} L',h^{\prime *} h^* L} \circ h^{'*} \circ h^* \\
&= \, \te^{\frac{\pi \ti}{4} [m+m' + \tau(L'',h^{\prime *}L', (h h')^* L)]} \,
F_{L'', (h h')^* L} \circ (h h')^* \\
&= \, U( (h,m) (h',m') )
\end{split}
\end{equation*}
This proves the first statement in (\ref{t:TQFT} (a)).

We prove now the second statement in (\ref{t:TQFT} (a)).
The orientation preserving diffeomorphism of compact oriented 3-manifolds $\Phi : X' \ra X$ induces an isomorphism $\Phi^* : \CM_X \ra \CM_{X'}$ between the moduli spaces of flat $\BT$-connections.
Its restriction $(\partial \Phi)^* : \CM_{\pX} \ra \CM_{\pX'}$ lifts to an isomorphism of prequantum line bundles $(\partial \Phi)^* : \CL_{\pX} \ra \CL_{\pX'}$.
Moreover we have $r_{X'} \circ \Phi^* = (\partial \Phi)^* \circ r_X$, where $r_X : \CM_X \ra \CM_{\pX}$ and $r_{X'} : \CM_{X'} \ra \CM_{\pX'}$ are the restriction maps to the boundaries.
Therefore $(\partial \Phi)^* L_X = L_{X'}$, so there is an induced bundle isomorphism $(\partial \Phi)^* : \CL_{\pX} \otimes \left| \text{Det} \CP_X^* \right|^{\frac{1}{2}}   \ra \CL_{\pX'} \otimes \left| \text{Det} \CP_{X'}^* \right|^{\frac{1}{2}}$.
From the construction in Sect.\ref{s:qth} of the standard vectors $Z_X$ and $Z_{X'}$ as sections of these line bundles it follows that 
\begin{equation*}
(\partial \Phi)^* (Z_X) \, = \, Z_{X'}
\end{equation*}
under the unitary map of Hilbert spaces $(\partial \Phi)^* : \CH(\pX,L_X) \ra \CH(\pX',L_{X'})$.
On the other hand to the e-2-morphism $(\partial \Phi,m) : (\pX',L') \ra (\pX,L)$ there corresponds the unitary operator
\begin{equation*}
U(\partial \Phi,m) \, = \, \te^{\frac{\pi \ti}{4} m} \, F_{L', (\partial \Phi)^* L} \circ (\partial \Phi)^* \, : \, \CH(\pX,L) \ra \CH(\pX',L')
\end{equation*}
According to the definition (\ref{e:Zb}) of vectors associated to e-3-manifolds we have
$Z_{(X,L,n)} = \te^{\frac{\pi \ti}{4} n} \, F_{L L_X}(Z_X)$ and 
$Z_{(X',L',n')} = \te^{\frac{\pi \ti}{4} n'} \, F_{L' L_{X'}}(Z_{X'})$.
Thus we find
\begin{equation*}
\begin{split}
U(\partial \Phi,m) Z_{(X,L,n)} \, &= \, 
\te^{\frac{\pi \ti}{4} (n+m)} \: F_{L', (\partial \Phi)^* L} \circ
 (\partial \Phi)^* \circ F_{L L_X} (Z_X) \\
&= \, \te^{\frac{\pi \ti}{4} (n+m)} \:
 F_{L', (\partial \Phi)^* L} \circ F_{ (\partial \Phi)^*L,  (\partial \Phi)^* L_X} \circ
 (\partial \Phi)^* (Z_X) \\
&= \, \te^{\frac{\pi \ti}{4} [n+m + \tau((\partial \Phi)^*L_X,L',(\partial \Phi)^*L)]} \:  F_{L', (\partial \Phi)^* L_X} (Z_{X'}) \\
&= \, \te^{\frac{\pi \ti}{4} [n+m + \tau(L_{X'},L',(\partial \Phi)^*L)]} \:
 F_{L' L_{X'}} (Z_{X'}) \\
&= \, \te^{\frac{\pi \ti}{4} n'} \, F_{L' L_{X'}} (Z_{X'}) \, = \, Z_{(X',L',n')}
\end{split}
\end{equation*}
This concludes the proof of (\ref{t:TQFT}(a)).

\nin (b) {\em Orientation}. Follows from (\ref{t:pline} (b)) and (\ref{t:CSsect} (b)).

\nin (c) {\em Disjoint union}.
To the disjoint union $\Si = \Si_1 \sqcup \Si_2$ there corresponds the product of symplectic spaces
\begin{equation*}
(\CM_{\Si}, \om_{\Si}) = (\CM_{\Si_1}, \om_{\Si_1}) \times (\CM_{\Si_2}, \om_{\Si_2}) = (\CM_{\Si_1} \times \CM_{\Si_2}, pr^*_{\Si_1} \om_{\Si_1} + pr^*_{\Si_2} \om_{\Si_2})
\end{equation*}
where $pr_{\Si_i} : \Si_1 \times \Si_2 \ra \Si_i$ are the natural projections.
If $L_i \subset H^1(\Si_i;\BR)$ are rational Lagrangian subspaces, then $L_1 \oplus L_2 \subset H^1(\Si_1 \sqcup \Si_2;\BR) = H^1(\Si_1;\BR) \oplus H^1(\Si_2;\BR)$ is again a rational Lagrangian subspace.
The Bohr-Sommerfeld orbits $\La$ of the polarization $P_{L_1 \oplus L_2}$ on $\CM_{\Si_1} \times \CM_{\Si_2}$ are Cartesian products $\La_1 \times \La_2$ of Bohr-Sommerfeld orbits $\La_i$ for the polarizations $\CP_{L_i}$ on $\CM_{\Si_i},\, i=1,2$.
The proof of (c) follows then from (\ref{t:pline} (c)) and (\ref{t:CSsect} (c)) and the definitions of the  Hilbert spaces associated to e-2-manifolds and of the vectors associated to e-3-manifolds.

\nin (d) {\em Cylinder axiom}.
Since $H^i(\Si \times I) \cong H^i(\Si)$, $i=0,1,2$, we obtain in particular that  $\CM_{\Si \times I} \cong \CM_{\Si}$.
As $\partial (\Si \times I) = (-\Si) \sqcup \Si$ we have $\CM_{\partial (\Si \times I)} = \CM_{\Si} \times \CM_{\Si}$ and the restriction map 
\begin{equation*}
r_{\Si \times I} : \CM_{\Si \times I} \cong \CM_{\Si} \lra \CM_{\partial (\Si \times I)} = \CM_{\Si} \times \CM_{\Si}
\end{equation*}
is the diagonal map $[\eta] \mapsto ([\eta],[\eta])$.
Hence the Lagrangian subspace $L_{\Si \times I} = L_{\De}$, with $L_{\De}$ denoting the diagonal in $H^1(-\Si;\BR) \oplus H^1(\Si;\BR)$ and the Lagrangian submanifold $\La_{\Si \times I} = \De$, with $\De$ the diagonal in $\CM_{\Si} \times \CM_{\Si}$.
According to the definitions and constructions of Sect.\ref{s:qth} we find that the Chern-Simons section over $\De$ of the prequantum line bundle is $\si_{\Si \times I} = 1$ and that the section $\mu_{\Si \times I}$ over $\De$ of the half-density bundle ${\left| \text{Det} \,{\mathcal P}^*_{\De} \right|}^{\frac{1}{2}} $  is identified  under the isomorphism $\De \cong \CM_{\Si}$ with the square root $(T_{\Si})^{\frac{1}{2}}$ of the R-torsion of $\Si$.
In \cite{Wi2} it is proved that the density $T_{\Si}$ on $\CM_{\Si}$ coincides with the density $\bigl| \frac{(\om_{\Si})^{g} }{g !} \bigr|$ determined by the symplectic form $\om_{\Si}$, where $g = \frac{1}{2} \dim H^1(\Si;\BR)$.
Having in view the definition (\ref{e:Zb}), we find therefore that the standard vector in $\CH(-\Si \sqcup \Si, L_{\De})$ associated to the 3-manifold $\Si \times I$ is
\begin{equation} \label{e:cylinder}
Z_{\Si \times I} \, = \, k^{\frac{1}{4} \dim H^1(\Si;\BR)} \, \si_{\Si \times I} \otimes \mu_{\Si \times I} 
\, = \, k^{\frac{1}{4} \dim H^1(\Si;\BR)} \, 1 \otimes (T_{\Si})^{\frac{1}{2}}\, .
\end{equation}
Let us consider now the e-3-manifold $(\Si \times I , L_{2 \Si},n)$, where $L_{2 \Si} = L_{\Si} \oplus L_{\Si}$ with $L_{\Si} \subset H^1(\Si;\BR)$ an arbitrary rational Lagrangian subspace.
To this e-3-manifold there corresponds the vector in $\CH(-\Si \sqcup \Si,L_{2\Si}) \cong \overline{\CH(\Si,L_{\Si})} \otimes \CH(\Si,L_{\Si})$ defined by 
\begin{equation*}
Z_{(\Si \times I , L_{2 \Si},n)} \, = \, \te^{\frac{\pi \ti}{4}   n} \, F_{L_{2\Si},L_{\De}} (Z_{\Si \times I}) \, .
\end{equation*}
Let us introduce a unitary basis $\{ v{\sst (\Si,L_{\Si})}_{\bq} \}$ for the Hilbert space $\CH(\Si,L_{\Si})$. Using the definition of the unitary operator 
$ F_{L_{2\Si},L_{\De}} : \CH(-\Si \sqcup \Si,L_{\De}) \lra  \CH(-\Si \sqcup \Si,L_{2\Si})$ in terms of the BKS pairing (\ref{e:BKS}), we find that 
\begin{equation} \label{e:ecyl}
Z_{(\Si \times I , L_{2 \Si},n)} \, = \,  \te^{\frac{\pi \ti}{4}   n} \, \sum\limits_{\bq} \overline{v{\sst (\Si,L_{\Si})}_{\bq}} \otimes v{\sst (\Si,L_{\Si})}_{\bq} \, .
\end{equation}
Thus $Z_{(\Si \times I , L_{2 \Si},0)} = Id \in \text{Hom}[\CH(\Si,L_{\Si}),\CH(\Si,L_{\Si})]$.

\nin (e) {\em Gluing}.
The moduli space of flat $\BT$-connections on $\pXC = \pX \sqcup (-\Si) \sqcup \Si$ is the symplectic manifold $(\CM_{\pXC},\om_{\pXC})$ with
\begin{align*}
\CM_{\pXC} \, &= \, \CM_{\pX} \times \CM_{\Si} \times \CM_{\Si} \\
\om_{\pXC} \, &= \, pr^*_{\pX}(\om_{\pX}) \, + \, pr^*_{-\Si}(-\om_{\Si}) \, + \, pr^*_{\Si}(\om_{\Si}) \, ,
\end{align*}
where $pr_{\pX}, pr_{-\Si}$ and $pr_{\Si}$ denote the projection maps from $ \CM_{\pX} \times \CM_{\Si} \times \CM_{\Si}$ onto the first, second and third factor, respectively.
Let $\De \subset \CM_{\Si} \times \CM_{\Si}$ denote the diagonal. 
We claim that the manifold
\begin{equation*}
C \, = \, \CM_{\pX} \times \De \, \subset \, \CM_{\pXC}
\end{equation*}
is a coisotropic submanifold of $(\CM_{\pXC},\om_{\pXC})$.
Let $T C^{\perp}$ denote the orthogonal complement of the tangent bundle $TC$ with respect to $\om_{\pXC}$, i.e.
\begin{equation*}
TC^{\perp} \, = \, \{ v \in T \CM_{\pXC} \, \mid \, \om_{\pXC}(v,w) =0, \, \text{ for all }\, w \in TC \} \, .
\end{equation*}
For any vector $v =(v_{\pX},v_{\Si},v_{\Si}') \in T \CM_{\pXC}$ such that $v \in T C^{\perp}$ we must have $0 = \om_{\pXC}(v,w) = \om_{\pX}(v_{\pX},w_{\pX}) -\om_{\Si}(v_{\Si},w_{\Si}) + \om_{\Si}(v_{\Si}',w_{\Si})$, for all vectors $w=(w_{\pX},w_{\Si},w_{\Si}) \in TC$.
This is true if and only if $v_{\pX} =0$ and $v_{\Si}= v_{\Si}'$.
Thus $TC^{\perp} \subset TC$ and $C \subset \CM_{\pXC}$ is coisotropic.
Moreover, $TC^{\perp}$ is the tangent bundle to an isotropic foliation $C^{\perp}$ of $C$.
The reduced symplectic manifold $C/C^{\perp}$ is identified with $(\CM_{\pX},\om_{\pX})$ and the quotient map $C \ra C/C^{\perp}$ with the natural projection $pr_{\pX}: \CM_{\pX} \times \De \ra \CM_{\pX}$.

With the usual notations let $\La_{\XC}$ be the image of $\CM_{\XC}$ under the restriction map $r_{\XC}: \CM_{\XC} \ra \CM_{\pXC}$.
According to (\ref{p:resmap}), $\La_{\XC}$ is a Lagrangian submanifold of  $(\CM_{\pXC},\om_{\pXC})$ and we consider its reduction \cite{We} relative to the coisotropic submanifold $C$.
The intersection $\La_{\XC} \cap C$ is a manifold and the tangent bundles satisfy $T(\La_{\XC} \cap C) = T\La_{\XC} \cap TC$. 
The reduction of $\La_{\XC}$ relative to $C$ is the image of $\La_{\XC} \cap C$ under the quotient map $C \ra C/C^{\perp}$.
It is a Lagrangian submanifold of $C/C^{\perp}$ which coincides with $\La_X \subset \CM_{\pX}$ under the identification of $C/C^{\perp}$ with $\CM_{\pX}$.
That is,
\begin{equation*}
\La_X \, = \, pr_{\pX}(\La_{\XC} \cap C) \, .
\end{equation*}
The kernel of the differential of $p_{\pX} : \La_{\XC} \cap C \ra \La_X$ is the vector bundle $T \La_{\XC} \cap TC^{\perp}$.

Now let $g : \XC \ra X$ denote the gluing map and $g^* : \CM_X \ra \CM_{\XC}$ the induced map between the moduli spaces of flat $\BT$-connections.
If $r_X : \CM_X \ra \CM_{\pX}$ and $r_{\Si} : \CM_X \ra \CM_{\Si}$ are the restriction maps from connections on $X$ to connections on $\pX$ and $\Si$, respectively, then we get the commutative diagram 
\begin{equation} \label{e:commut}
\begin{CD}
\CM_X @>{r_X \times r_{\Si}}>> \CM_{\pX} \times \CM_{\Si} \\
@V{g^*}VV @VV{id \times i_{\De}}V \\
\CM_{\XC} @>{r_{\XC}}>> \CM_{\pXC} = \CM_{\pX} \times \CM_{\Si} \times \CM_{\Si} \\
\end{CD}
\end{equation}
where $i_{\De} : \CM_{\Si} \ra \CM_{\Si} \times \CM_{\Si}$ is the diagonal map $i_{\De}([\eta]) = ([\eta],[\eta]) \in \De$.
We note that
\begin{equation} \label{e:lacut}
\La_{\XC} \cap C \, = \, (r_{\XC} \circ g^*) (\CM_X)\, .
\end{equation}
Recall that $\CM_{\XC} = \underset{p^{cut} \in \text{Tors} H^2(\XC;\BZ)}{\sqcup} \CM_{\XC,p^{cut}}$ and $\CM_X = \underset{p \in \text{Tors} H^2(X;\BZ)}{\sqcup} \CM_{X,p}$.
Thus for each $p^{cut}$ and each connected component $\big[\La_{\XC} \cap C]_{c}$ of $\La_{\XC} \cap C$ there is a unique $p$ such that
\begin{equation} \label{e:inters}
r^{-1}_{\XC} \big(\big[\La_{\XC} \cap C]_{c} \big) \cap \CM_{\XC,p^{cut}} \, = \, g^*(\CM_{X,p}) \, .
\end{equation}
Then (\ref{e:inters})  and (\ref{p:modtor} (ii)) imply the following relation:
\begin{equation} \label{e:count}
\big[ \# \pi_0(\La_{\XC} \cap C)] \cdot \big[ \# \text{Tors} \, H^2(\XC;\BZ) \big] \, = \, \big[ \# \text{Tors} H^2(X;\BZ) \big]
\end{equation}
With the usual notations let 
\begin{align*}
L_{\XC} \, &= \, \text{Im} \,\{ H^1(\XC;\BR) \ra H^1(\pXC;\BR) \} \\
L_X \, &= \, \text{Im} \,\{ H^1(X;\BR) \ra H^1(\pX;\BR) \} \notag
\end{align*}
and let us introduce the following Lagrangian subspaces of the symplectic vector space $H^1(\XC;\BR) = H^1(X;\BR) \oplus H^1(-\Si;\BR) \oplus H^1(\Si;\BR)$:
\begin{align} \label{e:Lagr}
L' \, &= \, L_X \oplus L_{\Si} \oplus L_{\Si} \\
L'' \, &= \, L \oplus  L_{\Si} \oplus L_{\Si} \, . \notag
\end{align}
The vector corresponding to the e-3-manifold $(\XC,L'',n^{cut})$ is the vector
\begin{equation} \label{e:zxcut1}
Z_{(\XC,L'',n^{cut})} \, = \, \te^{\frac{\pi \ti}{4} n^{cut}} \, F_{L'' L_{\XC}}(Z_{\XC}) \, \in \, \CH(\pXC,L'')\, ,
\end{equation}
obtained as the image of the standard vector 
\begin{equation} \label{e:zxcut}
Z_{\XC} \, = \, \frac{k^{m_{\XC}}}{[ \# \text{Tors} \, H^2(\XC;\BZ) ]} \; \si_{\XC} \otimes \mu_{\XC} \, \in \, \CH(\pXC,L_{\XC})
\end{equation}
defined in Sect.\ref{s:qth}. Using the composition law
\begin{equation*}
F_{L'' L_{\XC}} \, = \, \te^{\frac{\pi \ti}{4} \tau(L_{\XC},L',L'')} \, F_{L'' L'} \circ F_{L' L_{\XC}}
\end{equation*}
we can rewrite (\ref{e:zxcut1}) as
\begin{equation} \label{e:zxcut2}
Z_{(\XC,L'',n^{cut})} \, = \, \te^{\frac{\pi \ti}{4}[n^{cut}+ \tau(L_{\XC},L',L'')]}  \, F_{L'' L'} \circ F_{L' L_{\XC}} (Z_{\XC})
\end{equation}
Let us choose unitary bases $\{ v{\sst (\pX,L_X)}_{\bq} \}$ for $\CH(\pX,L_X)$, $\{ v{\sst (\pX,L)}_{\bq'} \}$ for $\CH(\pX,L)$ and
$\{ v{\sst (\Si,L_{\Si})}_{\bell} \}$ for $\CH(\Si,L_{\Si})$.
Then we have
\begin{equation} \label{e:FLL}
\begin{split}
F_{L'' L'} \circ F_{L' L_{\XC}} (Z_{\XC}) \, &= \, F_{L'' L'} \Big[ \sum\limits_{\bq,\bell,\bell'} \, M_{\bq \bell \bell'}\, v{\sst (\pX,L_X)}_{\bq} \otimes 
\overline{v{\sst (\Si,L_{\Si})}_{\bell}} \otimes v{\sst (\Si,L_{\Si})}_{\bell'} \Big] \\
&= \, \sum\limits_{\bq,\bell,\bell'} \, M_{\bq \bell \bell'} \,\Big[ F_{L L_X} \, v{\sst (\pX,L_X)}_{\bq} \Big] \, \otimes 
\overline{v{\sst (\Si,L_{\Si})}_{\bell}} \otimes v{\sst (\Si,L_{\Si})}_{\bell'}\, ,
\end{split}
\end{equation}
where $F_{L L_X} : \CH(\pX,L_X) \ra \CH(\pX,L)$ and 
\begin{equation*}
 M_{\bq \bell \bell'} \, = \, \Big\langle  v{\sst (\pX,L_X)}_{\bq} \otimes 
\overline{v{\sst (\Si,L_{\Si})}_{\bell}} \otimes v{\sst (\Si,L_{\Si})}_{\bell'} \, , \, F_{L' L_{\XC}}(Z_{\XC}) \Big\rangle_{\CH(\pXC,L')} \, .
\end{equation*}
The last equality in (\ref{e:FLL}) follows from the definition (\ref{e:BKSisom}) of the intertwining isomorphism $F_{L'' L'}$ through the BKS pairing (\ref{e:BKS}) which, having in view the special form (\ref{e:Lagr}) of the Lagrangian subspaces $L'$ and $L''$, leads to the expression $F_{L'' L'} = F_{L L_X} \otimes Id \otimes Id$.

Applying the contraction operator $\text{Tr}_{\Si}$ to $Z_{(\XC,L'',n^{cut})}$ and making use of the expression (\ref{e:FLL}), we get
\begin{equation*}
\text{Tr}_{\Si} \Big[ Z_{(\XC,L'',n^{cut})} \Big] \, = \, 
\te^{\frac{\pi \ti}{4}[n^{cut}+ \tau(L_{\XC},L',L'')]} \, \sum\limits_{\bq, \bell} \, M_{\bq \bell \bell} \, F_{L L_X} \big[ v{\sst (\pX,L_X)}_{\bq} \big]\, .
\end{equation*}
Now let $L_{\De} \subset H^1(-\Si;\BR) \oplus H^1(\Si;\BR)$ be the diagonal. 
It is obviously a rational Lagrangian subspace for the  symplectic form $ pr^*_{-\Si}(-\om_{\Si}) \, + \, pr^*_{\Si}(\om_{\Si})$ on $H^1(-\Si;\BR) \oplus H^1(\Si;\BR)$.
We introduce the notation $L_{2\Si} = L_{\Si} \oplus L_{\Si}$ and define $L_{\De}' = L_X \oplus L_{\De}$.
Recall that to the cylinder $\Si \times I$ there corresponds the standard vector $Z_{\Si \times I}$ in $\CH((-\Si) \sqcup \Si, L_{\De})$ given by the expression (\ref{e:cylinder}). 
Associated to the e-3-manifold $(\Si \times I, L_{2\Si}, 0)$ we have, according to (\ref{e:ecyl}), the vector 
\begin{equation*}
Z_{(\Si \times I, L_{2\Si}, 0)} \, = \, F_{L_{2\Si},L_{\De}}(Z_{\Si \times I}) \, = \, \sum\limits_{\bell} \overline{v{\sst (\Si,L_{\Si})}_{\bell}} \otimes v{\sst (\Si,L_{\Si})}_{\bell}
\end{equation*}
in $\CH((-\Si) \sqcup \Si, L_{2\Si})$.
Using this we can write the following:
\begin{equation*}
\begin{split}
\sum\limits_{\bell} \, M_{\bq \bell \bell} \, &= \, \sum_{\bell} \Big\langle  v{\sst (\pX,L_X)}_{\bq} \otimes 
\overline{v{\sst (\Si,L_{\Si})}_{\bell}} \otimes v{\sst (\Si,L_{\Si})}_{\bell} \, , \, F_{L' L_{\XC}}(Z_{\XC}) \Big\rangle_{\CH(\pXC,L')}\\
&= \,  \Big\langle  v{\sst (\pX,L_X)}_{\bq} \otimes Z_{(\Si \times I, L_{2\Si}, 0)}\, , \, F_{L' L_{\XC}}(Z_{\XC}) \Big\rangle_{\CH(\pXC,L')}\\
&= \,  \Big\langle  v{\sst (\pX,L_X)}_{\bq} \otimes F_{L_{2\Si},L_{\De}} (Z_{\Si \times I}) \, , \, F_{L' L_{\XC}}(Z_{\XC}) \Big\rangle_{\CH(\pXC,L')}\\
\end{split}
\end{equation*}
Now, since $F_{L' L_{\XC}} = \te^{\frac{\pi \ti}{4} \tau(L_{\XC}, L_{\De}',L')} \, F_{L' L_{\De}'} \circ F_{L_{\De}', L_{\XC}}$
and since the operator $F_{L' L_{\De}'} = Id \otimes F_{L_{2 \Si}, L_{\De}}$, we obtain
\begin{equation*}
\begin{split}
\sum\limits_{\bell}  M_{\bq \bell \bell}  &= \te^{\frac{\pi \ti}{4} \tau(L_{\XC}, L_{\De}',L')}\, \times \\ 
 &\times \, \Big\langle F_{L' L_{\De}'} \big[ v{\sst (\pX,L_X)}_{\bq} \otimes Z_{\Si \times I} \big] \, , \, F_{L' L_{\De}'} \circ F_{L_{\De}', L_{\XC}} (Z_{\XC}) \Big\rangle_{\CH(\pXC,L')} \\
&= \te^{\frac{\pi \ti}{4} \tau(L_{\XC}, L_{\De}',L')} \,
 \Big\langle  v{\sst (\pX,L_X)}_{\bq} \otimes Z_{\Si \times I}  \, , \, F_{L_{\De}', L_{\XC}} (Z_{\XC}) \Big\rangle_{\CH(\pXC,L_{\De}')} \\
&= \frac{ k^{m_{\XC} + \frac{1}{4} \dim H^1(\Si;\BR)} }{\big[ \# \text{Tors} \, H^2(\XC;\BZ) \big] } \; 
\te^{\frac{\pi \ti}{4} \tau(L_{\XC}, L_{\De}',L')} \, \times \\
& \, \times \int\limits_{\La_{\XC} \cap ((\La_{X})_{\bq} \times \De)} \,
\big( s{\sst (\pX)}_{\bq} \otimes \si_{\Si \times I} \, , \, \si_{\XC} \big) \; \big( \de_X \otimes \mu_{\Si \times I} \big) * \mu_{\XC}
\end{split}
\end{equation*}
The last equality follows from the BKS pairing formula  defining the operator $F_{L_{\De}',L_{\XC}}$ and the expressions (\ref{e:zxcut}) for $Z_{\XC}$ and (\ref{e:cylinder}) for $Z_{\Si \times I}$.
We also wrote $ v{\sst (\pX,L_X)}_{\bq} =  s{\sst (\pX)}_{\bq} \otimes \de_X$ where $\de_X$ is the invariant $\frac{1}{2}$-density on $\CP_X$ such that $\int\limits_{(\La_{X})_{\bq}} \de^2_X =1$, for any leaf $(\La_{X})_{\bq}$ of the polarization $\CP_X$ on $\CM_{\pX}$.

As previously shown we have the reduction map $\La_{\XC} \cap C \ra \La_X$ and, since $\La_X \cap (\La_{X})_{\bq} \neq \emptyset$ if and only if $ (\La_{X})_{\bq} = \La_X$, we conclude that
\begin{equation*}
\La_{\XC} \cap ((\La_{X})_{\bq} \times \De) = \emptyset \, , \quad \text{for all leaves } \,  (\La_{X})_{\bq} \neq \La_X \, .
\end{equation*}

Let us write $\mu_{\XC} = a_{\XC} \de_{\XC}$ with $\de_{\XC}$ the invariant $\frac{1}{2}$-density on $\CP_{\XC}$ satisfying $\int\limits_{\La_{\XC}} \de_{\XC}^2=1$.
We also recall from (\ref{e:cylinder}) that $\mu_{\Si \times I} = (T_{\Si})^{\frac{1}{2}} = \Bigl| \frac{(\om_{\Si})^g}{g!} \Bigr|^{\frac{1}{2}}$ and  that we have $\int\limits_{\De} \frac{(\om_{\Si})^g}{g!} = 1$, where $g = \frac{1}{2} \dim H^1(\Si;\BR)$.
Then, making use of the results in (\cite{Ma},\S 4), we find that for each connected component $\big[ \La_{\XC} \cap C \big]_c$ of the manifold $\La_{\XC} \cap C = \La_{\XC} \cap (\La_X \times \De)$ we have
\begin{equation*}
\begin{split}
\int\limits_{ \big[\La_{\XC} \cap C \big]_c}
\big( \de_X \otimes   T_{\Si}^{\frac{1}{2}} \big) * \de_{\XC} \, &= \, 
k^{ \frac{1}{2} [ \dim(\La_{\XC} \cap C) - \frac{1}{2} \dim H^1(\pXC;\BR)] } \, \times \\ 
& \times \, \big[ \# \pi_0 (\La_{\XC} \cap C) \big]^{-\frac{1}{2}}
\end{split}
\end{equation*}
Let us introduce the notation $ v{\sst (\pX,L_X)} =  v{\sst (\pX,L_X)}_{\bq}$ if $(\La_{X})_{\bq} = \La_X$. 
We also note that the function $(s{\sst (\pX)} \otimes \si_{\Si \times I} , \si_{\pXC})$ on $\La_{\pXC} \cap C$ is constant on each connected component.
Therefore, summarizing the previous results and observations and using the fact that $ \si_{\XC} = \sum\limits_{p^{cut} \in \text{Tors} H^2(\XC;\BZ)} \si_{\XC,p^{cut}}$, we can write
\begin{equation} \label{e:trexplicit}
\begin{split}
\text{Tr}_{\Si} \Big[ Z_{(\XC,L'',n^{cut})} \Big] \, &= \, 
\te^{\frac{\pi \ti}{4} [ n^{cut} + \tau(L_{\XC},L',L'') + \tau(L_{\XC},L_{\De}',L') ]}  \\
&\times \, \frac{a_{\XC}}{\big[ \# \text{Tors} H^2(\XC;\BZ) \big]} 
\:  \big[ \# \pi_0 (\La_{\XC} \cap C) \big]^{-\frac{1}{2}} \\
&\times \,  k^{ m_{\XC} + \frac{1}{4} \dim H^1(\Si;\BR) + \frac{1}{2} \dim (\La_{\XC} \cap C) - \frac{1}{4} \dim H^1(\pXC;\BR) } \\
& \times \Big[ \sum\limits_c \sum\limits_{p^{cut} \in \text{Tors} H^2(\XC;\BZ)}
 \big( s{\sst (\pX)} \otimes \si_{\Si\times I} , \si_{\XC,p^{cut}} \big)_c \Big] \\
&\times \,  F_{L L_X} (v{\sst (\pX,L_X)} )\, ,
\end{split}
\end{equation}
where $\big( s{\sst (\pX)} \otimes \si_{\Si\times I} , \si_{\XC,p^{cut}} \big)_c$ denotes the value of the function on the component $\big[\La_{\XC} \cap C]_{c}$.
Then let $[\Th] \in \CM_{X,p}$ with $[\eta] = [\Th \big|_{\Si}]$ and $[\partial \Th] = [ \Th \big|_{\pX}]$ and let $[\Th^{cut}] = [g^* \Th]$.
Using the classical gluing formula (\ref{e:tr}) and the definitions of $\si_{\Si \times I}, \, \si_{X,p}$ and $\si_{\XC,p^{cut}}$, we find that
\begin{equation*}
\begin{split}
\Big( \si_{X,p}([\partial \Th]) \otimes \si_{\Si \times I}([\eta] \sqcup [\eta]) \, &, \, \si_{\XC,p^{cut}}([\partial \Th] \sqcup [\eta] \sqcup [\eta]) \Big) \\
& = \, \Big( \te^{\pi \ti k S_{X,P}(\Th)}\, , \, \text{Tr}_{\eta} \big[
\te^{\pi \ti k S_{\XC,P^{cut}}(\Th^{cut}) } \big] \Big) = 1\, .
\end{split}
\end{equation*}
Since both $\si_{X,p}$ and $s{\sst (\pX)}$ are unitary and covariantly constant sections of $\CL_{\pX} \otimes \left| \text{Det} \CP_X^* \right|^{\frac{1}{2}}$ over $\La_X$, they differ by a constant phase factor.
In view of this observation, of the previous relation and of the relations (\ref{e:inters})-(\ref{e:count}), the expression (\ref{e:trexplicit}) becomes
\begin{equation} \label{e:trace}
\begin{split}
\text{Tr}_{\Si} \Big[ Z_{(\XC,L'',n^{cut})} \Big] \, &= \, 
\te^{\frac{\pi \ti}{4} [ n^{cut} + \tau(L_{\XC},L',L'') + \tau(L_{\XC},L_{\De}',L') ]} \\
&\times \, \frac{a_{\XC}}{ \big[ \# \text{Tors} H^2(\XC;\BZ) \big]^{\frac{1}{2}} \big[ \# \text{Tors} H^2(X;\BZ) \big]^{\frac{1}{2}} } \\  
&\times \,  k^{ m_{\XC} + \frac{1}{4} \dim H^1(\Si;\BR) + \frac{1}{2} \dim (\La_{\XC} \cap C) - \frac{1}{4} \dim H^1(\pXC;\BR) } \\
&\times \, \sum\limits_{p \in \text{Tors} H^2(X;\BZ)} \,
 F_{L L_X} (\si_{X,p} \otimes \de_X )\, .
\end{split}
\end{equation}

Let $a_X$ denote the real constant such that $\mu_X = a_X \de_X$.
We claim that the following relations hold:
\begin{equation} \label{e:ax}
a_{\XC}^2 \, = \, a_X^2 \, \frac{\big[ \# \text{Tors} H^2(\XC;\BZ) \big] }{\big[ \# \text{Tors} H^2(X;\BZ) \big] }
\end{equation}
and
\begin{equation} \label{e:mxeq}
m_{\XC} +  \frac{1}{4} \dim H^1(\Si;\BR) + \frac{1}{2} \dim (\La_{\XC} \cap C) - \frac{1}{4} \dim H^1(\pXC;\BR)  \, = \, m_X \, .
\end{equation}
We postpone for a moment their proof and note first that by  using the definition of $n^{cut}$ and the cocycle relation \cite{LV} for the Maslov-Kashiwara index $\tau$ we have:
\begin{equation}
\begin{split}
&n^{cut} + \tau(L_{\XC}, L',L'') + \tau(L_{\XC},L_{\De}',L') \pmod{8} \\
&= \, n - \tau(L'', L_{\XC}, L_{\De}') + \tau(L_{\XC}, L',L'') + 
\tau(L_{\XC},L_{\De}',L') \pmod{8}  \\
&= \, n + \tau(L',L'',L_{\De}') \pmod{8} .
\end{split}
\end{equation}
From the symplectic additivity property of $\tau$ we find
\begin{equation} \label{e:tau}
\begin{split}
\tau(L',L'',L_{\De}') \, &= \, \tau(L_X \oplus L_{\Si} \oplus L_{\Si}, L \oplus L_{\Si} \oplus L_{\Si}, L_X \oplus L_{\De}) \\
&= \, \tau(L_X,L,L_X) + \tau(L_{\Si} \oplus L_{\Si}, L_{\Si} \oplus L_{\Si}, L_{\De}) =0 \, .
\end{split}
\end{equation}
Putting the results (\ref{e:trace})-(\ref{e:tau}) together we obtain
\begin{equation*}
\begin{split}
\text{Tr}_{\Si} \Big[ Z_{(\XC,L'',n^{cut})} \Big] \, &= \, 
\te^{\frac{\pi \ti}{4} n} \,
\frac{k^{m_X}}{\big[ \# \text{Tors} H^2(X;\BZ) \big]} \; 
 F_{L L_X} (\si_X \otimes \mu_X ) \\
&= \, \te^{\frac{\pi \ti}{4} n} \, F_{L L_X} (Z_X) \, = \, Z_{(X,L,n)}
\end{split}
\end{equation*}
which proves the gluing property (\ref{t:TQFT} (e)).

Let us return now to the proof of the equality (\ref{e:mxeq}).
Recall that $m_X$ is defined by (\ref{e:mx}) and similarly $m_{\XC}$ is given by
\begin{align*}
m_{\XC} \, &= \,  \frac{1}{4} \big( \dim H^1(\XC;\BR) + \dim H^1(\XC,\pXC;\BR) \\
       & \quad - \dim H^0(\XC;\BR) - \dim H^0(\XC,\pXC;\BR) \big)\, .
\end{align*}
We have the following relations: \\
(i) The Mayer-Vietoris cohomology sequence for the spaces $X, \XC, \Si$
\begin{equation} \label{e:MVseq}
\cdots \ra  H^i(X;\BR) \ra H^i(\XC;\BR) \ra H^i(\Si;\BR) \ra H^{i+1}(X;\BR) \ra \cdots
\end{equation}
together with Poincar\'{e} duality imply the relation
\begin{equation} \label{e:i}
\begin{split}
\dim H^1(\XC;\BR) & - \dim H^1(\XC,\pXC;\BR)  \\
&\;  - \dim H^0(\XC;\BR) + \dim H^0(\XC,\pXC;\BR) \\
&= \dim H^1(X;\BR) - \dim H^1(X,\pX;\BR) - \dim H^0(X;\BR) \\
&\; + \dim H^0(X,\pX;\BR) + \dim H^1(\Si;\BR) - 2 \dim H^0(\Si;\BR)\,.
\end{split}
\end{equation}
(ii) The exact cohomology sequence for the pair of spaces $(\pX \sqcup \Si) \subset X$ and (\ref{e:commut})-(\ref{e:lacut}) give
\begin{align*}
0 \ra H^0(X;\BR) \ra H^0(\pX \sqcup \Si;\BR) & \ra H^1(X,\pX \sqcup \Si;\BR) \ra \\
& \ra H^1(X;\BR) \xrightarrow[]{\dot{r}_X \times (\dot{r}_{\Si}, \dot{r}_{\Si})}  T(\La_{\XC} \cap C) \ra 0 
\end{align*}
Together with the isomorphism $H^1(X,\pX \sqcup \Si;\BR) \cong H^1(\XC,\pXC;\BR)$ it implies that
\begin{equation} \label{e:ii}
\begin{split}
\dim (\La_{\XC} \cap C) \, = \, & \dim H^1(X;\BR) - \dim H^1(\XC,\pXC;\BR) \\
& + \dim H^0(\pX;\BR) + \dim H^0(\Si;\BR) - \dim H^0(X;\BR)\,.
\end{split}
\end{equation}
(iii) From the exact sequence (\ref{e:exfibr})
and the fact that $\dim \La_X = \frac{1}{2} \dim H^1(\pX;\BR)$ we obtain
\begin{equation} \label{e:iii}
\begin{split}
\frac{1}{2} \dim H^1(\pX;\BR) \, &= \, \dim H^1(X;\BR) -\dim H^1(X,\pX;\BR) \\
&\, + \dim H^0(\pX;\BR)  - \dim H^0(X;\BR) \\
&\, + \dim H^0(X,\pX;\BR) 
\end{split}
\end{equation}
The claimed formula (\ref{e:mxeq}) follows then from the relations (\ref{e:i}), (\ref{e:ii}) and (\ref{e:iii}).

Finally let us prove the relation (\ref{e:ax}).
Referring back to the definition (\ref{e:hdensid}) we have
\begin{align} \label{e:rel1}
\mu_X \, &= \, \Bigg[ \int\limits_{H^1(X,\pX;\BT)} \bw \Bigg] \, (T_X)^{\frac{1}{2}} \otimes w^{-1} \\
\mu_{\XC} \, &= \, \Bigg[ \int\limits_{H^1(\XC,\pXC;\BT)} \bw_{cut} \Bigg] \, (T_{\XC})^{\frac{1}{2}} \otimes w_{cut}^{-1} \, , \notag
\end{align}
with $ w \in \left| \text{Det} H^1(X,\pX;\BR)^* \right| $ and $w_{cut} \in \left| \text{Det} H^1(\XC,\pXC;\BR)^* \right| $.
On the other hand we introduced $a_X$ and $a_{\XC}$ such that
\begin{align} \label{e:rel2}
\mu_X \, &= \, a_X \de_X \\
\mu_{\XC} \, &= \, a_{\XC} \de_{\XC} \notag
\end{align}
The R-torsion $T_X$ of $X$ is related to the R-torsions $T_{\XC}$ of $\XC$ and $T_{\Si}$ of $\Si$ by the gluing formula \cite{V}
\begin{equation} \label{e:gluetors}
T_X \, = \, T_{\XC} \otimes (T_{\Si})^{-1} \,.
\end{equation}
The identification between the l.h.s. and the r.h.s. is made through the isomorphism of determinant lines
\begin{equation*}
\left| \text{Det} H^{\bullet}(X;\BR)^* \right| \cong 
\left| \text{Det} H^{\bullet}(\XC;\BR)^* \right| \otimes 
\left| \text{Det} H^{\bullet}(\Si;\BR) \right|
\end{equation*}
arising from the Mayer-Vietoris sequence (\ref{e:MVseq}).
Using the exact sequence (\ref{e:exfibr}) and, due to the normalization of the density $\de_X^2$ on $\La_X$, we see that
\begin{equation*}
\frac{w \otimes \de_X^2}{\Bigg[ \int\limits_{H^1(X,\pX;\BT)} \bw \Bigg]} \, \in
\, \left| \text{Det} H^1(X;\BR)^* \right|
\end{equation*}
defines an invariant density on the group $H^1(X;\BT)$, which gives $H^1(X;\BT)$ volume 1.
Similarly on $H^1(\XC;\BT)$ we have the density defined by
\begin{equation*}
\frac{w_{cut} \otimes \de_{\XC}^2}{\Bigg[ \int\limits_{H^1(\XC,\pXC;\BT)} \bw_{cut} \Bigg]} \, \in
\, \left| \text{Det} H^1(\XC;\BR)^* \right|
\end{equation*}
The Mayer-Vietoris sequence (\ref{e:MVseq}) and Poincar\'{e} duality induce the isomorphism of determinant lines
\begin{equation*}
\begin{split}
\left| \text{Det} H^1(X;\BR)^* \right| \, &\cong \,
\left| \text{Det} H^1(\XC;\BR)^* \right| \otimes 
\left| \text{Det} H^1(\XC,\pXC;\BR)^* \right| \\
&\, \otimes 
\left| \text{Det} H^1(\Si;\BR) \right| \otimes
\left| \text{Det} H^1(X,\pX;\BR) \right|
\end{split}
\end{equation*}
Thus we can define an invariant density $\rho_X$ on the group $H^1(X;\BT)$ by setting
\begin{equation*}
\rho_X = 
\frac{w_{cut} \otimes \de_{\XC}^2}{\Bigg[ \int\limits_{H^1(\XC,\pXC;\BT)} \bw_{cut} \Bigg]} \otimes
\frac{w_{cut}}{\Bigg[ \int\limits_{H^1(\XC,\pXC;\BT)} \bw_{cut} \Bigg]} \otimes
T_{\Si}^{-1} \otimes
\frac{w^{-1}}{\Bigg[ \int\limits_{H^1(X,\pX;\BT)} \bw \Bigg]^{-1}}
\end{equation*}
From the normalization of the densities on the r.h.s. and from the Mayer-Vietoris cohomology sequence with $\BT$-coefficients
for $X,\XC,\Si$
\begin{align*}
\cdots \ra H^1(X;\BT) & \ra H^1(\XC;\BT) \ra H^1(\Si;\BT) \\
&\ra H^2(X;\BT) \ra H^2(\XC;\BT) \ra H^2(\Si;\BT) \ra 0
\end{align*}
we see that, since $\pi_0 [ H^2(X;\BT)] \cong \pi_0 [ H^2(\XC;\BT)] = 0
$ and since, as shown in Prop. (\ref{p:modtor}),
$\pi_0 [ H^1(X,\pX;\BT)] \cong \text{Tors} H^2(X,\pX;\BZ) \cong \text{Tors} H^2(X;\BZ)$ and similarly
$\pi_0 [ H^1(\XC,\pXC;\BT)] \cong \text{Tors} H^2(\XC,\pXC;\BZ) \cong \text{Tors} H^2(\XC;\BZ)$,
the density $\rho_X$ gives $H^1(X;\BT)$ volume equal to $\frac{ \big[ \# \text{Tors} H^2(X;\BZ) \big]}{\big[ \# \text{Tors} H^2(\XC;\BZ) \big]}$.
Making use of the gluing formula (\ref{e:gluetors}) and the relations (\ref{e:rel1}) and (\ref{e:rel2}) we can write
\begin{equation*}
\begin{split}
\frac{w \otimes \de_X^2}{\Bigg[ \int\limits_{H^1(X,\pX;\BT)} \bw \Bigg]} \, &= \, \frac{1}{a_X^2} \, T_X \otimes \frac{w^{-1}}{\Bigg[ \int\limits_{H^1(X,\pX;\BT)} \bw \Bigg]^{-1}} \\  
&= \, \frac{1}{a_X^2} \, T_{\XC} \otimes (T_{\Si})^{-1} \otimes \frac{w^{-1}}{\Bigg[ \int\limits_{H^1(X,\pX;\BT)} \bw \Bigg]^{-1}} \\
&= \, \frac{a_{\XC}^2}{a_X^2} \, 
\frac{w_{cut}^2 \otimes \de_{\XC}^2}{\Bigg[ \int\limits_{H^1(\XC,\pXC;\BT)} \bw_{cut} \Bigg]^2} \otimes
(T_{\Si})^{-1} \otimes \frac{w^{-1}}{\Bigg[ \int\limits_{H^1(X,\pX;\BT)} \bw \Bigg]^{-1}} \\ 
& = \, \frac{a_{\XC}^2}{a_X^2} \, \rho_X
\end{split}
\end{equation*}
In view of the previous observation regarding the volume of $H^1(X;\BT)$ computed with $\rho_X$ and the fact that the density of the above equation gives $H^1(X;\BT)$ volume 1, the claimed relation (\ref{e:ax}) is proved.
\end{proof}

The mapping class group $\Ga_{\Si}$ of a closed oriented 2-manifold $\Si$ is the group $\text{Diff}_{+}(\Si)/ \text{Diff}_0(\Si)$ of isotopy classes of orientation preserving diffeomorphisms of $\Si$.
The mapping class group $\Ga_{(\Si,L)}$ of an e-2-manifold $(\Si,L)$ is a central extension by $\BZ/ 8\BZ$ of $\Ga_{\Si}$:
\begin{equation*} 
1 \lra \BZ/ 8\BZ \lra \Ga_{(\Si,L)} \lra \Ga_{\Si} \lra 1
\end{equation*}
As a set $\Ga_{(\Si,L)} = \Ga_{\Si} \times \BZ/ 8\BZ$ and the composition law in $\Ga_{(\Si,L)}$ is, according to (\ref{d:cemor2}), given by
\begin{equation*}
([h],m) ([h'],m') = \big( [h h'],\, m+m'+\tau(L'',h^{\prime *} L', (h h')^* L) \pmod{8} \big)
\end{equation*}

In the course of proving Theorem (\ref{t:TQFT} (a)) we  showed that, given a rational Lagrangian subspace $L \subset H^1(\Si;\BR)$, each element $[h]$ of the mapping class group $\Ga_{\Si}$ determines a unitary map of Hilbert spaces $h^* : \CH(\Si,L) \ra \CH(\Si, h^* L)$.
The composition of this map with the isomorphism $F_{L,h^* L} : \CH(\Si, h^* L) \ra \CH(\Si,L)$ induced by the BKS pairing is a unitary operator $U_L([h]) = F_{L,h^* L} \circ h^*$ on the Hilbert space $\CH(\Si,L)$.
The assignment 
\begin{equation*}
[h] \, \mapsto \, U_L([h])
\end{equation*}
defines a unitary projective representation of $\Ga_{\Si}$ on $\CH(\Si,L)$.
For the group $\Ga_{(\Si,L)}$, the assignment 
\begin{equation*}
([h],m) \, \mapsto \, U_L([h],m) = \te^{\frac{\pi \ti}{4} m} U_L([h])
\end{equation*}
determines a unitary representation of this group on the Hilbert space $\CH(\Si,L)$.

Since  the mapping class group $\Ga_{\Si}$ acts on $H^1(\Si;\BZ)$, there is a natural homomorphism from $\Ga_{\Si}$ to the group $Sp(\CZ)$ of symplectic transformations of $( H^1(\Si;\BR), \om_{\Si} )$ which preserve the integer lattice $\CZ= H^1(\Si;\BZ)$.
According to the results in (\cite{Ma},\S 9.1), for each rational Lagrangian subspace $L \subset H^1(\Si;\BR)$, there is a projective unitary representation of $Sp(\CZ)$ on the Hilbert space $\CH(\Si,L)$.
Let us choose an integer symplectic basis $(\bw; \bw^{\prime})$ for $H^1(\Si;\BR)$, with $w_1, \dots , w_g$ spanning $L$, where $g= \frac{1}{2} \dim H^1(\Si;\BR)$.
As shown in (\cite{Ma},\S 3), the choice of such a basis uniquely determines a unitary basis $\{ v{\sst (\Si,L)}_{\bq} \}_{\bq \in (\BZ/ k\BZ)^g}$ for $\CH(\Si,L)$.
The generators of $Sp(\CZ)$ are the elements with matrix form
\begin{equation*} 
\al = \begin{pmatrix} A & 0 \\ 0 & \,{ ^t A^{-1}} \end{pmatrix}
\quad
\be = \begin{pmatrix} I & B \\ 0 & I \end{pmatrix} \quad
\ga = \begin{pmatrix} 0 & I \\ -I & 0 \end{pmatrix}
\end{equation*}
with respect to the basis $(\bw; \bw^{\prime})$, where $A \in GL(g, \BZ)$ and $ B \in M(g, \BZ), \, \,{ ^t B} = B$.
Then the projective unitary representation of $Sp(\CZ)$ on $\CH(\Si,L)$ is described by the following operators representing the generators:
\begin{align} \label{e:unrep}
U_{L}(\al) \, v{\sst (\Si,L)}_{\bq} \, &=  \, v{\sst (\Si,L)}_{^t \! A^{-1} \bq} \notag \\
U_{L}(\be)  \, v{\sst (\Si,L)}_{\bq} \, &= \,  \te^{\frac{\pi \ti}{k} \,{ ^t \bq} B \bq} \, v{\sst (\Si,L)}_{\bq} \\
U_{L}(\ga) \, v{\sst (\Si,L)}_{\bq} \, &= \,   k^{- \frac{g}{2}} \sum_{\bq_1 \in (\BZ /k \BZ)^g} \, \te^{\frac{2\pi \ti}{k} \,{ ^t \bq} \bq_1} \, v{\sst (\Si,L)}_{\bq_1} \notag
\end{align}
For $g=1$ we have the torus mapping class group $SL(2,\BZ)$ with standard generators the matrices $S=\left( \begin{smallmatrix} 0 &1 \\ -1 &0 \end{smallmatrix} \right)$ and $T=\left( \begin{smallmatrix} 1 &1 \\ 0 &1 \end{smallmatrix} \right)$, subject to the relations $(ST)^3=I,\, S^4 =I$.
The projective unitary representation of $SL(2,\BZ)$ on the Hilbert space $\CH(T^2,L)$ associated to the torus $T^2$ and $L \subset H^1(T^2;\BR)$ is described by the operators:
\begin{align*}
U_{L}(T)_{q q'} \, &= \, \te^{\frac{\pi \ti}{k} \, q^2} \: \de_{q q'} \\
U_{L}(S)_{q q'} \, &= \,  k^{- \frac{1}{2}} \, \te^{\frac{2\pi \ti}{k} \,q q'} 
\end{align*}
It coincides with the representation of the modular group $\BP SL(2,\BZ)$ in two dimensional rational conformal field theory \cite{Ve}.
The projective representation (\ref{e:unrep}) of $Sp(\CZ)$ on $\CH(\Si,L)$ is the same projective representation as the one constructed in \cite{G} on the vector space of theta functions at level $k$ obtained through the quantization of the moduli space of flat $\BT$-connections on $\Si$ in a holomorphic polarization. 

\bsk


\section{The path integral approach} \label{s:pathint}

This section attempts to provide a motivation for the definition given at the end of Sect.\ref{s:qth} of the object $Z_X$ that we associate to a compact oriented 3-manifold $X$.
Following Witten's approach \cite{Wi1} we start by defining $Z_X$, for $X$ a closed oriented 3-manifold, as a partition function, that is, a Feynman type functional integral (path integral) over gauge equivalence classes of $\BT$-connections on $X$ with action the Chern-Simons functional.
For the definition of the path integral and its subsequent evaluation we rely almost entirely on the results in \cite{Sch1,Sch2}.
A similar derivation is presented in \cite{A2} in the context of the non-abelian version of the Chern-Simons theory.
In the abelian Chern-Simons theory the path integral is of Gaussian type and the so-called semi-classical or stationary phase approximation gives the exact result.
For the case of a 3-manifold $X$ with boundary, the path integral is defined over a space of connections which are fixed over the boundary. Hence it is a function of the boundary connections and leads, presumably, to an element in the Hilbert space associated to the boundary $\pX$.

\begin{rmk} \label{r:finite}
Let us begin with a brief description of a finite dimensional model which will serve as a prototype for the problem in infinite dimensions of defining $Z_X$ as a partition function.
The main references used are \cite{Sch1,Sch2}, to which we refer the reader for more details and proofs.

Let $G$ be a compact Lie group acting as a group of isometries of the Riemannian manifold $(E,g_E)$.
We assume that the Lie algebra $\text{Lie} \,G$ of $G$ is endowed with an invariant inner product which defines an invariant Riemannian metric on $G$.
Then $\text{Vol} \,G$ denotes the volume of $G$ with respect to this metric.
The action of $G$ on $E$ generates a homomorphism of $\text{Lie} \,G$ into the Lie algebra of vector fields on $E$.
This defines for each point $x \in E$ a linear map $\tau_x : \text{Lie} \,G \ra T_x E$.
Let $H_x$ denote the isotropy subgroup of $G$ at the point $x \in E$ and $\text{Vol} \,H_x$ the volume of $H_x$ with respect to the metric induced by the invariant metric in $G$.
We assume that the isotropy subgroups at all points $x \in E$ are conjugate to a fixed subgroup $H \subset G$.
Then  $\text{Vol} \,H_x =\text{Vol} \,H$, for all $x\in E$.
The $G$-invariant Riemannian metric $g_E$ on $E$ induces a Riemannian metric $g_{E/G}$ on the space $E/G$ of orbits of $G$ in $E$.
It is defined by setting for any $v, w \in T_{[x]} (E/G) \, : \,g_{E/G}(v,w) = g_E (\hat{v},\hat{w})$, where $\hat{v}, \hat{w} \in T_x E$ project onto $v, w$ and belong to the orthogonal complement of $\tau_x(\text{Lie} \,G)$ in $T_x E$.
We let $\mu_E$ and $\mu_{E/G}$ be the measures on $E$ and $E/G$ determined by these metrics.
Then, if $h$ is a $G$-invariant function on $E$, we have
\begin{equation} \label{e:intquot}
\int\limits_{E} h(x) \, \mu_E \, = \, \frac{\text{Vol} \,G}{\text{Vol} \,H} \, \int\limits_{E/G} h(x) \, \bigl[{\det}^{\prime} (\tau_x^* \tau_x) \bigr]^{\frac{1}{2}} \: \mu_{E/G} \, .
\end{equation}
The notation $\det^{\prime}$ refers to the regularized determinant of the operator, that is, the product of all the nonzero eigenvalues.
The $G$-invariant term $\frac{\text{Vol} \,G}{\text{Vol} \,H} \bigl[\det' (\tau_x^* \tau_x) \bigr]^{\frac{1}{2}}$ is the volume of the orbit through the point $x$. A proof of the above formula is given in \cite{Sch1}.

Now let us consider a $G$-invariant real valued function $f$ on $E$ such that the stationary points of $f$ form a $G$-invariant submanifold $F$ of $E$ and $f(x) = A$ for all $x \in F$.
The Hessian $(\text{Hess} f)_x$ of $f$ at the point $x \in E$ is a linear self-adjoint operator $(\text{Hess} f)_x : T_x E \ra T_x E$ and, since the function $f$ is $G$-invariant, $\tau_x(\text{Lie} \,G) \subset \text{Ker} (\text{Hess} f)_x$.
At every point $x$ belonging to the critical manifold $F$, the symmetric bilinear form on $T_x E$ defined by
\begin{equation} \label{e:hess}
(\text{Hess} f)_x (v,w) \, = \, g_E( (\text{Hess} f)_x v, w) \, = \, g_E( v, (\text{Hess} f)_x w) \, ,
\end{equation}
has the expression $(\text{Hess} f)_x (\frac{\partial}{\partial x^i} , \frac{\partial}{\partial x^j} ) \, = \, \frac{{\partial}^2 f}{\partial x^i \partial x^j} $ in local coordinates around $x$.
We assume that $T_x F = \text{Ker} (\text{Hess} f)_x$ at every point $x \in F$.
We also assume that the space $M= F/G$ of orbits of $G$ in $F$ is a manifold and we let $\mu_M$ denote the measure on $M$ determined by the natural metric on $M$ induced from the Riemannian metric of $F$.
The method of stationary phase (\cite{GS,Sch1}) together with the formula (\ref{e:intquot}) give for the asymptotic evaluation of the integral 
\begin{equation} \label{e:oscint}
\frac{1}{\text{Vol} \,G} \, \int\limits_{E} \, \te^{\pi \ti a f(x)} \, \mu_E \qquad  (a \in \BR_{+}) \,,
\end{equation}
in the limit of large $a$, the following expression \cite{Sch1,Sch2}:
\begin{equation} \label{e:mstp}
\frac{1}{ a^{(\dim E - \dim F)/2} } \frac{\te^{\pi \ti a A} }{\text{Vol} \,H} \; \int\limits_M \, \te^{\frac{\pi \ti}{4} \text{sgn} (\text{Hess} f)_x } \; \frac{ \bigl| {\det}^{\prime} (\tau_x^* \tau_x) \bigr|^{\frac{1}{2}}  } { \bigl| {\det}^{\prime} (\text{Hess} f)_x \bigr|^{\frac{1}{2}} } \: \mu_M \, .
\end{equation}
$\text{sgn} (\text{Hess} f)_x$ denotes the signature of the quadratic form in (\ref{e:hess}).
In the terminology of \cite{Sch1,Sch2}, the integral in (\ref{e:oscint}) is the {\em partition function} of the function(al) $f$.
\end{rmk}

In trying to apply the above result to the problem of defining $Z_X$ as a functional integral one needs to make sense of determinants and signatures of operators analogous to the ones in (\ref{e:mstp}), but this time in an infinite dimensional setting.
This is accomplished through the use of zeta-regularization of determinants \cite{RS1,Sch2} and regularization of signatures via eta-invariants \cite{APS, A2}. We recall the relevant definitions.

\begin{rmk} \label{r:zeta}
{\em Zeta-regularization of determinants}. 
We follow \cite{Sch2}.
A non-negative self-adjoint operator $B$ on a Hilbert space $\CH$ is called regular if $\text{Tr} [\te^{-t B} - \Pi(B)] = \sum_k \al_k(B) t^k$ as $t \ra +0$, where $k$ runs over a finite set of non-negative numbers. $\Pi(B)$ denotes the projection operator onto the kernel of $B$.
The zeta-function $\zeta_B(s)$ of a regular operator $B$ is defined for large $\text{Re} (s)$ by the expression
\begin{equation*}
\zeta_B(s) \, = \, \sum_{\la_j \neq 0} \la_j^{-s} \, = \, \frac{1}{\Ga(s)} \int\limits_0^{\infty} \text{Tr} [ \te^{-t B} - \Pi(B) ] \, t^{s-1} \, dt \, ,
\end{equation*}
where $\la_j$ are the eigenvalues of $B$.
The function $\zeta_B(s)$ admits a meromorphic continuation to $\BC$ and is analytic at $s=0$.
The regularized determinant $\det^{\prime} B$ is defined by the expression
\begin{equation} \label{e:regdet}
{\det}^{\prime} B \, = \, \te^{- \zeta^{\prime}_B(0)} 
\end{equation}
For an operator $B : \CH_0 \ra \CH_1$ between the Hilbert spaces $\CH_0$ and $\CH_1$, with adjoint operator $B^*$ and such that $B^* B$ is regular, the regularized determinant $\det^{\prime} B$ is defined by
\begin{equation} \label{e:detB}
{\det}^{\prime} B \, = \, \te^{-\frac{1}{2} \zeta^{\prime}_{B^* B}(0)} \, = \, \bigl[ {\det}^{\prime} (B^* B) \bigr]^{\frac{1}{2}}
\end{equation}
\end{rmk}

\begin{rmk} \label{r:det}
Let $S$ be a self-adjoint operator an a Hilbert space $\CH_1$.
Let $T : \CH_0 \ra \CH_1$ be an operator between the Hilbert spaces $\CH_0$ and $\CH_1$ with adjoint $T^* :  \CH_1 \ra \CH_0$ and such that $T(\CH_0) \subset \text{Ker} S$.
The operators $S^2$ and $T^* T$ are assumed regular.
By definition $\det^{\prime} S = (\det^{\prime}S^2)^{\frac{1}{2}}$ and $\det^{\prime} T = (\det^{\prime} T^* T)^{\frac{1}{2}}$.
Then we have \cite{Sch2}:
\begin{equation} \label{e:detS}
{\det}^{\prime} S \, = \, \frac{ \bigl[\det^{\prime} (S^2 + T T^*) \bigr]^{\frac{1}{2}} }
{\bigl[\det^{\prime} ( T T^*) \bigr]^{\frac{1}{2}} } \, = \, 
\frac{ \bigl[\det^{\prime} (S^2 + T T^*) \bigr]^{\frac{1}{2}} }{ \bigl[\det^{\prime} (T^* T) \bigr]^{\frac{1}{2}} } 
\end{equation}
\end{rmk}

\begin{rmk} \label{r:eta}
{\em Regularization of signatures}.
Let $B$ be a self-adjoint operator on a Hilbert space and define the function 
\begin{equation}
\eta_B (s) \, = \, \sum_{\la_j \neq 0} (\text{sign} \la_j) | \la_j |^{-s}
\end{equation}
with $\la_j$ the eigenvalues of $B$.
The function $\eta_B(s)$ can be meromorphically continued to $\BC$ and has a removable singularity at $s=0$.
The {\em eta-invariant} of the operator $B$ is defined as $\eta(B) = \eta_B(0) \,$ \cite{APS}.
If $B$ is a finite dimensional matrix, then $\eta_B(0)$ is the signature of $B$, i.e. the number of positive eigenvalues minus the number of negative eigenvalues.
\end{rmk}

\ssk

We return now to the problem of defining $Z_X$ as a functional integral.
We consider first a {\em closed} connected oriented 3-manifold $X$.
For each $p \in \text{Tors} H^2(X;\BZ)$, we choose a $\BT$-bundle $P$ on $X$ with $c_1(P) =p$.
The space $\CA_P$ of connections on $P$ is an affine space with vector space $2 \pi \ti \Om^1(X;\BR)$.
The group of gauge transformations $\CG_P \cong \CG_X$ acts on $\CA_P$ by (\ref{e:ggaction}) and, according to (\ref{p:ginv}), the Chern-Simons functional $S_{X,P}: \CA_P \ra \BR/\BZ$ defined in (\ref{d:CSnob}) is invariant under this action.
We define $Z_X$ by the expression
\begin{equation}
Z_X \, = \, \sum_{p \in \text{Tors} H^2(X;\BZ)} Z_{X,p} \, ,
\end{equation}
where $Z_{X,p}$ is the partition function of the Chern-Simons functional, that is,
\begin{equation} \label{e:partf1}
Z_{X,p} \, = \, \int\limits_{\CA_P / \CG_P}  [\CCD \Th] \,\te^{ \pi \ti k S_{X,P}(\Th)} \, .
\end{equation}
The functional integral on the r.h.s. of the above equation is just a formal expression meant to suggest that, pending the existence of a measure $[\CCD \Th]$, we integrate over the space of gauge equivalence classes of connections on $P$.
We are going to show in the following that this functional integral can be given a precise meaning.

Let us choose a Riemannian metric on $X$.
For each $q= 0,1, \dots , \dim X$, the metric determines an inner product on the space $\Om^q(X;\BR)$ of $q$-forms on $X$:
\begin{equation} \label{e:prodforms}
(\al, \be) \, = \, \int\limits_X \al \wedge \star \be
\end{equation}
In particular, it defines an inner product on the tangent space $T \CA_P \cong 2 \pi \ti \Om^1(X;\BR)$.
This makes $\CA_P$ a Riemannian manifold and the group $\CG_P$ acts on $\CA_P$ by isometries.
We let $\hat{\mu}$ denote the measure on $\CA_P$ determined by this metric and $\mu$ the induced measure on the quotient space $\CA_P /\CG_P$.
For each $\Th \in \CA_P$, let $\tau_{\Th}$ denote the differential of the map from $\CG_P$ to $\CA_P$ which defines by (\ref{e:ggaction}) the $\CG_P$-action.
The linear map $\tau_{\Th}$ sends $\text{Lie} \,\CG_P \cong 2 \pi \ti \Om^0(X;\BR)$ into $T_{\Th} \CA_P \cong 2 \pi \ti \Om^1(X;\BR)$.
Thus $\tau_{\Th}$ is identified with the exterior differential $\tau_{\Th} = d : \Om^0(X;\BR) \ra \Om^1(X;\BR)$.
The inner product in $\text{Lie} \,\CG_P \cong 2 \pi \ti \Om^0(X;\BR)$ defined by (\ref{e:prodforms}) induces an invariant metric on the group $\CG_P$.
The isotropy subgroup of $\CG_P$ at a point $\Th \in \CA_P$ is the group of constant maps from $X$ into $\BT$; hence, it is isomorphic to $\BT$.
We let $\text{Vol} \,\BT$ denote the volume of the isotropy subgroup, computed with respect to the metric induced from that of $\CG_P$.
Its evaluation gives
\begin{equation} \label{e:vol}
\text{Vol} \,\BT \, = \, \bigl[ \text{Vol} \,X \bigr]^{\frac{1}{2}} \, = \, \Bigl[ \int\limits_X \star 1 \Bigr]^{\frac{1}{2}} 
\end{equation}

Thus, besides the fact that we are dealing with infinite dimensional spaces, the setting is the same  as the one in Remark (\ref{r:finite}).
We have the identifications $E= \CA_P, \, G=\CG_P$ and $f= k S_{X,P}$.
According to (\ref{p:statnb}) the stationary points of the Chern-Simons functional $S_{X,P}$ are the flat connections; hence the critical manifold $F$ is the subspace $\CA_P^f = \{ \Th \in \CA_P \mid F_{\Th} = d \Th =0 \}$ and the constant value of the functional $S_{X,P}$ on $\CA_P^f$ is the Chern-Simons invariant of flat connections on $P$.
If $\Th_P$ is a flat connection in $\CA_P$, any $\Th \in \CA_P$ can be written as $\Th = \Th_P + 2 \pi \ti A$, for some $A \in \Om^1(X;\BR)$.
Then, since $\pX = \emptyset$ and $F_{\Th_P} = 0$, we obtain from (\ref{l:diffCS}) that
\begin{align*}
S_{X,P}(\Th) \, &= \, S_{X,P}(\Th_P) + \int\limits_X \langle 2 \pi \ti A \wedge d(2 \pi \ti A) \rangle \; \pmod{1} \\
&= \, S_{X,P}(\Th_P)  - (A , \star d A) \; \pmod{1} \notag
\end{align*}
Thus the Hessian of the functional $f= k S_{X,P}$ is the linear operator $\text{Hess} f = - k \star d : \Om^1(X;\BR) \ra \Om^1(X;\BR)$.
The quotient manifold $\CM_P = \CA_P^f/ \CG_P$ is isomorphic to the torus $H^1(X;\BR)/ H^1(X;\BZ)$.
The inner product (\ref{e:prodforms}) on 1-forms determines a natural inner product on $H^1(X;\BR)$ through the Hodge-deRham isomorphism $H^1(X;\BR) \cong \CH^1(X)$ with the space of harmonic 1-forms.
Since the tangent space $T_{[\Th]} \CM_P \cong H^1(X;\BR)$, at any $[\Th] \in \CM_P$, this inner product defines a measure $\nu$ on $\CM_P$.

In view of these observations and of the results (\ref{e:intquot}) and (\ref{e:mstp}) from the model of Remark (\ref{r:finite}), we {\em formally} write:
\begin{align} \label{e:partf2}
Z_{X,p} \, &= \, \frac{1}{\text{Vol} \,\CG_P} \, \int\limits_{\CA_P} \te^{\pi \ti k S_{X,P}(\Th)} \hat{\mu} \\
&= \, \frac{1}{\text{Vol} \,\BT} \, \int\limits_{\CA_P/ \CG_P} \te^{\pi \ti k S_{X,P}(\Th)}  \bigl[ {\det}^{\prime} (\tau^*_{\Th} \tau_{\Th}) \bigr]^{\frac{1}{2}} \, \mu  \tag{\ref{e:partf2} $'$} \\
&= \, \frac{ \te^{\pi \ti k S_{X,P}(\Th_P)} }{\text{Vol} \,\BT} \, \int\limits_{\CM_P} \te^{\frac{\pi \ti}{4} \text{sgn}(- \star d) } 
\frac{ \bigl[ \det^{\prime} (\tau^*_{\Th} \tau_{\Th}) \bigr]^{\frac{1}{2}} }{ \bigl[ \det^{\prime} (- k \star d) \bigr]^{\frac{1}{2}} } \, \nu \tag{\ref{e:partf2} $''$}
\end{align}
The  equality  (\ref{e:partf2}$''$) is the result of the fact that, since $S_{X,P}$ is a quadratic functional, the stationary phase method which gives for an oscillatory integral of the type (\ref{e:oscint}) the asymptotic evaluation (\ref{e:mstp}), produces in this case the exact result.
The expression (\ref{e:partf2}$'$) will be taken as the formal definition for $Z_{X,p}$ in (\ref{e:partf1}).
Although the r.h.s. of (\ref{e:partf2}) is meaningless, the expression in (\ref{e:partf2}$''$) has rigorous mathematical meaning if the determinants and signatures of the operators therein are regularized according to the definitions in Remark (\ref{r:zeta}) and (\ref{r:eta}). 
The signature of the operator $- \star d$ on $\Om^1(X;\BR)$ is regularized via the eta-invariant, that is we take $\text{sgn} (- \star d) = \eta(- \star d)$.
If $\De_q = d^* d + d d^*$ is the Laplacian on $\Om^q(X;\BR)$ and if we let $\zeta_q = \zeta_{\De_q}$, then for any real number $\la > 0$
\begin{equation} \label{e:scaledet}
{\det}^{\prime} (\la \De_q) \, = \, \la^{\zeta_q(0)} \, {\det}^{\prime} (\De_q)
\end{equation}
and \cite{Mu}
\begin{equation} \label{e:zetazero}
\zeta_q(0) \, = \, - \dim \text{Ker} \De_q \, = \, -\dim H^q(X;\BR)
\end{equation}
Now, let us make use of the results stated in remark (\ref{r:det}) to evaluate the regularized determinant $\det'(- k \star d)$.
Thus, we let $S \, = \, - k \star d$ acting on $\Om^1(X;\BR)$ and $T=k d : \Om^0(X;\BR) \ra \Om^1(X;\BR)$.
Then $S^2 + T T^* = k^2 (d^* d + d d^*) = k^2 \De_1$ and $T^* T = k^2 d^* d = k^2 \De_0$.
Using (\ref{e:detS}) and (\ref{e:scaledet})-(\ref{e:zetazero}) we get
\begin{equation} \label{e:detd}
{\det}^{\prime} (- k \star d) = \frac{ \bigl[ \det^{\prime} (k^2 \De_1) \bigr]^{\frac{1}{2}} }{ \bigl[ \det^{\prime} (k^2 \De_0) \bigr]^{\frac{1}{2}}  } = \frac{k^{- \dim H^1(X;\BR)}}{k^{- \dim H^0(X;\BR)}} \, 
\frac{ \bigl[ \det^{\prime} \De_1 \bigr]^{\frac{1}{2}} }{ \bigl[ \det^{\prime} \De_0 \bigr]^{\frac{1}{2}} }
\end{equation}
We also note that $\tau_{\Th}^* \tau_{\Th} = d^* d = \De_0$. Inserting the above results into the expression (\ref{e:partf2}$''$) we obtain
\begin{equation} \label{e:partf3}
Z_{X,p} \, = \, k^{m_X} \, \te^{\pi \ti k S_{X,P}(\Th_P)}  \te^{\frac{\pi \ti}{4} \eta(- \star d) } \, \int\limits_{\CM_P}  \frac{1}{\bigl[ \text{Vol} X \bigr]^{\frac{1}{2}} }
\frac{ \bigl[ \det^{\prime} \De_0 \bigr]^{\frac{3}{4}} }{ \bigl[ \det^{\prime} \De_1 \bigr]^{\frac{1}{4}} } \: \nu 
\end{equation}
where $m_X = \frac{1}{2} \big( \dim H^1(X;\BR) - \dim H^0(X;\BR) \big)$.
The expression under the integral sign in the above formula is related to the Ray-Singer analytic torsion of $X$. 

The analytic torsion of a closed Riemannian manifold was introduced in \cite{RS1,RS2} as a norm on the determinant line of the deRham cohomology of the manifold.
Thus, for the closed 3-manifold $X$ endowed with a Riemannian metric, the Ray-Singer analytic torsion is a density
\begin{align}
T_X^a \in  | \Det H^{\bullet}(X;\BR)^* |  &=
| \Det H^0(X;\BR) |  \otimes | \Det H^1(X;\BR)^* |  \\ &\otimes | \Det H^2(X;\BR) | \otimes | \Det H^3(X;\BR)^* | \notag 
\end{align}
It is defined by the expression \cite{RS2}:
\begin{equation} \label{e:antors}
T_X^a  = \de_{| \Det H^{\bullet}(X;\BR) | } \cdot \exp \Bigl[ \frac{1}{2} \sum\limits_{q=0}^{\dim X} (-1)^q q \zeta^{\prime}_q(0) \Bigr]
\end{equation}
The inner product on the space $\CH^q(X)$ of harmonic $q$-forms defined by the product (\ref{e:prodforms}) on forms determines an inner product on $H^q(X;\BR) \cong \CH^q(X)$ and thus, a density $\de_{| \Det H^{\bullet}(X;\BR) |  }$ on $| \Det H^{\bullet}(X;\BR) |$.
If $b_q =\dim \CH^q(X)$ and $\nu^{(q)}_1, \dots , \nu^{(q)}_{b_q}$ is any orthonormal basis for $\CH^q(X)$, then
\begin{equation}
\de_{| \Det H^{\bullet}(X;\BR) |  } \, = \, \overset{\dim X}{\underset{q=0}{\otimes}} \, 
| \nu^{(q)} |^{(-1)^q} \, ,
\end{equation}
where $\nu^{(q)} = \nu^{(q)}_1 \wedge \cdots \wedge \nu^{(q)}_{b_q}$.
According to the definition of regularized determinants we have $\det^{\prime} \De_q = \te^{- \zeta_q^{\prime}(0)}$.
Thus the Ray-Singer torsion of the closed connected 3-manifold $X$ can be represented as
\begin{align}
T_X^a \, &= \, \bigl[ ({\det}^{\prime} \De_0)^0 ({\det}^{\prime} \De_1)^1 ({\det}^{\prime} \De_2)^{-2} ({\det}^{\prime} \De_3)^3 \bigr]^{\frac{1}{2}}  \times \\
& \: \times | \nu^{(0)} |  \otimes | \nu^{(1)} |^{-1} \otimes |\nu^{(2)} | \otimes | \nu^{(3)} |^{-1} \notag
\end{align}
The Hodge $\star$-operator determines isomorphisms $\De_q \cong \De_{3-q}$ and by Poincar\'{e} duality $H^q(X;\BR) \cong H^{3-q}(X;\BR)^*$.
Moreover, any orthonormal basis $\nu^{(0)}$ of $\CH^0(X) \cong \BR$ is a constant such that $| \nu^{(0)} |= \bigl[ \text{Vol} X \bigr]^{-\frac{1}{2}}$.
Hence the square root of the analytic torsion $T_X^a$ can be regarded as a density $(T_X^a)^{\frac{1}{2}} \in | \Det H^1(X;\BR)^* |$ which is given by the expression
\begin{equation} \label{e:tors}
(T_X^a)^{\frac{1}{2}} \, = \, \frac{1}{ \bigl[ \text{Vol} X \bigr]^{\frac{1}{2}} } \, \frac{ \bigl[ \det^{\prime} \De_0 \bigr]^{\frac{3}{4}} }{ \bigl[ \det^{\prime} \De_1 \bigr]^{\frac{1}{4}} } \: \nu \, ,
\end{equation}
where $\nu = |\nu^{(1)}|^{-1} = | \nu^{(1)}_1 \wedge \cdots \wedge \nu^{(1)}_{b_1}|^{-1}$ for any orthonormal basis $\nu^{(1)}_1 ,\cdots , \nu^{(1)}_{b_1}$ of $\CH^1(X)$.
In view of the isomorphism $\CM_P \cong H^1(X;\BR)/ H^1(X;\BZ)$, the square root $(T_X^a)^{\frac{1}{2}}$ of the analytic torsion defines an invariant density on the moduli space of flat connections $\CM_P$.
Using (\ref{e:tors}), the expression (\ref{e:partf3}) reads:
\begin{equation} \label{e:partf4}
Z_{X,p} \, = \, k^{m_X} \, \te^{\pi \ti k S_{X,P}(\Th_P)} \, \te^{\frac{\pi \ti}{4} \eta(- \star d) } \, \int\limits_{\CM_P} (T_X^a)^{\frac{1}{2}} \end{equation}

It is proved in \cite{RS2} that the analytic torsion of a closed manifold is independent of the metric. Thus $T_X^a$ is a manifold invariant.
The only metric dependence in (\ref{e:partf4}) comes from the phase factor $\te^{\frac{\pi \ti}{4} \eta(- \star d)}$.
Since, after all, {\em we} choose {\em how} to define $Z_{X,p}$, we are going to modify our original definition (\ref{e:partf2}) by multiplying that expression  with $\te^{-\frac{\pi \ti}{4} \eta(- \star d)}$. This makes the resulting expression metric independent.
We remark that in the non-abelian Chern-Simons theory \cite{Wi1}, Witten compensates the metric dependence of the path integral by adding to the original action a counterterm, the gravitational Chern-Simons action.

Starting with a heuristic functional integral formula, we have thus succeeded to associate to the closed connected 3-manifold $X$ the topological invariant
\begin{equation} \label{e:partf5}
Z_X \, = \, k^{m_X} \sum_{p \in \text{Tors} H^2(X;\BZ)} \te^{\pi \ti k S_{X,P}(\Th_P)}  \, \int\limits_{\CM_P} (T_X^a)^{\frac{1}{2}} \end{equation}
We note that the factor $k^{m_X}$ agrees with the one appearing in \cite{FG}, in the analogous expression of the closed 3-manifold invariant for the non-abelian version of the Chern-Simons theory.
\ssk

We treat now the case of a compact connected oriented 3-manifold $X$ {\em with boundary} $\pX$.
We fix a trivializable $\BT$-bundle $Q$ over $\pX$ and a section $s : \pX \ra Q$.
Then, for each $p \in \text{Tors} H^2(X;\BZ)$, we choose a $\BT$-bundle $P \ra X$ with $c_1(P) =p$ and fix a bundle isomorphism $\phi_P : \partial P \ra Q$.
The section $s : \pX \ra Q$ determines a section $s_P = \phi_P^{-1} \circ s : \pX \ra \partial P$.
If $\eta$ is a connection on $Q$, then we define the space
\begin{equation*}
\CA_P(\eta) \, = \, \{ \Th \in \CA_P \mid \partial \Th = \phi^*_P \eta \}
\end{equation*}
of connections on $P$ whose restriction $\partial \Th = \Th \bigr|_{\pX}$ to the boundary of $X$ is identified with $\eta$ under the given bundle isomorphism over $\pX$.
The space $\CA_P(\eta)$ is an affine space with vector space $2 \pi \ti \Om^1_{1,tan}(X;\BR)$, where $\Om^1_{1,tan}(X;\BR)$ denotes the space of 1-forms $A$ on $X$ whose geometrical restriction to $\pX$ is zero (i.e. the pullback $i^* A =0$ under the inclusion map $i : \pX \ra X$).

Let $\CG_P(e) \cong \CG_X(e)$ denote the subgroup in $\CG_P \cong \CG_X$ of gauge transformations which restrict to the identity map over $\pX$.
The group $\CG_P(e)$ acts on $\CA_P(\eta)$ through (\ref{e:ggaction}) and  the action is free.
Given the section $s_P : \pX \ra \partial P$, the Chern-Simons functional $S_{X,P}(s_P, \cdot) : \CA_P(\eta) \ra \BR/\BZ$ defined by (\ref{d:CSb}) is invariant under the $\CG_P(e)$-action on $\CA_P(\eta)$.
As for the closed manifold case previously discussed, we set
\begin{equation} \label{e:partfb1}
Z_{X,p}(\eta) \, = \, \int\limits_{\CA_P(\eta) / \CG_P(e)} \, [\CCD \Th] \,
\te^{\pi \ti k S_{X,P}(s_P,\Th) } 
\end{equation}
and using similar methods aim to show in what follows that one can make sense of such a functional integral.

For this purpose let us choose a Riemannian metric on $X$.
The inner product defined on the tangent space by the metric determines, at each point of $\pX$, a normal vector to the boundary.
Therefore, a differential form $\al$ on $X$ can be decomposed into a tangential and a normal component $\al = \al_{tan} + \al_{norm}$ and one has $(\star \al)_{norm} = \star (\al_{tan})$, $\, (d \al)_{tan} = d(\al_{tan})$, $\, (d^* \al)_{norm} = d^*(\al_{norm})$.
We introduce the subspaces in $\Om^q(X;\BR)$ of $q$-forms on $X$ satisfying {\em relative} boundary conditions \cite{RS1}: 
\begin{equation*}
\Om^q_{tan}(X;\BR)  = \{ \al \in \Om^q(X;\BR) \mid \al_{tan}=(d^* \al)_{tan}=0 \, \text{ on } \pX \} \, ,
\end{equation*}
and {\em absolute} boundary conditions:
\begin{equation*} 
\Om^q_{norm}(X;\BR) = \{ \be \in \Om^q(X;\BR) \mid
\be_{norm}=(d \be)_{norm}=0 \, \text{ on } \, \pX \} \,.
\end{equation*}
Then we let $\De^{tan}_q$ denote the Laplace operator $d^* d + d d^*$ acting on $\Om^q_{tan}(X;\BR)$ and $\De^{norm}_q$ the Laplace operator on $\Om^q_{norm}(X;\BR)$.
If $i^*$ is the pullback under the inclusion $i: \pX \ra X$, the conditions $i^* \al =0$ and $i^*(\star \be)=0$ are equivalent to the conditions $\al_{tan}=0$ on $\pX$ and $\be_{norm}=0$ on $\pX$, respectively.
The Hodge-deRham theory for manifolds with boundary gives the isomorphisms \cite{RS1,Mu}
\begin{align*}
H^q(X,\pX;\BR) \, & \cong \CH^q_{tan}(X) \\
H^q(X;\BR) \, & \cong \CH^q_{norm}(X) \, , \notag
\end{align*}
where
\begin{align*}
\CH^q_{tan}(X)  &= \{ \al \in \Om^q(X;\BR) \mid \al_{tan}=0 \, \text{ on } \, \pX  , \, d \al = d^* \al =0 \, \text{ on } \, X  \} \\
\CH^q_{norm}(X) &= \{ \be \in \Om^q(X;\BR) \mid \be_{norm}=0 \, \text{ on } \, \pX  , \, d \be = d^* \be =0 \, \text{ on } \, X  \} \notag
\end{align*}
are the spaces of harmonic forms on $X$ with relative and absolute boundary conditions.

Let $\hat{\mu}$ be the measure on $\CA_P(\eta)$ determined by the inner product (\ref{e:prodforms}) defined on each tangent space $T_{\Th} \CA_P(\eta) \cong \Om^1_{1,tan}(X;\BR)$ by the metric on $X$ and let $\mu$ be the induced measure on the quotient space $\CA_P(\eta)/ \CG_P(e)$.
For each $\Th \in \CA_P(\eta)$, we have the linear map $\tau_{\Th} =d : \text{Lie} \,\CG_P(e) \cong \Om^0_{tan}(X;\BR) \ra T_{\Th} \CA_P(\eta) \cong \Om^1_{1,tan}(X;\BR)$ determined by the $\CG_P(e)$-action on $\CA_P(\eta)$.

The stationary points of the Chern-Simons functional $S_{X,P}(s_P, \cdot) : \CA_P(\eta) \ra \BR/\BZ$ form the subspace of flat connections $\CA_P^f(\eta) = \{ \Th \in \CA_P(\eta) \mid F_{\Th} = d \Th =0 \}$.
We assume that the connection $\eta$ over $\pX$ is flat and extends to flat connections over $X$.
Otherwise the spaces $\CA_P^f(\eta)$ are empty.
Let us fix an arbitrary flat connection $\Th_P \in \CA_P^f(\eta)$.
Then any connection $\Th \in \CA_P(\eta)$ can be expressed as $\Th = \Th_P + 2 \pi \ti A$, for some $A \in \Om^1_{1,tan}(X;\BR)$, and it follows from (\ref{l:diffCS}) that
\begin{equation*}
S_{X,P}(s_P,\Th) \, = \, S_{X,P}(s_P, \Th_P) - (A, \star d A) \;\pmod{1}
\end{equation*}
Proceeding in the same way as for the closed manifold case, we mimic the finite dimensional model of Remark (\ref{r:finite}) and {\em formally} define the partition function of the quadratic functional $S_{X,P}(s_P,\cdot)$ as
\begin{align} \label{e:partfb2}
Z_{X,p}(\eta) \, &= \, \frac{1}{\text{Vol} \,\CG_P(e)} \, \int\limits_{\CA_P(\eta)} \te^{\pi \ti k S_{X,P}(s_P,\Th)} \,\hat{\mu} \\
&= \, \int\limits_{\CA_P(\eta)/ \CG_P(e)} \te^{\pi \ti k S_{X,P}(s_P,\Th)}  \bigl[ {\det}^{\prime} (\tau^*_{\Th} \tau_{\Th}) \bigr]^{\frac{1}{2}} \, \mu  \tag{\ref{e:partfb2} $'$} \\
&= \,  \te^{\pi \ti k S_{X,P}(s_P,\Th_P)}  \, \int\limits_{\CM_P(\eta)} \te^{\frac{\pi \ti}{4} \text{sgn}(- \star d) } 
\frac{ \bigl[ \det^{\prime} (\tau^*_{\Th} \tau_{\Th}) \bigr]^{\frac{1}{2}} }{ \bigl[ \det^{\prime} (- k \star d) \bigr]^{\frac{1}{2}} } \: \nu_t \tag{\ref{e:partfb2} $''$}
\end{align}
The deRham theory for manifolds with boundary gives the following identification 
$H^1(X,\pX;\BR) = (\text{Ker } d \bigr|_{\Om^1_{1,tan}(X;\BR)}) / (\text{Im } d \bigr|_{\Om^0_{tan}(X;\BR)})$. 
Thus the quotient space $\CM_P(\eta) = \CA_P^f(\eta) /\CG_P(e)$ is isomorphic to the torus $H^1(X,\pX;\BR)/H^1(X,\pX;\BZ)$.
The measure $\nu_t$ on $\CM_P(\eta)$ appearing in (\ref{e:partfb2}$''$) is the measure induced by the inner product determined on $H^1(X,\pX;\BR)$ through the isomorphism with $\CH^1_{tan}(X)$.
The expression (\ref{e:partfb2}$'$) will be taken as the formal definition for $Z_{X,p}(\eta)$ in (\ref{e:partfb1}).
We are going to show in what follows that (\ref{e:partfb2}$''$) has rigorous mathematical meaning if one takes the regularized determinants and signatures of the operators therein.

To evaluate $\det^{\prime}(- k \star d)$ we are going to use the results in Remark (\ref{r:det}).
Let $S = - k \star d$ acting on the space $\Om^1_{tan}(X;\BR)$ and $\tilde{S} = \left( \begin{smallmatrix} 0 & S \\ S^* & 0 \end{smallmatrix} \right)= \left( \begin{smallmatrix} 0 & - k \star d \\ - k \star d & 0 \end{smallmatrix} \right)$ acting on $\Om^1_{tan}(X;\BR) \oplus \Om^1_{norm}(X;\BR)$.
We also let $\tilde{T}$ be the operator $\tilde{T} = \left( \begin{smallmatrix} k d & 0 \\ 0 & k d \end{smallmatrix} \right)$ on the space $\Om^0_{tan}(X;\BR) \oplus \Om^0_{norm}(X;\BR)$.
Then we find that $\tilde{S}^2 + \tilde{T} \tilde{T}^* = k^2 ( \De_1^{tan} \oplus \De_1^{norm})$ and $\tilde{T}^* \tilde{T} = k^2 ( \De_0^{tan} \oplus \De_0^{norm})$.
According to the definition (\ref{e:detB}), we have $\det^{\prime} S = 
| \det^{\prime} (S^* S) |^{\frac{1}{2}} = | \det^{\prime} \tilde{S} |^{\frac{1}{2}}$.
Using this and the relation (\ref{e:detS}), we can write
\begin{align*}
{\det}^{\prime} (- k \star d) \, &= \, | {\det}^{\prime} \tilde{S} |^{\frac{1}{2}}
\, = \, \frac{ \bigl| {\det}^{\prime} (\tilde{S}^2 + \tilde{T} \tilde{T}^*) \bigr|^{\frac{1}{4}} }{\bigl| \det^{\prime} (\tilde{T}^* \tilde{T}) \bigr|^{\frac{1}{4}} } \\
&= \, \frac{ \bigl| \det^{\prime} (k^2 \De_1^{tan}) \, \det^{\prime} (k^2 \De_1^{norm}) \bigr|^{\frac{1}{4}} }{  \bigl| \det^{\prime} (k^2 \De_0^{tan}) \, \det^{\prime} (k^2 \De_0^{norm}) \bigr|^{\frac{1}{4}} } \notag \\
&= \, \frac{ k^{-\frac{1}{2} [\dim H^1(X;\BR) + \dim H^1(X,\pX;\BR)]} }{ k^{-\frac{1}{2} [\dim H^0(X;\BR) + \dim H^0(X,\pX;\BR)]} } \, 
\frac{ \bigl| \det^{\prime} \De_1^{tan} \, \det^{\prime} \De_1^{norm} \bigr|^{\frac{1}{4}} }{  \bigl| \det^{\prime}  \De_0^{tan} \, \det^{\prime} \De_0^{norm} \bigr|^{\frac{1}{4}} } \notag 
\end{align*}
Inserting the above result together with $\det^{\prime} (\tau_{\Th}^* \tau_{\Th}) = \det^{\prime} \De_0^{tan}$ into (\ref{e:partfb2}$''$), we obtain
\begin{align} \label{e:partfb3}
Z_{X,p}(\eta) \, &= \, k^{m_X} \, \te^{\pi \ti k S_{X,P}(s_P,\Th_P)}  \, \te^{\frac{\pi \ti}{4} \eta(- \star d) } \times \\
&\, \times \int\limits_{\CM_P(\eta)} 
\frac{ \bigl| \det^{\prime} \De_0^{tan} \bigr|^{\frac{5}{8}} \, 
\bigl| \det^{\prime} \De_0^{norm} \bigr|^{\frac{1}{8}} }
{ \bigl| \det^{\prime}  \De_1^{tan} \bigr|^{\frac{1}{8}} \, 
\bigl| \det^{\prime} \De_1^{norm} \bigr|^{\frac{1}{8}} }  \: \nu_t \notag 
\end{align}
with $m_X$ defined as in (\ref{e:mx}). 
We are going to relate the term inside the integral sign in the expression (\ref{e:partfb3}) to the Ray-Singer analytic torsion of $X$.

The same definition (\ref{e:antors}) of the Ray-Singer analytic torsion as a density $T_X^a$ on $| \Det H^{\bullet}(X;\BR) |$ applies when $X$ is a Riemannian manifold with boundary \cite{Mu,V}.
In this case the zeta-function $\zeta_q$ in (\ref{e:antors}) is that of the Laplace operator $\De_q^{norm}$ on $\Om^q_{norm}(X;\BR)$ and the density $\de_{| \Det H^{\bullet}(X;\BR) |}$ on $| \Det H^{\bullet}(X;\BR) |$ corresponds to the natural inner product on $H^{\bullet}(X;\BR) \cong \CH_{norm}^{\bullet}(X)$.
Thus the Ray-Singer analytic torsion of the compact connected 3-manifold $X$ with nonempty boundary $\pX$ is the density
$T_X^a \in | \Det H^0(X;\BR) | \otimes | \Det H^1(X;\BR)^* | \otimes | \Det H^2(X;\BR) |$ given by
\begin{align*}
T_X^a \, &= \, \bigl[ ({\det}^{\prime} \De_0^{norm})^0 ({\det}^{\prime} \De_1^{norm})^1 ({\det}^{\prime} \De_2^{norm})^{-2} ({\det}^{\prime} \De_3^{norm})^3 \bigr]^{\frac{1}{2}}  \times \\
&  \times \, | \nu^{(0)}_{norm} |  \otimes | \nu^{(1)}_{norm} |^{-1} \otimes |\nu^{(2)}_{norm} | \notag
\end{align*}
where $\nu^{(q)}_{norm} = \nu^{(q)}_1 \wedge \cdots \wedge \nu^{(q)}_{b_q}$, with $b_q = \dim \CH^q_{norm}(X)$ and $\nu^{(q)}_1 , \dots , \nu^{(q)}_{b_q}$ an orthonormal basis for $\CH^q_{norm}(X)$.
An orthonormal basis $\nu^{(0)}_{norm}$ of $\CH^0(X) \cong \BR$ is such that $| \nu^{(0)}_{norm} | = \bigl[ \text{Vol} X \bigr]^{-\frac{1}{2}}$.
Poincar\'{e} duality gives the isomorphisms $H^q(X;\BR) \cong H^{3-q}(X,\pX;\BR)^*$. 
Also, the spaces $\Om^q_{norm}(X;\BR)$ and $\Om^{3-q}_{tan}(X;\BR)$ are isomorphic under the Hodge $\star$-operator and $\det^{\prime} \De_q^{norm} = \det^{\prime} \De_{3-q}^{tan}$. 
In consequence, the square root of the analytic torsion can be regarded as  a half-density $(T_X^a)^{\frac{1}{2}}$ in $| \Det H^1(X;\BR)^* |^{\frac{1}{2}} \otimes | \Det H^1(X,\pX;\BR)^* |^{\frac{1}{2}}$ which is given by the expression
\begin{equation} \label{e:btors}
(T_X^a)^{\frac{1}{2}} \, = \, \frac{1}{ \bigl[ \text{Vol} X \bigr]^{\frac{1}{4}} } \, \frac{ \bigl[ \det^{\prime} \De_0^{tan} \bigr]^{\frac{3}{4}} \, \bigl[ \det^{\prime} \De_1^{norm} \bigr]^{\frac{1}{4}}  }{ \bigl[ \det^{\prime} \De_1^{tan} \bigr]^{\frac{1}{2}} } \: (\nu_n)^{\frac{1}{2}} \otimes (\nu_t)^{\frac{1}{2}} 
\end{equation}
Here $\nu_n = | \nu^{(1)}_{norm} |^{-1} = | \nu^{(1)}_1 \wedge \cdots \wedge \nu^{(1)}_{b_1} |^{-1}$, for any orthonormal basis $\nu^{(1)}_1 , \dots , \nu^{(1)}_{b_1}$ of harmonic forms of $\CH^1_{norm}(X)$. Similarly, $\nu_t = | \nu^{(1)}_{tan} |^{-1} = | \tilde{\nu}^{(1)}_1 \wedge \cdots \wedge \tilde{\nu}^{(1)}_{l_1} |^{-1}$, for any orthonormal basis $\tilde{\nu}^{(1)}_1 , \dots , \tilde{\nu}^{(1)}_{l_1}$ of  $\CH^1_{tan}(X)$, where $l_1= \dim \CH^1_{tan}(X)$.
Having in view of the isomorphism $\CM_P(\eta) \cong H^1(X,\pX;\BR) /H^1(X,\pX;\BZ)$, we note that 
\begin{equation} \label{e:hde}
\int\limits_{\CM_P(\eta)} (T_X^a)^{\frac{1}{2}} \, \equiv \,
\Big[ \int\limits_{\CM_P(\eta)}  \nu_t \Big] \:
(T_X^a)^{\frac{1}{2}} \otimes (\nu_t)^{-1}
\end{equation}
is a half-density in $| \Det H^1(X;\BR)^* |^{\frac{1}{2}} \otimes | \Det H^1(X,\pX;\BR) |^{\frac{1}{2}}$.

As in the case of a closed manifold, it was proven in \cite{V} that the analytic torsion $T_X^a$ of a manifold $X$ with boundary is independent of the metric chosen on $X$ (the metric is assumed to be a direct product metric near $\pX$).
On the other hand the partition function $Z_{X,P}(\eta)$ in (\ref{e:partfb3}) is a metric dependent complex number.
In the light of the previous observations on the analytic torsion, we remark that we can get rid of this metric dependence if we multiply $Z_{X,P}(\eta)$ by the half-density $\Upsilon \in | \Det H^1(X;\BR)^*|^{\frac{1}{2}} \otimes | \Det H^1(X,\pX;\BR) |^{\frac{1}{2}}$ defined by
\begin{equation} \label{e:Ups}
\Upsilon  =  2^{- \frac{\chi(\pX)}{2}} \frac{ \te^{-\frac{\pi \ti}{4} \eta(- \star d) }  }{\bigl[ \text{Vol} X \bigr]^{\frac{1}{4}} } \,
 \frac{ \bigl[ \det^{\prime} \De_0^{tan} \bigr]^{\frac{1}{8}} \, \bigl[ \det^{\prime} \De_0^{norm} \bigr]^{-\frac{1}{8}}  }{ \bigl[ \det^{\prime} \De_1^{tan} \bigr]^{\frac{3}{8}} \, \bigl[ \det^{\prime} \De_1^{norm} \bigr]^{-\frac{3}{8}} } \: 
(\nu_n)^{\frac{1}{2}} \otimes (\nu_t)^{-\frac{1}{2}} \, ,
\end{equation}
where $\chi(\pX)$ is the Euler characteristic of $\pX$.
Thus, multiplying (\ref{e:partfb3}) by (\ref{e:Ups}) and having (\ref{e:hde}) in view, we obtain
\begin{equation*}
Z_{X,p}(\eta) \, \Upsilon \, = \, k^{m_X} \,
\te^{\pi \ti k S_{X,P}(s_P,\Th_P)}   \int\limits_{\CM_P(\eta)} 
2^{- \frac{\chi(\pX)}{2}} (T_X^a)^{\frac{1}{2}}
\end{equation*}
which is independent of the metric on $X$.

Hence to a compact connected 3-manifold $X$ with boundary $\pX$ the path integral approach outlined above associates, for every flat connection $\eta$ over $\pX$ extending to flat connections over $X$, the half-density 
\begin{equation*}
Z_{X}(\eta) \, \in \, | \Det H^1(X;\BR)^* |^{\frac{1}{2}} \otimes | \Det H^1(X,\pX;\BR) |^{\frac{1}{2}}
\end{equation*}
given by the expression
\begin{align} \label{e:partfb4}
Z_{X}(\eta) \, &= \, \sum_{p \in \text{Tors} H^2(X;\BZ)} Z_{X,p}(\eta) \,\Upsilon \\ &= \, k^{m_X} \sum_{p \in \text{Tors} H^2(X;\BZ)} \te^{\pi \ti k S_{X,P}(s_P,\Th_P)}  
\int\limits_{\CM_P(\eta)}  2^{- \frac{\chi(\pX)}{2}} (T_X^a)^{\frac{1}{2}} \notag
\end{align}
According to the remarks made in Sect.\ref{s:moduli} regarding the group $H^1(X,\pX;\BT)$, we have $\pi_0(H^1(X,\pX;\BT)) \cong \text{Tors} H^2(X,\pX;\BZ) \cong \text{Tors} H^2(X;\BZ)$ and each component of $H^1(X,\pX;\BT)$ is isomorphic to the torus $H^1(X,\pX;\BR)/ H^1(X,\pX;\BZ)$.
Since
\begin{align*}
\int\limits_{\CM_P(\eta)} (T_X^a)^{\frac{1}{2}} \, &= \, 
\int\limits_{H^1(X,\pX;\BR)/ H^1(X,\pX;\BZ)} (T_X^a)^{\frac{1}{2}} \\
&= \, \frac{1}{[ \# \pi_0(H^1(X,\pX;\BT)) ]} \, \int\limits_{H^1(X,\pX;\BT)} (T_X^a)^{\frac{1}{2}} \, ,\notag
\end{align*}
for all $P$, we can rewrite (\ref{e:partfb4}) as
\begin{equation} \label{e:partfb5}
Z_{X}(\eta) \, = \, \frac{k^{m_X}}{[ \# \text{Tors} H^2(X;\BZ) ]}
\sum_{p \in \text{Tors} H^2(X;\BZ)} \te^{\pi \ti k S_{X,P}(s_P,\Th_P)} 
 \int\limits_{H^1(X,\pX;\BT)} 2^{- \frac{\chi(\pX)}{2}}  (T_X^a)^{\frac{1}{2}}
\end{equation}

Ray and Singer conjectured \cite{RS1} that for a closed manifold $X$ the analytic torsion norm $T_X^a$ and the Reidemeister torsion norm $T_X$ on $| \Det H^{\bullet}(X;\BR) |$ coincide
\begin{equation} \label{e:tt}
T_X^a \, = \, T_X .
\end{equation}
A proof of this fact can be found in \cite{V}. Thus, in view of (\ref{e:tt}), the expression (\ref{e:partf5}) for the invariant $Z_X$ associated to the closed connected oriented 3-manifold $X$ through the path integral approach coincides with the expression (\ref{e:Znb}) from the geometric quantization approach.

For a compact manifold $X$ with nonempty boundary $\pX$, \cite{V} proves that the two norms on $| \Det H^{\bullet}(X;\BR) |$, the analytic torsion norm $T_X^a$ and the Reidemeister torsion norm $T_X$,  are related by
\begin{equation} \label{e:ttb}
T_X^a \, = \, 2^{\frac{\chi(\pX)}{2}} \, T_X .
\end{equation}
Referring to the definition of the Chern-Simons section in Sect.\ref{s:line} and to the discussion in Sect.\ref{s:qth}, we note that the terms $\te^{\pi \ti k S_{X,P}(s_P,\Th_P)}$ in (\ref{e:partfb5}) give rise to sections of the prequantum line bundle $\CL_{\pX}$ over the Lagrangian image $\La_X \subset \CM_{\pX}$ of the moduli space of flat $\BT$-connections $\CM_X$.
In view of (\ref{e:hde}) and of the isomorphism (\ref{e:tangLa}), the half-density 
$\int\limits_{H^1(X,\pX;\BT)} (T_X^a)^{\frac{1}{2}}$ can be interpreted as a section of the bundle of half-densities on $\La_X$.
Hence, with the identification (\ref{e:ttb}), the expression (\ref{e:partfb5}) is seen to correspond to the vector $Z_X$ in (\ref{e:Zb}).
It is interesting that the pure path integral (\ref{e:partfb3}) and the half-density factor (\ref{e:Ups}) introduced to compensate the metric dependence of (\ref{e:partfb3}) combine together to produce a result which agrees with the one obtained in the geometric quantization approach.

\bsk


\nin {\bf Acknowledgements.}
I would like to thank Prof. Dan Freed for introducing me to the Chern-Simons theory, as well as for the helpful suggestions and advice during the development of this project.

\bsk



\end{document}